\documentclass[a4paper, 11pt]{amsart}
\usepackage[left=4cm,right=4cm,top=3cm,bottom=3cm]{geometry}
\usepackage[french, english]{babel}

\usepackage{amsfonts, amsthm, amssymb, amsmath, stmaryrd}
\usepackage[T1]{fontenc}
\usepackage{blindtext}

\usepackage{mathrsfs,array}
\usepackage{eucal,color,times,enumerate,accents}
\usepackage{tikz-cd}
\usepackage{url}
\usepackage{bbm}
\usepackage[utf8]{inputenc}
\usepackage{pgfplots}
\usepackage{lipsum}
\usepackage{graphicx}
\usepackage{rotating}
\usepackage{enumitem}
\usepackage{enumerate}
\usepackage{footmisc}
\newenvironment{forcedcentertikzcd}
 {\begin{lrbox}{\forcedcentertikzcdbox}\begin{tikzcd}}
 {\end{tikzcd}\end{lrbox}\makebox[0pt]{\usebox{\forcedcentertikzcdbox}}}
\newsavebox{\forcedcentertikzcdbox}
\usepackage{textcomp}
\usepackage[bookmarksnumbered,bookmarksopen]{hyperref}
\hypersetup{colorlinks, linkcolor=black, citecolor=black}
\usepackage[colorinlistoftodos]{todonotes}
\usepackage[left=4cm,right=4cm,top=3cm,bottom=3cm]{geometry}
\usepackage{mathtools}
\newtheorem{innercustomthm}{Corollary}
\newenvironment{customthm}[1]
{\renewcommand\theinnercustomthm{#1}\innercustomthm}
  {\endinnercustomthm}
\newcommand{\C}{\mathbb{C}}
\newcommand{\calM}{\mathcal{M}}

\newcommand{\F}{\mathbb{F}}
\newcommand{\Q}{\mathbb{Q}}
\newcommand{\Z}{\mathbb{Z}}
\newcommand{\N}{\mathbb{N}}

\newcommand{\Foi}{\Fisoc^\textrm{\dag}}
\newcommand{\oi}{\Isoc^\textrm{\dag}}
\newcommand{\Zl}{\mathbb{Z}_{\ell}}
\newcommand{\Ql}{\mathbb{Q}_{\ell}}

\newcommand{\Dim}{\mathrm{Dim}}
\newcommand{\topo}{\mathrm{top}}
\newcommand{\Ogus}{\mathrm{Ogus}}
\newcommand{\Leray}{\mathrm{Leray}}
\newcommand{\GL}{\mathrm{GL}}
\newcommand{\SL}{\mathrm{SL}}
\newcommand{\Br}{\mathrm{Br}}
\newcommand{\Lie}{\mathrm{Lie}}
\newcommand{\geo}{\mathrm{geo}}
\newcommand{\red}{\mathrm{red}}
\newcommand{\Fr}{\mathrm{Fr}}
\newcommand{\conver}{\mathrm{conv}}
\newcommand{\cristalline}{\mathrm{crys}}
\newcommand{\Rango}{\mathrm{Rank}}

\DeclareMathAlphabet{\mathpzc}{OT1}{pzc}{m}{it}

\newcommand{\Ker}{\mathrm{Ker}}
\newcommand{\Image}{\mathrm{Im}}
\newcommand{\coker}{\mathrm{coker}}

\newcommand{\Spec}{\mathrm{Spec}}

\newcommand{\an}{\mathrm{an}}

\newcommand{\rig}{\mathrm{rig}}

\newcommand{\NS}{\mathrm{NS}}
\newcommand{\Pic}{\mathrm{Pic}}
\newcommand{\speci}{\mathrm{sp}}

\newcommand{\Fisoc}{\mathbf{F\textrm{-}Isoc}}

\newcommand{\calF}{\mathcal{F}}

\newcommand{\calY}{\mathcal{Y}}

\newcommand{\ar}{\mathrm{ar}}

\newcommand{\Strat}{\mathbf{Strat}}

\newcommand*{\sExt}{\mathcal{E}\kern -.5pt xt}
\newcommand*{\sHom}{\mathcal{H}\kern -.5pt om}

\newcommand{\Isoc}{\mathbf{Isoc}}

\newcommand{\Mod}{\mathrm{Mod}}
\begin{document}
\bibliographystyle{alpha}
\newtheorem{theorem}[equation]{Theorem}
\newtheorem{hope}[equation]{Hope}
\newtheorem{proposition}[equation]{Proposition}
\newtheorem{lemma}[equation]{Lemma}
\newtheorem{claim}[equation]{Claim}
\newtheorem{corollary}[equation]{Corollary}
\newtheorem{fact}[equation]{Fact}
\theoremstyle{definition}
\newtheorem{definition}[equation]{Definition}
\newtheorem{question}[equation]{Question}
\newtheorem{conjecture}[equation]{Conjecture}
\newtheorem{answer}[equation]{Answer}
\newtheorem{remark}[equation]{Remark}
\newtheorem{example}[equation]{Example}
\newtheorem{warning}[equation]{Warning}
\newtheorem{notation}[equation]{Notation}
\newtheorem{construction}[equation]{Construction}
\newtheorem*{corollarynoname}{Corollary}

\numberwithin{equation}{subsection}
\title[Specialization of Néron-Severi groups in positive characteristic]{Specialization of Néron-Severi groups in positive characteristic}
\author{Emiliano Ambrosi}
\begin{abstract}
Let $k$ be an infinite finitely generated field of characteristic $p>0$. Fix a separated scheme $X$ smooth, geometrically connected, and of finite type over $k$ and a smooth proper morphism $f:Y\rightarrow X$. The main result of this paper is that there are ``lots of" closed points $x\in X$ such that the fibre of $f$ at $x$ has the same geometric Picard rank as the generic fibre. If $X$ is a curve we show, under a minimal technical assumption, that this is true for all but finitely many $k$-rational points. In characteristic zero, these results have been proved by Andr\'e (existence) and Cadoret-Tamagawa (finiteness) using Hodge theoretic methods. To extend the argument in positive characteristic we use the variational Tate conjecture in crystalline cohomology, the comparison between various $p$-adic cohomology theories and independence techniques. The result has applications to the Tate conjecture for divisors, uniform boundedness of Brauer groups, proper families of projective varieties and to the study of families of hyperplane sections of smooth projective varieties. 
\end{abstract}
\maketitle
\tableofcontents

\section{Introduction}

\subsection{Conventions}
For a field $k$, a $k$-variety is a reduced scheme separated and of finite type over $k$. For a $k$-variety $X$, write $|X|$ for the set of closed points. If $x\in X$, write $k(x)$ for its residue field and $\overline x$ for a geometric point over $x$. If $Y\rightarrow X$ is a morphism and $x\in X$ write $i_{x}:Y_x\rightarrow Y$ for the natural inclusion of the fibre $Y_x$ at $x$ in $Y$. We use $\twoheadrightarrow$ and $\hookrightarrow$ to denote surjective and injective maps respectively. If $\F_q$ is a finite field, write $\F$ for its algebraic closure. If $\mathcal C$ is an abelian category write $\mathcal C\otimes \Q$ for its isogeny category and $\otimes \Q:\mathcal C\rightarrow \mathcal C\otimes \Q$ for the canonical functor.
\subsection{Summary}\label{sectionsummeryNS}
Let $k$ be a finitely generated field of characteristic $p>0$, $\ell\neq p$ a prime, $X$ a smooth geometrically connected $k$-variety, and $f:Y\rightarrow X$ a smooth proper morphism. In first approximation, the main result of this paper is a version of the variational Tate conjecture for divisors in the generic case: for $x\in |X|$, if $H^2(Y_{\overline x},\Ql(1))$ has no more Galois invariants than the generic fibre, then $Y_{\overline x}$ has no more divisors than the generic fibre. When $k$ is a field of characteristic zero, this has been proved by Andr\'e as a consequence of Lefschetz (1,1)-theorem and the Hodge theory in \cite{HodgeII}; see Section \ref{char0NS} for more details.

The starting point of our proof is to replace Hodge theory with crystalline cohomology, since a variational form of the Tate conjecture (Fact \ref{thmvariationalTate}) is known in this setting. The main difficulty to overcome is to transfer the information about the Galois invariants of the $\ell$-adic lisse sheaf $R^2f_*\Ql(1)$ to the crystalline local system ($F$-isocrystal) $R^2f_{\cristalline,*}\mathcal O_{Y/K}(1)$. This is the main new contribution of this paper (Theorem \ref{ltocrys}). More precisely, since the $F$-isocrystal $R^2f_{\cristalline,*}\mathcal O_{Y/K}(1)$ has a behaviour which is quite different from $R^2f_*\Ql(1)$ (for example, in general its cohomology is not finite dimensional), this comparison cannot be done directly. The idea is then to show (Theorem \ref{berthelotconjecture}) that $R^2f_{\cristalline,*}\mathcal O_{Y/K}(1)$ is coming from a smaller and better behaved category of $p$-adic local systems: the category of overconvergent $F$-isocrystals. As it has been understood that overconvergent $F$-isocrystals share many properties with lisse sheaves (\cite{Crew}, \cite{finiteness}, \cite{AbeCaro}), the idea is to compare first $R^2f_{\cristalline,*}\mathcal O_{Y/K}(1)$ with its overconvergent incarnation $R^2f_{*}\mathcal O^{\dagger}_{Y/K}(1)$ via various $p$-adic comparison theorems and then $R^2f_{*}\mathcal O^{\dagger}_{Y/K}(1)$ with $R^2f_*\Ql(1)$ via the theory of weights (\cite{Weil2}, \cite{KatzMessing}).

However, the theory of weights allows us to transfer only information readable on characteristic polynomials of the Frobenii, that is to compare $R^2f_{*}\mathcal O^{\dagger}_{Y/K}(1)$ and $R^2f_*\Ql(1)$ only up to semi-simplification. The way to grasp the missing information is to use Tannakian techniques: instead of considering only $R^2f_*\Ql(1)$, we consider all the possible tensor constructions and sub quotients arising from them, obtaining an algebraic group $G_{\ell}$.  Since $G_{\ell}$ identifies with the Zariski closure of the image of $\pi_1(X,\overline x)$ acting on $(R^2f_*\Ql(1))_{\overline x}\simeq H^2(Y_{\overline x},\Ql(1))$, instead of asking that $H^2(Y_{\overline x},\Ql(1))$ has no more Galois invariants than the generic fibre, we ask that the Zariski closure $G_{\ell,x}$ of the image of $\pi_1(x,\overline x)$ acting on $H^2(Y_{\overline x},\Ql(1))$ identifies with $G_{\ell}$. Then, the theory of weights, combined now with  some algebraic group theory, allows us to compare the $\ell$-adic and the $p$-adic worlds.

Behind this is the idea that, while $R^2f_{*}\mathcal O^{\dagger}_{Y/K}(1)$ and $R^2f_*\Ql(1)$ should be different incarnations of the same motives, each of them contains some specific feature: $R^2f_*\Ql(1)$ can be studied via $\ell$-adic Lie groups theory, while $R^2f_{*}\mathcal O^{\dagger}_{Y/K}(1)$ is an overconvergent incarnation of $R^2f_{\cristalline,*}\mathcal O_{Y/K}(1)$, which, in turn, contains information on the deformations of cycles.
\subsection{Galois-generic points}\label{introduction}
Let $k$ be a field of characteristic $p>0$ with algebraic closure $\overline k$, $X$ a smooth geometrically connected $k$-variety with generic point $\eta$ and $f:Y\rightarrow X$ a smooth proper morphism of $k$-varieties. For $x\in X$, fix an étale path from $\overline x$ to $\overline \eta$. For every $\ell\neq p$, by smooth proper base change $R^2f_*\Ql(1)$ is a lisse sheaf on $X$ and the choice of the étale path gives equivariant isomorphisms
\begin{center}
\begin{tikzcd}[
 row sep=small,
  ar symbol/.style = {draw=none,"\textstyle#1" description,sloped},
  isomorphic/.style = {ar symbol={\cong}},
  ]
    H^2(Y_{\overline \eta},\Ql(1)) \arrow[loop below]\arrow[phantom]{r}{\simeq}& R^2f_*\Ql(1)_{\overline \eta} \arrow[phantom]{r}{\simeq} & R^2f_*\Ql(1)_{\overline x}\arrow[phantom]{r}{\simeq} &  H^2(Y_{\overline x},\Ql(1)) \arrow[loop below]\\
    \pi_1(X,\overline \eta)\arrow[phantom]{r}{\simeq} &  \pi_1(X,\overline x)& & \pi_1(x,\overline x)\arrow{ll}.
\end{tikzcd}
\end{center}
\begin{definition}
A point $x\in X$ is $\ell$-Galois-generic (resp. strictly $\ell$-Galois-generic) for $f:Y\rightarrow X$ if the image of $\pi_1(x,\overline x)\rightarrow \pi_1(X,\overline \eta)\rightarrow \GL(H^2(Y_{\overline \eta},\Ql(1)))$ is open (resp. coincides with) in the image of $\pi_1(X,\overline \eta)\rightarrow \GL(H^2(Y_{\overline \eta},\Ql(1)))$.
\end{definition}
By \cite[Theorem 1.1]{uniformadelic}\footnote{Recall that finite fields are in particular $\ell$-non-Lie semisimple for every $\ell$ different from the characteristic, so that \cite[Theorem 1.1]{uniformadelic} applies in our setting.}, $x$ is $\ell$-Galois-generic for one $\ell\neq p$ if and only if $x$ is $\ell$-Galois-generic for every $\ell\neq p$. So one simply says that $x$ is Galois-generic for $f$. This is not true for strictly Galois-generic points, and one says that $x$ is strictly Galois-generic if there exists an $\ell\neq p$ such that $x$ is strictly $\ell$-Galois-generic.
\numberwithin{equation}{subsection}
 \subsection{N\'eron-Severi generic points}\label{nerongeneric}
 \numberwithin{equation}{subsubsection}
 \subsubsection{Tate conjecture for divisors}
The geometric N\'eron-Severi group $\NS(Z_{\overline k})$ of a smooth proper $k$-variety $Z$ is a finitely generated abelian group such that $\NS(Z_{\overline k})\otimes \Q$ identifies with the image of the cycle class map for $\ell$-adic cohomology
$$c_{Z_{\overline k}}:\Pic(Z_{\overline k})\otimes \Q\rightarrow H^2(Z_{\overline k},\Ql(1)).$$
Since $\NS(Z_{\overline k})$ is a finitely generated abelian group, $\pi_1(k)$ acts on it through a finite quotient and hence $\NS(Z_{\overline k})\subseteq H^2(Z_{\overline k},\Ql(1))$ is fixed under the action of the connected component $G^0_{\ell}$ of the Zariski closure $G_{\ell}$ of the image $\Pi_{\ell}$ of $\pi_1(k)$ acting on $H^2(Z_{\overline k},\Ql(1))$. Recall that the $\ell$-adic Tate conjecture for divisors (\cite{tateconj}) predicts the following:
\begin{conjecture}[$T(Z,\ell)$]\label{tateconjectureconj}
Let $k$ be a finitely generated field and $Z$ a smooth proper $k$-variety. Then the map
$c_{Z_{\overline k}}:\Pic(Z_{\overline k})\otimes \Ql\rightarrow H^2(Z_{\overline k},\Ql(1))^{G^0_{\ell}}$
is surjective.
\end{conjecture}
\subsubsection{Specialization morphisms}
Retain the notation and the assumptions of Section \ref{introduction}.
For every $x\in X$, there is an injective specialization homomorphism (see e.g. \cite[Proposition 3.6.]{poonen})
$$\speci_{\eta,x}:\NS(Y_{\overline \eta})\otimes \Q\rightarrow \NS(Y_{\overline x})\otimes \Q$$
compatible with the cycle class map, in the sense that the following diagram commutes:
\begin{center}
\begin{tikzcd}
\Pic(Y_{\eta})\otimes \Q\arrow{d}\arrow[bend right=70, swap]{dd}{c_{Y_{\eta}}} & \Pic(Y)\otimes \Q\arrow[swap]{l}{i_{\eta}^*}\arrow{r}{i_{x}^*} &  \Pic(Y_{x})\otimes \Q \arrow{d}\arrow[bend left=70]{dd}{c_{Y_x}}\\
\NS(Y_{\overline{\eta}})\otimes \Q\arrow[hook]{rr}{\speci_{\eta,x}}\arrow[hook]{d} && \NS(Y_{\overline x})\otimes \Q\arrow[hook]{d}\\
 H^2(Y_{\overline{\eta}},\Ql(1))\arrow[phantom]{rr}{\simeq} && H^2(Y_{\overline{x}},\Ql(1))
\end{tikzcd}
\end{center}
Since the N\'eron-Severi group is invariant under extensions of algebraically closed fields (see e.g. \cite[Proposition 3.1]{poonen}), the map $\speci_{\eta,x}$ is well defined, independently of the choice of the geometric points $\overline \eta$ over $\eta$ and $\overline x $ over $x$. 

The abelian group $\NS(Y_{\overline{\eta}})\otimes \Q$ is a $\pi_1(X,\overline \eta)$-module and hence the group $\pi_1(x,\overline x)$ acts on $\NS(Y_{\overline{\eta}})\otimes \Q$ by restriction through the morphism $\pi_1(x,\overline x)\rightarrow \pi_1(X,\overline x)\simeq \pi_1(X,\overline \eta)$. Since the map $\speci_{\eta,x}$ is $\pi_1(x,\overline x)$-equivariant with respect to the natural action of $\pi_1(x,\overline x)$ on $\NS(Y_{\overline x})\otimes \Q$, one constructs an injective specialization map
$$\speci^{\ar}_{\eta,x}:\NS(Y_{\eta})\otimes \Q\subseteq (\NS(Y_{\overline \eta})\otimes \Q)^{\pi_1(x,\overline x)}\xrightarrow{\speci_{\eta,x}} \NS(Y_{x})\otimes \Q,$$
where for a smooth proper $k$-variety $Z$ one writes $\NS(Z)\otimes \Q:=(\NS(Z_{\overline k})\otimes \Q)^{\pi_1(k)}$.
\begin{definition}
One says that $x$ is $\NS$-generic (resp. arithmetically $\NS$-generic) for $f:Y\rightarrow X$ if $\speci_{\eta,x}$ (resp. $\speci^{\ar}_{\eta,x}$) is an isomorphism.
\end{definition}
Conjecture \ref{tateconjectureconj} predicts that every (strictly) Galois-generic point is (arithmetically) $\NS$-generic.
Our main result is that this holds (without assuming Conjecture \ref{tateconjectureconj}), at least when $f:Y\rightarrow X$ is projective.
\begin{theorem}\label{main}
Let $k$ be a finitely generated field and $f:Y\rightarrow X$ a smooth projective morphism, where $X$ is a smooth and geometrically connected $k$-variety. If $x\in X$ is Galois-generic (resp. strictly Galois-generic) for $f:Y\rightarrow X$ then it is $\NS$-generic (resp. arithmetically $\NS$-generic) for $f:Y\rightarrow X$. If $f:Y\rightarrow X$ is smooth and proper, the same is true for all $x$ in a dense open subset of $X$.
\end{theorem}
\numberwithin{equation}{subsection} 
\subsection{Proof in characteristic zero}\label{char0NS}
When $k$ is a field of characteristic zero Theorem \ref{main} is due to Andr\'e (\cite{AndreIHES}; see also \cite[Corollary 5.4]{MC} and \cite[Proposition 3.2.1]{brauer}) and it holds for $f:Y\rightarrow X$ proper. Since it is the starting point for our proof we briefly recall the argument when $k\subseteq \mathbb C$ and $x$ is a closed point. 
Fix a smooth compactification $Y\subseteq \overline Y$ of $Y$. The commutative diagram of $k$-varieties
\begin{center}
\begin{tikzcd}
Y_{x}\arrow{r}\arrow{d}\arrow[phantom]{dr}{\Box} & Y\arrow[hook]{r}\arrow{d}& \overline Y\\
k(x)\arrow{r}{x} & X
\end{tikzcd}
\end{center}
induces a commutative diagram:
\begin{center}
\begin{tikzcd}
H^0(X_{\mathbb C},R^2f_*\Q(1))\arrow{d}{\simeq} & H^{2}(Y_{\mathbb C},\Q(1))\arrow{l}{\Leray} & \NS(Y_{\overline \eta})\otimes \Q\arrow[bend left=100]{dd}{\speci_{\eta,x}}\\
H^{2}(Y_{\overline x},\Q(1))^{\pi^{\topo}_1(X_{\mathbb C})}\arrow[hook]{d} & H^2(\overline Y_{\C}, \Q(1))\arrow{l}\arrow{u} & \NS(\overline Y_{\C})\otimes \Q\arrow{d}\arrow{u}\arrow[hook]{l}\\
H^{2}(Y_{\overline x},\Q(1))& & \NS(Y_{\overline x})\otimes \Q\arrow[hook]{ll}
\end{tikzcd}
\end{center}
where $\Leray$ is the edge map in the Leray spectral sequence attached to $f:Y\rightarrow X$. 
Take any $z_x\in \NS(Y_{\overline x})\otimes \Q$. Since $z_x$ is fixed by an open subgroup of $\pi_1(x)$ and $x$ is Galois-generic, up to replacing $X$ with a finite étale cover one can assume that $z_x$ is fixed by $\pi_1(X_{\C})$. By the comparison between the étale and the singular sites, $z_x$ is fixed by $\pi^{\topo}_1(X_{\C})$. 
By Deligne's fixed part theorem (\cite[Theoreme 4.1.1]{HodgeII}) the map
$$H^2(\overline Y_{\C}, \Q(1))\rightarrow H^{2}(Y_{\overline x},\Q(1))^{\pi^{\topo}_1(X_{\mathbb C})}$$
is surjective. By semisemplicity, the map $H^2(\overline Y_{\C}, \Q(1))\rightarrow H^{2}(Y_{\overline x},\Q(1))^{\pi^{\topo}_1(X_{\mathbb C})}$ splits in the category of polarized $\Q$-Hodge structures, so that $z_x$ is the image of an element $z$ in $H^{0,0}(\overline Y_{\C},\Q(1))$. By the Lefschetz (1,1) theorem, $z$ lies in $\NS(\overline Y_{\mathbb C})\otimes \Q$. One concludes the proof observing that, by construction, the restriction $z_{\eta}$ of $z$ to $\NS(Y_{\overline \eta})\otimes \Q$ is an element such that $\speci_{\eta,x}(z_{\eta})=z_x$.
\numberwithin{equation}{subsection} 
\subsection{Strategy in positive characteristic}
In characteristic zero the main ingredients are the combination of Deligne's fixed part theorem and the Lefschetz-(1,1) theorem (what is called the variational Hodge conjecture for divisors; see e.g. \cite[Conjecture 9.6, Remark 9.7]{poonen}) and the comparison between the étale and the singular sites. To try and make the argument of Section \ref{char0NS} works in positive characteristic the idea is to replace Betti cohomology with crystalline cohomology. The main reason for this is that the variational Tate conjecture for projective morphisms (Fact \ref{thmvariationalTate}), that we now recall, is known in this setting.
\numberwithin{equation}{subsubsection}
\subsubsection{Crystalline variational Tate conjecture}
Let $\F_q$ be the finite field with $q=p^s$ elements, $\mathcal X$ a connected smooth $\F_q$-variety and $\mathfrak f:\mathcal Y\rightarrow \mathcal X$ a smooth proper morphism of $\F_q$-varieties (in our application $\mathfrak f:\mathcal Y\rightarrow \mathcal X$ is a model for $f:Y\rightarrow X$). Write respectively  $\Mod(\mathcal X|W)$, $\Mod(\mathcal Y|W)$ for the categories of $\mathcal O_{\mathcal X|W}$, $\mathcal O_{\mathcal Y|W}$ modules in the crystalline site of $\mathcal X$, $\mathcal Y$ over $W:=W(\F_q)$ (\cite[Appendix A.1]{Varmor}). Then there is a higher direct image functor  
$$R^{i}\mathfrak f_{\cristalline,*}: \Mod(\mathcal Y|W)\rightarrow \Mod(\mathcal X|W)$$
and, for every $\mathfrak t\in \mathcal X(\F_q)$, a commutative diagram 
\begin{center}
\begin{tikzcd}
H^2_{\cristalline}(\mathcal Y)\arrow{d}{\Leray}\arrow{rd}{i^*_{\mathfrak t}}& \Pic(\mathcal Y)\otimes \Q \arrow[swap]{l}{c_{\mathcal Y}}\arrow{rd}{i^*_{\mathfrak t}}\\
H^0(\mathcal X,R^2\mathfrak f_{\cristalline, *}\mathcal O_{\mathcal Y/W})\otimes \Q\arrow{r} & H^{2}_{\cristalline}(\mathcal Y_{\mathfrak t}) & \Pic(\mathcal Y_{\mathfrak t})\otimes \Q\arrow{l}{c_{\mathcal Y_{\mathfrak t}}}
\end{tikzcd}
\end{center}
where $H^2_{\cristalline}(\mathcal Y)$ and $H^{2}_{\cristalline}(\mathcal Y_{\mathfrak t})$ are the (rational) crystalline cohomology of $\mathcal Y$ and $\mathcal Y_{\mathfrak t}$ respectively, Leray is the edge map in the Leray spectral sequence attached to $\mathfrak f:\mathcal Y\rightarrow \mathcal X$ and $c_{\mathcal Y},c_{\mathcal Y_{\mathfrak t}}$ are the crystalline cycle class maps. Write $F$ for the $s$-power of the absolute Frobenius of $\mathcal X$ and recall that the images of $c_{\mathcal Y}$ and $c_{\mathcal Y_{\mathfrak t}}$ lie in $H^2_{\cristalline}(\mathcal Y)^{F=q}$ and $H^2_{\cristalline}(\mathcal Y_{\mathfrak t})^{F=q}$, respectively. Then we have the variational Tate conjecture in crystalline cohomology:
\begin{fact}\label{thmvariationalTate}\cite[Theorem 1.4]{Varmor}
If $\mathfrak f:\mathcal Y\rightarrow \mathcal X$ is projective, for every $z_{\mathfrak t}\in \Pic(\mathcal Y_{\mathfrak t})\otimes \Q$ the following are equivalent:
\begin{enumerate}
\item There exists $z\in \Pic(\mathcal Y)\otimes \Q$ such that $c_{\mathcal Y_{\mathfrak t}}(z_{\mathfrak t})=i^*_{\mathfrak t}(c_{\mathcal Y}(z))$;
\item $c_{\mathcal Y_{\mathfrak t}}(z_{\mathfrak t})$ lies in $H^0(\mathcal X,R^{2}\mathfrak f_{\cristalline, *}\mathcal O_{\mathcal Y/W})\otimes \Q$;
\item $c_{\mathcal Y_{\mathfrak t}}(z_{\mathfrak t})$ lies in $H^0(\mathcal X,R^{2}\mathfrak f_{\cristalline, *}\mathcal O_{\mathcal Y/W})^{F=q}\otimes \Q$.
\end{enumerate}
\end{fact}
However, to apply Fact \ref{thmvariationalTate} in our setting, there are two difficulties to overcome:
\begin{enumerate}
\item Crystalline cohomology works well only over a perfect field, while our base field $k$ is not perfect;
\item There is no direct way to compare the étale and the crystalline sites, so that one has to find a different way to transfer the Galois-generic assumption to the crystalline setting. 
\end{enumerate}
\numberwithin{equation}{subsubsection}
\subsubsection{Spreading out}\label{introspreadingNS}
To overcome (1) one uses a spreading out argument, so that our morphism $f:Y\rightarrow X$ will appear as the generic fibre of a smooth proper morphism $\mathfrak f:\mathcal Y\rightarrow \mathcal X$, where $\mathcal X$ is a smooth geometrically connected $\F_q$-variety. The idea is then to lift an element $\varepsilon_x\in \NS(Y_{\overline x})\otimes \Q$ to $\NS(Y_{\overline \eta})\otimes \Q$ by specializing it first to an element $\varepsilon_{\mathfrak t}\in \NS(\mathcal Y_{\overline{\mathfrak t}})\otimes \Q $ of the geometric fibre $\mathcal Y_{\overline{\mathfrak t}}$ of the morphism $\mathfrak f:\mathcal Y\rightarrow \mathcal X$ over a closed point $\mathfrak t$ and then to try and lift $\varepsilon_{\mathfrak t}$ to an element $\varepsilon\in \Pic(\mathcal Y)\otimes \Q$, via the crystalline variational Tate conjecture over $\F_q$.
\numberwithin{equation}{subsubsection}
\subsubsection{From \texorpdfstring{$\ell$}- to p}\label{seccrysttol}
In order to show that $\varepsilon_{\mathfrak t}\in \NS(Y_{\overline{\mathfrak t}})\otimes \Q $ satisfies the assumption of Fact \ref{thmvariationalTate}, one has to transfer the $\ell$-adic information that $x$ is Galois-generic to crystalline cohomology. For this the key ingredient is Theorem \ref{ltocrys} below.
Assume that $\mathcal Z$ is a smooth geometrically connected  $\F_q$-variety admitting an $\F_q$-rational point $\mathfrak t$ and that there is a map $\mathfrak g:\mathcal Z\rightarrow \mathcal X$ (in our application $\mathfrak g:\mathcal Z\rightarrow \mathcal X$ is a model for $x:k(x)\rightarrow X$). The cartesian square
\begin{center}
\begin{tikzcd}
\mathcal Y_{\mathcal Z}\arrow{r}\arrow{d}{\mathfrak f_{\mathcal Z}}\arrow[phantom]{rd}{\Box} & \mathcal Y\arrow{d}{\mathfrak f}\\
\mathcal Z\arrow{r}{\mathfrak g} & \mathcal X
\end{tikzcd}
\end{center}
induces representations
$$\pi_1(\mathcal Z,\overline{\mathfrak t})\rightarrow \pi_1(\mathcal X,\overline{\mathfrak t})\rightarrow \GL(H^{i}(\mathcal Y_{\overline{\mathfrak t}},\Ql(j)).$$
\begin{theorem}\label{ltocrys}
Assume that the image of  $\pi_1(\mathcal Z,\overline{\mathfrak t})\rightarrow \pi_1(\mathcal X,\overline{\mathfrak t})\rightarrow \GL(H^{i}(\mathcal Y_{\overline{\mathfrak t}},\Ql(j))$ is open in the image of $\pi_1(\mathcal X,\overline{\mathfrak t})\rightarrow \GL(H^{i}(\mathcal Y_{\overline{\mathfrak t}},\Ql(j))$ and that the Zariski closures of the images of $\pi_1(\mathcal X,\overline{\mathfrak t})$ and $\pi_1(\mathcal X_{\mathbb F},\overline{\mathfrak t})$ acting on $H^{i}(\mathcal Y_{\overline{\mathfrak t}},\Ql(j))$ are connected.
Then the base change map
$$H^{0}(\mathcal X,R^{i}\mathfrak f_{\cristalline,*}\mathcal O_{\mathcal Y/W})^{F=q^j}\otimes \Q\rightarrow H^{0}(\mathcal Z,R^{i}\mathfrak f_{\mathcal Z, crys,*}\mathcal O_{\mathcal Y_{\mathcal Z}/W})^{F=q^j}\otimes \Q$$
is an isomorphism.
\end{theorem}
As mentioned in Section \ref{sectionsummeryNS}, the subtle point in the proof of Theorem \ref{ltocrys} is to compare the category of $F$-isocrystals, where the crystalline variational Tate conjecture holds, with the category of $\ell$-adic lisse sheaves. These categories behaves differently. For example, if $\mathfrak f:\mathcal Y\rightarrow \mathcal X$ is a non-isotrivial family of ordinary elliptic curves, $R^1\mathfrak f_{\cristalline,*} \mathcal O_{\mathcal Y/W}\otimes \Q$ carries a two step filtration, reflecting the decomposition of the $p$-divisible groups of the generic fibre of $\mathfrak f:\mathcal Y\rightarrow \mathcal X$ into étale and connected parts, while $R^1\mathfrak f_{*} \Ql$ is irreducible. This leads to consider the smaller category of overconvergent $F$-isocrystals, whose behaviour is closer to the one of $\ell$-adic lisse sheaves. Then the proof of Theorem \ref{ltocrys} decomposes as follows:
\begin{enumerate}
\item We prove that $R^{i}\mathfrak f_{\cristalline,*}\mathcal O_{\mathcal Y/W}\otimes \Q$ and $R^{i}\mathfrak f_{\mathcal Z, crys,*}\mathcal O_{\mathcal Y_{\mathcal Z}/W}\otimes \Q$ are overconvergent $F$-isocrystals (Theorem \ref{Shiho}, which uses a technical result proved in the Part 3, building on the work of Shiho on relative log convergent cohomology and relative rigid cohomology \cite{ShihoI}, \cite{Shiho});
\item We use that one doesn't lose information passing from overconvergent $F$-isocrystals to $F$-isocrystals (Fact \ref{kedlaya});
\item The assumption implies that the Zariski closures $G_{\ell}$ and $G_{\mathcal Z,\ell}$ of the image of $\pi_1(\mathcal X,\overline{\mathfrak t})$ and $\pi_1(\mathcal Z,\overline{\mathfrak t})$ acting on $H^{i}(\mathcal Y_{\overline{\mathfrak t}},\Ql(j))$ are equal.
\item To show that $(3)$ implies Theorem \ref{ltocrys}, one uses (1),(2), the theory of Frobenius weights and some algebraic group theory.
\end{enumerate}
\begin{remark}
In Theorem \ref{ltocrys}, the assumptions that $\mathcal Z$ has a $\F_q$ rational point and that the Zariski closure of the image of $\pi_1(\mathcal X_{\mathbb F},\overline{\mathfrak t})$ is connected are not necessary, but the proof without these assumptions requires the more elaborated formalism of $\overline \Q_p$-overconvergent isocrystals. In our application to Theorem \ref{main} one can reduce to the case where these assumptions are satisfied, so that we did not include the proof of the general form of Theorem \ref{ltocrys}.
\end{remark}
\begin{remark}
In characteristic 0, the proof sketched in Section \ref{char0NS} shows that the variational Hodge conjecture implies the $\ell$-adic variational Tate conjecture over fields of characteristic zero. 
In positive characteristic, our method does not show that the crystalline variational Tate conjecture implies the $\ell$-adic one. The issue comes from the fact that one does not know how to compare the $\ell$-adic and the crystalline cycle class maps.
\end{remark}
\numberwithin{equation}{subsection} 
\subsection{Applications}
\numberwithin{equation}{subsubsection}
\subsubsection{Existence of $\NS$-generic points}
Let $k$ be a field of transcendence degree $\geq 1$ over $\F_p$ and $X$ a smooth geometrically connected $k$-variety with generic point $\eta$. Before explaining our applications, we recall some definition. We start recalling  (\cite[Definition 8.1]{poonen}) the definition of sparse set.
\begin{definition}
A subset $S$ of $|X|$ is said to be sparse if there exists a dominant and generically finite morphism $g:X'\rightarrow X$ with $X'$ an irreducible $k$-variety such that for each $x\in S$ the fiber $g^{-1}(x)$ is either empty or contains more than one closed point.
\end{definition}
If $k$ is Hilbertian (in particular, by Hilbert irreducibility theorem, if $k$ is finitely generated) and $S\subseteq |X|$ is sparse, then there exists an integer $d\geq 1$, such that $|X|-S$ contains infinitely many points of degree $\leq d$ (see e.g. \cite[Remark 8.8]{poonen} or \cite[Lemma 1.2.2.5.2]{miotesi}). 

Assume now that $k$ is finitely generated over $\F_p$. Write $\F_q$ (resp. $\F$) for the algebraic closure of $\F_p$ in $k$ (resp. $\overline k$) and set $k\F\subseteq \overline k$ for the field generated by $k$ and $\F$. 
Let $f:Y\rightarrow X$ be a smooth proper morphism of $k$-varieties. Consider the normal inclusions 
$$\pi_1(X_{\overline k},\overline \eta)\subseteq \pi_1(X_{k\F},\overline \eta)\subseteq \pi_1(X,\overline \eta),$$
and write
$$\Pi:=\Image(\pi_1(X_{\overline k},\overline \eta)\rightarrow \pi_1(X,\overline \eta)\rightarrow \GL(H^2(Y_{\overline{\eta}},\mathbb Q_{\ell}(1))))$$
$$\widetilde \Pi:=\Image(\pi_1(X_{k\F},\overline \eta)\rightarrow \pi_1(X,\overline \eta)\rightarrow \GL(H^2(Y_{\overline{\eta}}),\mathbb Q_{\ell}(1)))),$$
so that $\Pi\subseteq \widetilde \Pi$ is a normal subgroup.
Following \cite{tamagawa}, we say that $f:Y\rightarrow X$ has the LU-property (Lie-Unrelated\footnote{The terminology comes from the fact that the LU-property is equivalent to ask that the Lie algebras $\Lie(\Pi)$ and $\Lie(\widetilde \Pi/\Pi)$ of $\Pi$ and $\widetilde \Pi/\Pi$ have no nontrivial common quotient.} property) in degree 2, if for all open subgroups $U\subseteq \Pi$, $V\subseteq \widetilde \Pi/\Pi$ all the common quotients of $U$ and $V$ are finite. Then we have the following:
\begin{fact}\label{sisagiaspec}
Assume that $k$ is finitely generated. Then:
\begin{enumerate}
\item  The subset of non-strictly $\ell$-Galois-generic closed points for $f:Y\rightarrow X$ is sparse. In particular there exists an integer  $d\geq 1$ such that there are infinitely many $x\in |X|$ with $[k(x):k]\leq d$ that are strictly $\ell$-Galois-generic for $f:Y\rightarrow X$. 
\item If $X$ is a curve and $f:Y\rightarrow X$ has the LU-property in degree 2, then all but finitely many $x\in X(k)$ are Galois-generic for $f:Y\rightarrow X$.
\end{enumerate}
\end{fact}
\proof
\begin{enumerate}
\item[]
\item This follows form the arguments in \cite[Section 10.6]{Serreweil} (see for example \cite[Section 1.2.2.5]{miotesi}). For the reader convenience, we quickly recall the argument. 

Let $X_{\ell}^{\mathrm{nsGg}}$ be the set of $x\in \vert X\vert$ that are not strictly $\ell$-Galois generic and let $\Phi$ the Frattini subgroup (i.e. the intersection of all maximal proper closed subgroups) of the image $\Pi_{\eta}$ of $$\pi_1(X,\overline \eta)\rightarrow \GL(H^2(Y_{\overline{\eta}},\mathbb Q_{\ell}(1))).$$
Since $\Pi_{\eta}$ is an $\ell$-adic Lie group, $\Phi\subseteq \Pi_{\eta}$ is an open subgroup by \cite[Proposition page 148]{Serreweil}. Hence there exist only finitely many proper open subgroups $H_1,\dots, H_n$ of $\Pi_{\eta}$ containing $\Phi$, and for any proper closed subgroup of $\Pi_{\eta}$ there exists some $H_i$ containing it. For $1\leq i\leq n$, let $g_i:X_i\rightarrow X$ be the finite connected étale cover of $X$ corresponding to the preimage of $H_i$ along 
$$\pi_1(X,\overline \eta)\rightarrow \GL(H^2(Y_{\overline{\eta}},\mathbb Q_{\ell}(1))).$$
Then, by the Galois formalism (see e.g. \cite[Lemma 1.2.2.3.1]{miotesi}), one has 
$$X_{\ell}^{\mathrm{nsGg}}=\bigcup_{1\leq i\leq n}\bigcup_{[k':k]< +\infty}\Image(X_i(k')\rightarrow X(k')).$$
Since for each $1\leq i\leq n$ the set
$$\bigcup_{[k':k]< +\infty}\Image(X_i(k')\rightarrow X(k'))$$ is sparse by definition, the conclusion follows from the fact that a finite union of sparse sets is sparse (see e.g. \cite[Proposition 8.5 (b)]{poonen}).
\item This is proved in \cite{tamagawa}.\endproof
\end{enumerate}
\bigskip
So we get the following corollary.
\begin{corollary}\label{existence}
Assume that $k$ is finitely generated. Then:
\begin{enumerate}
\item The subset of closed non-arithmetically $\NS$-generic (resp. non-$\NS$-generic) points for $f:Y\rightarrow X$ is sparse. In particular there exists an integer $d\geq 1$ such that there are infinitely many $x\in |X|$ with $[k(x):k]\leq d$ that are arithmetically $\NS$-generic (resp. $\NS$-generic) for $f:Y\rightarrow X$.
\item  If $X$ is a curve and $f:Y\rightarrow X$ has the LU-property in degree 2, all but finitely many $x\in X(k)$ are $\NS$-generic for $f:Y\rightarrow X$.
\end{enumerate}
\end{corollary}
\proof  
\begin{enumerate}
\item[]
	\item Since if $S\subseteq |X|$ is a subset and $U\subseteq X$ is a dense open subscheme such that $U\cap S\subseteq U$ is sparse then  $S\subseteq |X|$ is again sparse (\cite[Proposition 8.5 (a)]{poonen}) and every strictly Galois generic point is Galois generic, both the statements (for arithmetically $\NS$-generic points and for $\NS$-generic points) follow from Theorem \ref{main}, together with Fact \ref{sisagiaspec}(1).
	\item This follows Theorem \ref{main} and Fact \ref{sisagiaspec}(2)\endproof
\end{enumerate}
\begin{remark}\label{counterexemple}
	\begin{enumerate}
\item[]
\item In a first version of this paper, Fact \ref{sisagiaspec}(2) and Corollary \ref{existence}(2) were stated without requiring the  LU-property. The proof of Corollary \ref{existence}(2) was based on the use of \cite[Theorem 1.3.2]{mioUOIp}, which in turn was based on the use of \cite[Theorem 3]{mordel}, for which Akio Tamagawa recently found a counterexample. 
\item Without the LU-property Corollary \ref{existence}(2) does not hold, as the following counterexample, suggested to us by Jordan Ellenberg, shows. Assume that $k$ is infinite and finitely generated field. Let $Y\rightarrow X$ be the moduli space of elliptic curves over $k$. After possibly replacing $k$ with a finite field extension, we can chose an $x_0\in X(k)$ such that $Y_{\overline{x}_0}$ has no complex multiplication. Let $f:Y\times_k Y_{x_0}\rightarrow X$ be the product family. Then the set of non-$\NS$-generic point for $f:Y\times_k Y_{x_0}\rightarrow X$ contains the points $x\in X(k)$ such that $Y_{\overline x_0}$ and $Y_{\overline x}$ are isogenous. 
Recall that $Y_{x_0}$ is isogenous to its $p^n$-power Frobenius twists $Y^{(n)}_{x_0}$ for every $1\leq n\in \mathbb N$ and that, since $Y_{\overline x_0}$ has no complex multiplication, $Y_{\overline x_0}$ is not isomorphic to $Y^{(n)}_{x_0,\overline k}$. This shows that the set of non-$\NS$-generic points for $f:Y\times_k Y_{x_0}\rightarrow X$ contains the point $x\in X(k)$ corresponding to $Y^{(n)}_{x_0}$, so that it is infinite. Observe that the family $f:Y\times_k Y_{x_0}\rightarrow X$ has not the LU-property in degree 2. Indeed, $\widetilde \Pi$ is an open subgroup of $\SL_2(\Zl)\times \SL_2(\Zl)$, $\Pi$ is an open subgroup of $\SL_2(\Zl)$ and the inclusion $\Pi\subseteq \widetilde{\Pi}$ is induced by the inclusion in the second component $\SL_2(\Zl)\rightarrow \SL_2(\Zl)\times \SL_2(\Zl)$. In particular, $\widetilde{\Pi}/\Pi$ and $\Pi$ identify both with open subgroups of $\SL_2(\Zl)$, so that they have an infinite open subgroup in common.
\item On the other hand, the LU-property is usually satisfied in universal families. It is for example satisfied when $\Pi$ is an open subgroup of the derived subgroup of the image of $$\pi_1(X,\overline \eta)\rightarrow \GL(H^2(Y_{\overline{\eta}},\mathbb Q_{\ell}(1))).$$
Hence the LU-property is satisfied for example when $Y\rightarrow X$ is the universal family over the moduli space of abelian varieties, of K3 surfaces, etc...
\end{enumerate}
\end{remark}
Via a spreading out argument, from Corollary \ref{existence}(1) one deduces the following extension of the main result of \cite{poonen} to positive characteristic:
\begin{corollary}\label{miciomicio miao miao}
If $k$ is a field of transcendence degree $\geq 1$ over $\F_p$, then $X$ has a closed $\NS$-generic point.
\end{corollary}
\begin{remark}
Atticus Christensen (\cite[Theorem 1.0.1]{Chr}) has independently proved Corollary \ref{miciomicio miao miao}. His proof is very different from ours, since his approach is inspired from the analytic approach in \cite{poonen}, while ours is inspired from the Hodge theoretic approach in \cite{AndreIHES}. On the other hand, it seems that Corollary \ref{existence} (that will be used to prove Corollaries \ref{exinstencehyperplane},\ref{corollaryTate} and \ref{corollarybrauer}) can not be obtained via his method, that gives different information on the set of $\NS$-generic points (\cite[Theorems 1.0.3, 1.0.4.]{Chr}).
\end{remark}
From Corollary \ref{miciomicio miao miao} one easily deduces the following results on the behaviour of the Tate conjecture in families:
\begin{corollary}\label{Tate family}
If $T(Y_x,\ell)$ holds for all $x\in |X|$, then $T(Y_{\eta},\ell)$ holds.
\end{corollary}
\begin{remark}
Corollary \ref{Tate family} together with a spreading out argument can be used to reduce the Tate conjecture for smooth proper varieties over arbitrary finitely generated fields of characteristic $p$, to fields of transcendence degree one over $\F_p$, extending results from \cite{Mori}, specific to abelian schemes, to arbitrary families of varieties. 
\end{remark}
The argument in \cite[Theorem 7.1.]{poonen} shows that Corollary \ref{miciomicio miao miao} is enough to prove the following:
\begin{corollary}\label{properprojective}
Assume furthermore that $Y_{x}$ is projective for every $x\in |X|$. Then there exists a dense open subscheme $U\subseteq X$ such that the base change  $f_U:U\times_X Y\rightarrow U$ of $f:Y\rightarrow X$ through $U\subseteq X$ is projective.
\end{corollary} 
\numberwithin{equation}{subsubsection}
\begin{remark}
Whether the analogue of Corollary \ref{properprojective} holds over fields algebraic over $\F_p$ is not known.
The problem over this kind of fields is that it is not true in general that there exists a $\NS$-generic closed point (as the example of a family of abelian surfaces such that the generic fibre has not complex multiplication shows). 
\end{remark}
\numberwithin{equation}{subsubsection}
\subsubsection{Hyperplane sections}
From now on, assume that $k$ is finitely generated.
Assume that $Z$ is a smooth projective $k$-variety of dimension $\geq 3$ and let $Z\subseteq \mathbb P_k^n$ be a projective embedding. One can ask whether there exists a smooth hyperplane section $D$ of $Z$ such that the canonical map 
$$\NS(Z_{\overline k})\otimes \Q\rightarrow \NS(D_{\overline k})\otimes \Q$$
is an isomorphism. This is not true in general (see Example \ref{counter}), but one can apply Theorem \ref{main} to obtain the following arithmetic variant (see Section \ref{sechyperplane}):
\numberwithin{equation}{subsubsection}
\begin{corollary}\label{exinstencehyperplane}
If $\Dim(Z)\geq 3$ there are infinitely many smooth $k$-rational hyperplane sections $D\subseteq Z$ such that
the canonical map 
$$\NS(Z)\otimes \Q\rightarrow \NS(D)\otimes \Q$$ is an isomorphism.
\end{corollary}
\numberwithin{equation}{subsubsection}
As already mentioned in Section \ref{nerongeneric}, Conjecture \ref{tateconjectureconj} implies Theorem \ref{main}.
The $\ell$-adic Tate conjecture for divisors conjecture is still widely open except for some special classes of varieties like Abelian varieties and K3 surfaces. Using Corollary \ref{exinstencehyperplane}, one can enlarge the class of varieties for which it holds:
\begin{corollary}\label{corollaryTate}
Let $Z$ be a smooth projective $k$-variety of dimension $\geq 3$ and choose a projective embedding $Z\subseteq \mathbb P^n_k$. If $T(D,\ell)$ holds for the smooth hyperplane sections $D\subseteq Z$, then $T(Z,\ell)$ holds.
\end{corollary}
\begin{remark}
Corollary \ref{corollaryTate} can be used to reduce the $\ell$-adic Tate conjecture for divisors on smooth proper $k$-varieties to smooth projective $k$-surfaces, extending an unpublished result (\cite{djtate}) of de Jong (whose proof has been simplified in \cite[Theorem 4.3]{Varmor}) to infinite finitely generated fields.
\end{remark}
\numberwithin{equation}{subsubsection}
\subsubsection{Uniform boundedness of Brauer groups}
Combining Theorem \ref{main} with the main result of \cite{mioUOIp} and the arguments of \cite{brauer}, one gets the following application to uniform boundedness for the $\ell$-primary torsion of the cohomological Brauer group in smooth proper families of $k$-varieties (see Section \ref{secbrauer}). 
\begin{corollary}\label{corollarybrauer}
Let $X$ be a smooth geometrically connected $k$-curve and let $f:Y\rightarrow X$ be a smooth proper morphism of $k$-varieties. If $T(Y_{x},\ell)$ holds for all $x\in |X|$ and $f:Y\rightarrow X$ has the LU-property in degree 2, then there exists a constant $C:=C(Y\rightarrow X, \ell)$ such that $$|\Br(Y_{\overline x})[\ell^{\infty}]^{\pi_1(x,\overline{x})}|\leq C $$ for all $x\in X(k)$.
\end{corollary}
Corollary \ref{corollarybrauer} extends to positive characteristic the main result of \cite{brauer} and gives some evidence for a positive charateristic version of the conjectures on the uniform boundedness of Brauer group in \cite{VAV}. Elaborating the argument in the proof of Corollary \ref{corollarybrauer}, one  gets also an unconditional variant of Corollary \ref{corollarybrauer} (Corollary \ref{uncondintionalbrauer}) and a result on the specialization of the $p$-adic Tate module of the Brauer group (Corollary \ref{corptate}).
\numberwithin{equation}{subsection} 
\subsection{Organization of the Paper}
The paper is divided in three parts. Part \ref{part1} is devoted to the proof of Theorem \ref{main}: In Section \ref{strategy} we prove Theorem \ref{ltocrys} and in Section \ref{implication} we show Theorem \ref{main}. Part \ref{secapp} is devoted to applications of Theorem \ref{main}: in Section \ref{sechyperplane} we prove Corollary \ref{exinstencehyperplane}, and in Section \ref{secbrauer} we give the proof of Corollary \ref{corollarybrauer}. Finally, in Part \ref{appendix}, we prove the overconvergence of the higher direct image in crystalline cohomology (Theorem \ref{berthelotconjecture}), which is used in the proof of Theorem \ref{ltocrys}.
\subsection{Acknowledgements}
This paper is part of the author’s Ph.D. thesis under the supervision of Anna Cadoret.
He thanks her for suggesting the problem and several sessions of discussion. He is also very grateful for her careful re-readings of this paper and her constructive suggestions. The author is very grateful to Akio Tamagawa to point out a counterxemple to \cite[Theorem 3]{mordel} and to send him a preliminary version of \cite{tamagawa}. He would like also to thank Jordan Ellenberg for suggesting the counterexemple in Remark \ref{counterexemple}(2). The author is very grateful to the referees for their accurate comments, which helped clarify the exposition of the paper and the proof of Theorem \ref{ltocrys}. Several hours of discussion with Marco D'Addezio on overconvergent isocrystals and their monodromy groups has been of great help. The author would like to thanks Tomoyuki Abe for having pointed out that there might be a gap in \cite{berthelotconjecture} and for suggesting the use of Shiho's work on relative rigid cohomology to get the overconvergence of the higher direct image in crystalline cohomology. He thanks also Atsushi Shiho for kindly answering some technical questions on his work and for some comments on a preliminary draft of this paper. This circle of ideas has been presented in a small course at RIMS, a Joint Usage/Research Center located in Kyoto University. The author would like to thanks the participants for their interest and Akio Tamagawa for the invitation. A discussion around Example \ref{example} with Matilde Manzaroli has been of great help. Finally the author thanks Vincent de Daruvar and Gregorio Baldi for their interest and for some help to improve the exposition.
\part{Proof of the main theorem}\label{part1}
\section{Proof of Theorem \texorpdfstring{\ref{ltocrys}}-}\label{strategy}
This section is devoted to the proof of Theorem \ref{ltocrys}. In Section \ref{tannkianref}, after recalling the various categories of isocrystals needed in our argument, we reformulate Theorem \ref{ltocrys} in terms of overconvergent $F$-isocrystals. In Section \ref{conclusion}, we use independence techniques to prove Theorem \ref{ltocrys}.
\numberwithin{equation}{subsection} 
\subsection{Overconvergent reformulation of Theorem \texorpdfstring{\ref{ltocrys}}-}\label{tannkianref}
\numberwithin{equation}{subsubsection}
\subsubsection{Overconvergent isocrystals}\label{sectionNSover}
Let $\mathcal X$ be a smooth geometrically connected $\F_q$-variety with $q=p^s$ and write $F$ for $s$-power of the absolute Frobenius on $\mathcal X$. Write $W:=W(\F_q)$ for the Witt ring of $\F_q$, $K$ for its fraction field and   $\Mod(\mathcal X|W)$ for the category of $\mathcal O_{\mathcal X|W}$-modules in the crystalline site of $\mathcal X$.
Consider the following categories:
\begin{center}
\begin{table}[h]
\begin{tabular}{|l|l|l|l|}
\hline
Category & Notation & Name       & Reference   \\
\hline

$\Isoc(\mathcal X|W)$     & $\Isoc(\mathcal X)$  & Isocrystals  & \cite[Appendix A.1]{Varmor}        \\
\hline

$\oi(\mathcal X|K)$   & $\oi(\mathcal X)$  & Overconvergent Isocrystals & \cite[Definition 2.3.6]{berthelot}       \\
\hline

$\oi(\mathcal X,\mathcal X|K)$   & $\oi(\mathcal X,\mathcal X)$ & Convergent isocrystals & \cite[Definition 2.3.2]{berthelot} \\
\hline 
\end{tabular}
\end{table}
\end{center}
and their enriched version with Frobenius structure: $\Fisoc(\mathcal X)$, $\Foi(\mathcal X)$ and $\Foi(\mathcal X,\mathcal X)$.
They fit into the following commutative diagram:
\begin{center}
\begin{tikzcd}
& \Fisoc(\mathcal X)\arrow[r,"(1)","\simeq "']\arrow[swap]{d}{(-)^{\geo}}& \Foi(\mathcal X,\mathcal X)& \arrow[bend right,swap]{ll}{(3)}\arrow[swap]{l}{(2)}\arrow{ddl}{(4)} \Foi(\mathcal X)\arrow{d}{(-)^{\geo}}\\
& \Isoc(\mathcal X)\arrow[swap]{rd}{(5)} & & \oi(\mathcal X)\\
& & \Mod(\mathcal X|W)\otimes \Q 
\end{tikzcd}
\end{center}
where $(-)^{\geo}$ are the forgetful functors, $(1)$ is the equivalence of categories constructed in \cite[Theoreme 2.4.2]{berthelot} and $(2)$ is the obvious functor.
Write $$R^i\mathfrak f_{\cristalline,*}\mathcal O_{\mathcal Y/K}:=R^i\mathfrak f_{\cristalline,*}\mathcal O_{\mathcal Y/W}\otimes \Q\in \Mod(\mathcal X|W)\otimes \Q $$ and recall that, by \cite[Proposition A.7]{Varmor}, $R^i\mathfrak f_{\cristalline,*}\mathcal O_{\mathcal Y/K}$ is in the essential image of $(5)\circ (-)^{\geo}$. So from now on we will consider 
$R^i\mathfrak f_{\cristalline,*}\mathcal O_{\mathcal Y/K}$ as an object in $\Fisoc(\mathcal X)$.

Recall also the following fact:
\begin{fact}\cite[Theorem 1.1]{Kedfull}\label{kedlaya}
The functor $(3)$ is fully faithful.
\end{fact}
The following result, which give us an overconvergent incarnation of $R^i\mathfrak f_{\cristalline,*}\mathcal O_{\mathcal Y/K}$, is a consequence of the main result (Theorem \ref{berthelotconjecture}) of Part 3 of this paper, building on the work of Shiho on relative rigid cohomology (\cite{ShihoI}, \cite{Shiho}).
\begin{theorem}\label{Shiho}
Let $\mathfrak f:\mathcal Y\rightarrow \mathcal X$ be a smooth proper morphism. Then $R^i\mathfrak f_{\cristalline,*}\mathcal O_{\mathcal Y/K}$ in $\Fisoc(\mathcal X)$ lies in the essential image of $(3)$.
\end{theorem}
\proof
Under the equivalence $(1)$, $R^i\mathfrak f_{\cristalline,*}\mathcal O_{\mathcal Y/K}$ is sent to the Ogus higher direct image $R^if_{\Ogus,*}\mathcal O_{\mathcal Y/K}$, see \cite[Section 3, Theorem 3.1]{Ogus} and \cite[Corollary A.13]{Varmor}.
One concludes by Theorem \ref{berthelotconjecture}, which says that $R^if_{\Ogus,*}\mathcal O_{\mathcal Y/K}$ is in the image of an overconvergent $F$-isocrystal.
\endproof
Write $R^i\mathfrak f_{*}\mathcal O_{\mathcal Y|K}^{\dagger}$ for the (unique up to isomorphism) object of $\Foi(\mathcal X)$ lifting $R^i\mathfrak f_{\cristalline,*}\mathcal O_{\mathcal Y/K}$. 
\subsubsection{Overconvergent reinterpretation of Theorem \texorpdfstring{\ref{ltocrys}}-}
We now retain the notation and assumption of Theorem \ref{ltocrys}. With the notation of Theorem \ref{Shiho}, write
$$\mathcal F_p:=R^i\mathfrak f_*\mathcal O^{\dagger}_{\mathcal Y/K}(j)$$
$$\calF_{\mathcal Z, p}:=R^i\mathfrak f_{\mathcal Z,*}\mathcal O^{\dagger}_{\mathcal Y_{\mathcal Z}|K}(j)\simeq \mathfrak g^*\calF_{p}$$
where the isomorphism comes from smooth proper base change in crystalline cohomology (e.g. \cite[Proposition A.7]{Varmor}) and Fact \ref{kedlaya}.
By Fact \ref{kedlaya} one has a commutative diagram
\begin{center}
\begin{tikzcd}
H^0(\mathcal X,R^i\mathfrak f_{\cristalline,*}\mathcal O_{\mathcal Y/K})^{F=q^j}\arrow{d}\arrow{r}{\simeq} & H^0(\mathcal X,\calF^{\geo}_{p})^{F=q^j}\arrow[hook]{d}\\
H^0(\mathcal Z,R^i\mathfrak f_{\mathcal Z,\cristalline,*}\mathcal O_{\mathcal Y_{\mathcal Z}/K})^{F=q^j}\arrow{r}{\simeq} & H^0(\mathcal Z,\calF^{\geo}_{\mathcal Z, p})^{F=q^j}
\end{tikzcd}
\end{center}
where the horizontal arrows are isomorphisms.
Hence to prove Theorem \ref{ltocrys} is it enough to show that
\begin{equation}\label{claimone}
 \text{the natural injective map } 
H^0(\mathcal X,\calF^{\geo}_{p})\hookrightarrow H^0(\mathcal Z,\calF^{\geo}_{\mathcal Z, p}) \text{ is an isomorphism},
\end{equation}
since Theorem \ref{ltocrys} would then follow taking the eigenspace relative to $q^j$ of $F$ acting\footnote{Recall that $F$ is the $s$-power Frobenius, so that its action on $\oi(\mathcal X)$ is $K$-linear.} on $H^0(\mathcal X,\calF^{\geo}_{p})$ and $H^0(\mathcal Z,\calF^{\geo}_{\mathcal Z, p})$.
\subsection{Overconvergent F-isocrystals and \texorpdfstring{$\ell$-}-adic sheaves }\label{conclusion}
\subsubsection{Compatibility}
For $\ell\neq p$ write 
$$\mathcal F_{\ell}:=R^i\mathfrak f_*\mathcal \Ql(j), \quad \calF_{\mathcal Z, \ell}:=R^i\mathfrak f_{\mathcal Z,*}\Ql(j)$$
and $\calF_{\ell}^{\geo}$ (resp. $\calF_{\mathcal Z,\ell}^{\geo}$) for the restriction of $\calF_{\ell}$ (resp. $\calF_{\mathcal Z,\ell}$) to $\mathcal X_{\F}$ (resp. $\mathcal Z_{\F}$).

We now show that to prove (\ref{claimone}) it is enough to show that 
\begin{equation}\label{claimone2}
 \text{the natural injective map } H^0(\mathcal X_{\F},\calF_{\ell}^{\geo})\hookrightarrow H^0(\mathcal Z_{\F},\calF_{\mathcal Z,\ell}^{\geo}) \text{ is an isomorphism}.
\end{equation}
For this, it this is enough to show that one has 
$$\Dim(H^0(\mathcal X_{\F},\calF_{\ell}^{\geo}))=\Dim(H^0(\mathcal X,\calF_{p}^{\geo}))\quad \text{and} \quad  \Dim(H^0(\mathcal Z_{\F},\calF_{\mathcal Z, \ell}^{\geo}))=\Dim(H^0(\mathcal Z,\calF_{\mathcal Z, p}^{\geo})).$$
Recall the following.
\begin{fact}\label{compatible}\cite{Weil2}, \cite{KatzMessing},\cite{purityproper}
$\calF_{p}$,$\calF_{\ell}$ (resp. $\calF_{\mathcal Z,p}$,$\calF_{\mathcal Z,\ell}$) is a $\Q$-rational compatible system on $\mathcal X$ (resp. $\mathcal Z$)  pure of weight $i+2j$.
\end{fact}
We just prove that 
\begin{equation}\label{independencedimensioni}
\Dim(H^0(\mathcal X_{\F},\calF_{\ell}^{\geo}))=\Dim(H^0(\mathcal X,\calF_{p}^{\geo}))
\end{equation}
since the other equality is entirely similar. Let $?\in \{\ell, p\}$. Since $\mathcal F_{?}$ is pure by Fact \ref{compatible}, by the Grothendieck-Lefschetz fixed point formula (\cite[Theorem 10.5.1, page 603]{Leifu} if $?=\ell$ and \cite[Theorem 6.3]{Trace} if $?=p$) the left and the right hand sides are the number of poles, counted with multiplicity, with absolute value $q^{w/2}$ in the L-function of $\mathcal F_{?}^\vee(d)$ (see \cite[Corollaire VI.3]{lafforgue} if $?=\ell$ and \cite[Proposition 4.3.3]{abe} $?=p$), where $d$ is the dimension of $\mathcal X$. Since $\mathcal F_{\ell}$ and $\mathcal F_{p}$ are compatible, the $L$-function of $\mathcal F_{?}^\vee(d)$ does not depend on $?$, so that (\ref{independencedimensioni}) is proved. 
\subsubsection{Geometric monodromy}
We now prove (\ref{claimone2}). Let $G(\calF_{\mathcal Z,\ell},\mathfrak t)$ (resp. $G(\mathcal F_{\ell},\mathfrak t)$) denotes the Zariski closure of the image of $\pi_1(\mathcal Z,\overline{\mathfrak t})$ (resp. $\pi_1(\mathcal X,\overline{\mathfrak t})$) acting on $H^{i}(\mathcal Y_{\overline t},\Ql(j))$. Since $\pi_1(\mathcal Z,\overline{\mathfrak t})\rightarrow \pi_1(\mathcal X,\overline{\mathfrak t})\rightarrow \GL(H^{i}(\mathcal Y_{\overline{\mathfrak t}},\Ql(j)))$ is open in the image of $\pi_1(\mathcal X,\overline{\mathfrak t})\rightarrow \GL(H^{i}(\mathcal Y_{\overline{\mathfrak t}},\Ql(j)))$ and $G(\mathcal F_{\ell},\mathfrak t)$ is connected, one has   
$$G(\calF_{\mathcal Z,\ell},\mathfrak t)=G(\mathcal F_{\ell},\mathfrak t).$$

Let $G^{\geo}(\calF_{\ell},\mathfrak t)$ (resp. $G^{\geo}(\mathcal F_{\mathcal Z,\ell},\mathfrak t)$) for the Zariski closure of the image of $\pi_1(\mathcal X_{\mathbb F},\mathfrak t)$ (resp. $\pi_1(\mathcal Z_{\mathbb F},\mathfrak t)$) acting on $H^{i}(\mathcal Y_{\overline t},\Ql(j))$. To prove (\ref{claimone2}), it is enough to show that $G^{\geo}(\calF_{\mathcal Z,\ell},\mathfrak t)=G^{\geo}(\calF_{\ell},\mathfrak t)$. Recall the following:
\begin{fact}\label{reductiontogeometric}
The groups $G^{\geo}(\calF_{\ell},\mathfrak t)^0$ and $G^{\geo}(\calF_{\mathcal Z,\ell},\mathfrak t)^0$ are semisimple algebraic groups.
\end{fact}
\proof
Since $\mathcal F_{\ell}$ and $\mathcal F_{\mathcal Z,\ell}$ are pure, $\mathcal F^{\geo}_{\ell}$ and $\mathcal F^{\geo}_{\mathcal Z,\ell}$ are semisimple by \cite[Theorem 3.4.1(iii)]{Weil2}. So one apply \cite[Corollaire 1.3.9]{Weil2} to conclude.
\endproof
By assumption $G(\mathcal F_{\ell},\mathfrak t)=G(\mathcal F_{\mathcal Z,\ell},\mathfrak t)$ and $G^{\geo}(\mathcal F_{\ell},\mathfrak t)$ are connected, so that it is enough to show that 
$$Q:=G^{\geo}(\mathcal F_{\ell},\mathfrak t)/G^{\geo}(\mathcal F_{\mathcal Z,\ell},\mathfrak t)$$ is finite. 
The fundamental exact sequences for $\mathcal X$ and $\mathcal Z$ induces a commutative diagram with exact rows
 \begin{center}
\begin{tikzcd}
0\arrow{r} & \pi_1(\mathcal Z_{\F},\overline{\mathfrak t})\arrow{r}\arrow{d} &\pi_1(\mathcal Z,\overline{\mathfrak t})\arrow{r}\arrow{d} & \pi_1(\F_q)\arrow[equal]{d}\arrow{r}& 0\\
0\arrow{r} & \pi_1(\mathcal X_{\F},\overline{\mathfrak t})\arrow{r} &\pi_1(\mathcal X,\overline{\mathfrak t})\arrow{r} & \pi_1(\F_q)\arrow{r}& 0,
\end{tikzcd}
\end{center}
hence we get a commutative diagram with exact rows
 \begin{center}
\begin{tikzcd}
0\arrow{r} & G^{\geo}(\mathcal F_{\mathcal Z,\ell},\mathfrak t)\arrow{r}\arrow[hook]{d} &G(\mathcal F_{\mathcal Z,\ell},\mathfrak t)\arrow{r}\arrow[equal]{d} & Q_{\mathcal Z}\arrow[two heads]{d}\arrow{r}& 0\\
0\arrow{r} &  G^{\geo}(\mathcal F_{\ell},\mathfrak t)\arrow{r} &G(\mathcal F_{\ell},\mathfrak t)\arrow{r} & Q_{\mathcal X}\arrow{r}& 0,
\end{tikzcd}
\end{center}
where $Q_{\mathcal Z}$ and $Q_{\mathcal X}$ are abelian algebraic groups. This shows that $G^{\geo}(\mathcal F_{\mathcal Z,\ell},\mathfrak t)\subseteq G(\mathcal F_{\ell},\mathfrak t)$ is a normal subgroup, hence $G^{\geo}(\mathcal F_{\mathcal Z,\ell},\mathfrak t)\subseteq G^{\geo}(\mathcal F_{\ell},\mathfrak t)$ is also a normal subgroup. The snake lemma shows then that
$$Q=\Ker(Q_{\mathcal Z}\rightarrow Q_{\mathcal X}),$$
hence that $Q$ is an abelian algebraic group. But since, by Fact \ref{reductiontogeometric}, $G^{\geo}(\mathcal F_{\ell},\mathfrak t)^0$ is a semisimple algebraic group, this implies that $Q$ is finite and concludes the proof of Theorem \ref{ltocrys}.
\endproof
\section{Proof of Theorem \texorpdfstring{\ref{main}}- }\label{implication}
\numberwithin{equation}{subsection} 
In Section \ref{prelim}, we collect some preliminary remarks. The proof when $f:Y\rightarrow X$ is proper is a technical elaboration (involving alteration and the trace formalism) of the proof when $f:Y\rightarrow X$ is projective. To clarify the exposition we carry out the proof when $f:Y\rightarrow X$ is projective in Section \ref{problectionsection} and turn to the general case in Section \ref{seccrys}.  
\subsection{Preliminary remarks}\label{prelim}
\subsubsection{Strictly generic vs generic}
Observe that the assertion for Galois-generic points implies the assertion for strictly Galois-generic points. Indeed, strictly Galois-generic implies Galois-generic, hence for a strictly Galois-generic point $x\in X$ the specialization morphism
$$\speci_{\eta,x}:\NS(Y_{\overline \eta})\otimes \Q\rightarrow  \NS(Y_{\overline x})\otimes \Q$$
is an isomorphism. Recall that, as explained in \ref{nerongeneric}, the map $\speci_{\eta,x}$ is $\pi_1(x,\overline x)$-equivariant. Since $\pi_1(x,\overline x)$ and $\pi_1(X,\overline \eta)\simeq \pi_1(X,\overline x)$ acting on $H^2(Y_{\overline \eta},\Ql(1))\simeq H^2(Y_{\overline x},\Ql(1))$ have the same image $\Pi_{\ell}$ (since $x$ is strictly Galois-generic), taking $\Pi_{\ell}$-invariants in $\speci_{\eta,x}$, one deduces the statement for strictly Galois-generic points. So, from on, we focus on the assertion for Galois-generic points. To simplify, in this section, we omit base points in our notation for the étale fundamental group. 
\subsubsection{Finite cover}\label{etale cover}
If $X'\rightarrow X$ is a surjective finite morphism of smooth connected $k$-varieties, the map $\pi_1(X')\rightarrow \pi_1(X)$ has open image. So $x\in X$ is Galois-generic (resp. $\NS$-generic) for $f:Y\rightarrow X$ if and only if any lifting $x'\in X'$ of $x$ if Galois-generic (resp. $\NS$-generic) for the base change $f_{X'}:Y'\times_XX'\rightarrow X'$ of $f:Y\rightarrow X$ along $X'\rightarrow X$. As a consequence we can freely replace $X$ with $X'$ during the proof.
\subsection{Proof when \texorpdfstring{$f$}- is projective}\label{problectionsection}
Let $f:Y\rightarrow X$ be smooth projective. For the general strategy of the proof see Section \ref{introspreadingNS}.
\numberwithin{equation}{subsubsection}
\subsubsection{Step 1: Spreading out}\label{spreading}
Replacing $k$ with a finite field extension (\ref{etale cover}), one can assume that there exists a finite field $\F_q$, smooth and geometrically connected $\F_q$-varieties $\mathcal K$, $\mathcal Z$, $\mathcal X$ with generic points $\zeta:k\rightarrow \mathcal K$, $\beta:k(x)\rightarrow \mathcal Z$ and a commutative diagram in which the squares indicated with a $\Box$ are cartesian 
 \begin{center}
\begin{tikzcd}
\mathcal Y_{\mathfrak t} \arrow{r}\arrow{d}\arrow[phantom]{rd}{\Box}&\mathcal Y_{\mathcal Z}\arrow{r}\arrow{d}{\mathfrak f_{\mathcal Z}}\arrow[phantom]{rd}{\Box}& \mathcal Y\arrow{d}{\mathfrak f}\arrow[phantom]{rd}{\Box} & Y\arrow{l}\arrow{d}{f}\arrow[phantom]{rd}{\Box}&Y_x\arrow{l}\arrow{d}{f_x}\arrow[bend right]{lll}\\

\F_q \arrow{r}{\mathfrak t}\arrow[equal]{dr}&\mathcal Z\arrow{d}\arrow{r}\arrow{rd}&\mathcal X\arrow[phantom]{rd}{\Box}\arrow{d}& X\arrow{l}\arrow{d}&k(x)\arrow[swap]{l}{x}\arrow[bend left,near start]{lll}{\beta}\\
              & \F_q & \arrow{l}\mathcal K & k\arrow{l}{\zeta}
\end{tikzcd}
\end{center}
where $\mathfrak f:\mathcal Y\rightarrow \mathcal X$ is a smooth projective morphism and the base change of $f_{\mathcal Z}:\mathcal Y_{\mathcal Z}\rightarrow \mathcal Z$ along $\beta:k(x)\rightarrow \mathcal Z$ identifies with $f_x:Y_x\rightarrow k(x)$. 
Replacing $X$ with a finite étale cover (\ref{etale cover}) one can also assume that
\begin{enumerate}
\item $\NS(Y_{\overline x})\otimes \Q=\NS(Y_{x})\otimes \Q$ and $\NS(\mathcal Y_{\overline{\mathfrak t}})\otimes \Q=\NS(\mathcal Y_{\mathfrak t})\otimes \Q$;
\item the Zariski closures of the images of $\pi_1(\mathcal X)\rightarrow \GL(H^{2}(Y_{\overline \eta},\Ql(1)))$ and $\pi_1(\mathcal X_{\F})\rightarrow  \GL(H^{2}(Y_{\overline \eta},\Ql(1)))$ are connected.
\end{enumerate}
Note that, by smooth proper base change, one has the following factorization
\begin{center}
\begin{tikzcd}
\pi_1(x)\arrow{r}\arrow[two heads]{d} &\pi_1(X)\arrow[two heads]{d}\arrow{r} &\GL(H^{2}(Y_{\overline \eta},\Ql(1)))\simeq \GL(H^{2}(\mathcal Y_{\overline t},\Ql(1)))\\
\pi_1(\mathcal Z)\arrow{r}&\pi_1(\mathcal X)\arrow{ru}.
\end{tikzcd}
\end{center}
In particular, since $x$ is Galois-generic with respect to $f:Y\rightarrow X$, the image of $\pi_1(\mathcal Z)\rightarrow \pi_1(\mathcal X)\rightarrow \GL(H^{2}(\mathcal Y_{\overline t},\Ql(1)))$ is open in the image of $\pi_1(\mathcal X)\rightarrow \GL(H^{2}(\mathcal Y_{\overline t},\Ql(1)))$.
Hence by \ref{spreading}(2) and Theorem \ref{ltocrys} the base change map
$$H^0(\mathcal X,R^2\mathfrak f_{\cristalline,*}\mathcal O_{\mathcal Y/W})^{F=q}\otimes \Q\rightarrow H^0(\mathcal Z,R^2\mathfrak f_{\mathcal Z,\cristalline,*}\mathcal O_{\mathcal Y_{\mathcal Z}/W})^{F=q}\otimes \Q$$
is an isomorphism.
\numberwithin{equation}{subsubsection}
\subsubsection{Step 2: Using the variational Tate conjecture}\label{diagram}
Since $\mathfrak t$ is a specialization of $x$ (in $\mathcal Z$) and $x$ is a specialization of $\eta$ (in $\mathcal X$), there is a canonical commutative diagram
\begin{center}
\begin{tikzcd}
z\in \Pic(\mathcal Y)\otimes \Q \arrow[bend right=65,swap]{dddd}{c_{\mathcal Y}}\arrow{rr}\arrow{d}\arrow{dr}{i_{\mathfrak t}^*}& & \Pic(\calY_{\mathcal Z})\otimes \Q\ni z_x\arrow[two heads]{d}{(i)}\arrow[bend left=65]{dddd}{c_{\mathcal Y_{\mathcal Z}}} \arrow[swap]{dl}{i_{\mathfrak t}^*}\\
\Pic(Y_{\eta})\otimes \Q\arrow{d}& \Pic(\mathcal Y_{\mathfrak t})\otimes \Q\arrow[near start,two heads]{dd}{(ii)}\arrow[bend left=100]{dddd}{c_{\mathcal Y_{\mathfrak t}}}\ni z_t & \Pic(Y_x)\otimes \Q\arrow[two heads]{d}{(ii)}\\
\varepsilon_{\eta} \in \NS(Y_{\overline{\eta}})\otimes \Q\arrow[hook, near start]{rr}{\speci_{\eta,x}}\arrow[hook]{dr}{\speci_{\eta,\mathfrak t}}& & \NS(Y_{\overline x})\otimes \Q \ni \varepsilon_x \arrow[hook,swap]{dl}{\speci_{x,\mathfrak t}}\\
& \NS(\mathcal Y_{\mathfrak t})\otimes \Q=\NS(\mathcal Y_{\overline{\mathfrak t}})\arrow[hook]{dd}\otimes \Q \ni \varepsilon_{\mathfrak t} &\\

H^2_{\cristalline}(\mathcal Y)^{F=q}\arrow{rd}{i_{\mathfrak t}^*}& & H^2_{\cristalline}(\mathcal Y_{\mathcal Z})^{F=q}\arrow[near end, swap]{ld}{i_{\mathfrak t}^*}\\
& H^2_{\cristalline}(\mathcal Y_{\mathfrak t})^{F=q}, &
\end{tikzcd}
\end{center}
where the arrow $(i)$ is surjective, since an open immersion of smooth varieties induces a surjection on the Picard groups, and the arrows $(ii)$ are surjective by \ref{spreading}$(1)$. Take an $\varepsilon_x$ in $\NS(Y_{\overline x})\otimes \Q$ with lifting $z_x\in \Pic(\mathcal Y_{\mathcal Z})\otimes \Q$ and write 
$$z_{\mathfrak t}:=i_{\mathfrak t}^*(z_x)\in \Pic(\mathcal Y_{\mathfrak t})\otimes \Q \quad \text{and} \quad \varepsilon_{\mathfrak t}=\speci_{x,\mathfrak t}(\varepsilon_x)\in \NS(\mathcal Y_{\mathfrak t})\otimes \Q.$$
Since $\speci_{x,\mathfrak t}: \NS(Y_{\overline{\eta}})\otimes \Q\rightarrow \NS(Y_{\overline{t}})\otimes \Q$ is injective, it is enough to show that there exists a $z\in \Pic(\mathcal Y)\otimes \Q$  such that $i_{\mathfrak t}^*c_{\mathcal Y}(z)=c_{\mathcal Y_{\mathfrak t}}i_{\mathfrak t}^*(z_x)$. Indeed, if  $\varepsilon_{\eta}$ denotes the image of $z$ in $\NS(Y_{\overline{\eta}})\otimes \Q$, one would then have 
$$\speci_{x,\mathfrak t}(\speci_{\eta,x}(\varepsilon_{\eta}))=\speci_{\eta,\mathfrak t}(\varepsilon_{\eta})=\varepsilon_{\mathfrak t}=\speci_{x,\mathfrak t}(\varepsilon_x)$$
hence $\speci_{\eta,x}(\varepsilon_{\eta})=\varepsilon_{x}$.

To construct such $z\in \Pic(\mathcal Y)\otimes \Q$, consider the commutative diagram
\begin{center}
\begin{tikzcd}
z\in \Pic(\mathcal Y)\otimes \Q \arrow[swap]{dd}{c_{\mathcal Y}}\arrow{rr}\arrow[swap]{dr}{i_{\mathfrak t}^*}& & \Pic(\calY_{\mathcal Z})\otimes \Q\ni z_x\arrow{dd}{c_{\mathcal Y_{\mathcal Z}}} \arrow{dl}{i_{\mathfrak t}^*}\\
& \Pic(\mathcal Y_{\mathfrak t})\otimes \Q\arrow{dd}{c_{\mathcal Y_{\mathfrak t}}}\ni z_t & \\
H^2_{\cristalline}(\mathcal Y)^{F=q}\arrow{rd}{i_{\mathfrak t}^*}\arrow[swap]{dd}{\Leray} & & H^2_{\cristalline}(\mathcal Y_{\mathcal Z})^{F=q}\arrow{dd}{\Leray}\arrow[near end, swap]{ld}{i_{\mathfrak t}^*}\\
& H^2_{\cristalline}(\mathcal Y_{\mathfrak t})^{F=q} &\\
H^0(\mathcal X,R^2\mathfrak f_{\cristalline,*}\mathcal O_{\mathcal Y/W})^{F=q}\otimes \Q\arrow{rr}{(iii)}\arrow[hook]{ru}& &H^0(\mathcal Z,R^2\mathfrak f_{\mathcal Z,\cristalline,*}\mathcal O_{\mathcal Y_{\mathcal Z}/W})^{F=q}\otimes \Q\arrow[hook]{lu}
\end{tikzcd}
\end{center}
where the arrow $(iii)$ is an isomorphism by Theorem \ref{ltocrys}. By construction $c_{\mathcal Y_{\mathfrak t}}i_{\mathfrak t}^*(z_x)$ is in the image of $$(iii):H^0(\mathcal X,R^2\mathfrak f_{\cristalline,*}\mathcal O_{\mathcal Y/W})^{F=q}\otimes \Q\xrightarrow{\sim} H^0(\mathcal Z,R^2\mathfrak f_{\mathcal Z,\cristalline,*}\mathcal O_{\mathcal Y_{\mathcal Z}/W})^{F=q}\otimes \Q . $$ 
By Fact \ref{thmvariationalTate} applied to $\mathfrak f:\mathcal Y\rightarrow \mathcal Z$ and $\mathfrak t$, there exists $z\in \Pic(\mathcal Y)\otimes \Q$ such that $i_{\mathfrak t}^*c_{\mathcal Y}(z)=c_{\mathcal Y_{\mathfrak t}}i_{\mathfrak t}^*(z_x)$. This concludes the proof of Theorem \ref{main} when $f:Y\rightarrow X $ is projective.\endproof
\subsection{Proof when \texorpdfstring{$f$}- is proper}\label{seccrys}
Assume now that $f:Y\rightarrow X $ is only proper. Since Fact \ref{thmvariationalTate} is only available when $f$ is projective, we can no longer apply it directly to $\mathfrak f:\mathcal Y\rightarrow \mathcal X$. To overcome this difficulty we proceed as follows. Using de Jong's alteration theorem and replacing $X$ with a finite cover of a dense open subset, one first constructs a commutative diagram
\begin{center}
\begin{tikzcd}
\widetilde Y\arrow{rr}{g}\arrow[swap]{dr}{\widetilde f}&& Y\arrow{dl}{f}\\
& X
\end{tikzcd}
\end{center}
with $\widetilde f$ smooth projective and $g$ dominant and generally finite. While every $x\in X$ which is $\NS$-generic for $\widetilde f$ is $\NS$-generic for $f$ (as the argument in \ref{trace} will show), the hypothesis of being Galois-generic for $f$ does not transfer to $\widetilde f$ in general, so that one can not reduce directly the assertion for the (proper) morphism $f:Y\rightarrow X$ to the assertion for the (projective) morphism $\widetilde f:\widetilde Y\rightarrow X$. However, the trace formalism is functorial enough to allow us to transfer information from $f:Y\rightarrow X$ to $\widetilde f:\widetilde Y\rightarrow X$ for cohomology classes coming from $Y$.
\subsubsection{Step 1: de Jong's alterations theorem}
First one reduces to the situation where$f$ has geometrically connected fibres (this hypothesis is used in \ref{trace} to apply Poincar\'e duality). By \cite[X, Proposition 1.2]{SGA1} and replacing $X$ with a finite étale cover (\ref{etale cover}), one can assume that $f:Y\rightarrow X$ decomposes in a disjoint union of morphisms $f_i:Y_i\rightarrow X$ with geometrically connected fibres.
Since for every (not necessarily closed) point $x\in X$ there are natural decompositions
\begin{center}
\begin{tikzcd}
\NS(Y_{\overline x})\otimes\Q\arrow{r}\arrow{d}{\simeq}&H^2(Y_{\overline x},\Ql(i))\arrow{d}{\simeq}\\
\oplus_i \NS(Y_{i,\overline x})\otimes\Q\arrow{r}&\oplus_i H^2(Y_{i,\overline x},\Ql(i))
\end{tikzcd}
\end{center}
one may work with each $f_i:Y_i\rightarrow X$ separately and hence assume that $f:Y\rightarrow X$ has geometrically connected fibres.

By de Jong's alterations theorem (\cite{alteration}) for $Y_{\overline \eta}$ over $\overline{k(\eta)}$, there exists a proper, surjective and generically finite morphism $\widetilde Y_{\overline \eta}\rightarrow Y_{\overline \eta}$, where $\widetilde Y_{\overline \eta}$ is a connected, smooth and projective $\overline{k(\eta)}$-variety. By descent and spreading out, there exists a commutative diagram of connected smooth $k$-varieties:
\begin{center}
\begin{tikzcd}
\widetilde Y_{\overline \eta}\arrow{r}\arrow{d}\arrow[phantom]{rd}{\Box} & \widetilde Y_{\eta'}\arrow[phantom]{rd}{\Box}\arrow{r}\arrow{d} &\widetilde Y\arrow{d}{g}\arrow[bend right, swap,near start]{dd}{\widetilde f}\\
Y_{\overline \eta}\arrow{r}\arrow{d}\arrow[phantom]{rd}{\Box} & Y_{\eta'}\arrow{r}\arrow{d}\arrow[phantom]{rd}{\Box}&Y_{U'}\arrow{r}\arrow{d}{f_{U'}}\arrow{d}\arrow[phantom]{rd}{\Box} & Y_U\arrow{r}\arrow{d}{f_U}\arrow[phantom]{rd}{\Box} & Y\arrow{d}{f}\\
\overline{k(\eta)} \arrow{r} & k(\eta')\arrow{r}{\eta'} & U'\arrow{r}{j}&U\arrow{r}{i}& X
\end{tikzcd}
\end{center}
where $\eta':k(\eta')\rightarrow U'$ is the generic point of $U'$, $i:U\rightarrow X$ is a open immersion with dense image, $j:U'\rightarrow U$ is a finite surjective morphism, $\widetilde f:\widetilde Y\rightarrow U'$ is smooth, projective with geometrically connected fibres and $g:\widetilde Y\rightarrow Y_{U'}$ is proper, surjective and generically finite. In conclusion, replacing $X$ with $U'$ (\ref{etale cover}), one can assume that there exists a diagram
\begin{center}
\begin{tikzcd}
\widetilde Y\arrow{rr}{g}\arrow[swap]{dr}{\widetilde f}&& Y\arrow{dl}{f}\\
& X
\end{tikzcd}
\end{center}
where $\widetilde f:\widetilde Y\rightarrow X$ is smooth projective with geometrically connected fibres, $f:Y\rightarrow X$ is smooth proper with geometrically connected  fibres and $g:\widetilde Y\rightarrow Y$ is generically finite and dominant.
\numberwithin{equation}{subsubsection}
\subsubsection{Step 2: Spreading out}
Now one spreads out to finite fields. Up to replacing $k$ with a finite field extension (\ref{etale cover}), there exists a finite field $\F_q$, smooth and geometrically connected $\F_q$-varieties $\mathcal K$,$\mathcal Z$ with generic points $\zeta:k\rightarrow \mathcal K$, $\beta:k(x)\rightarrow \mathcal Z$ and a commutative diagram in which the squares indicated with a $\Box$ are cartesian:
 \begin{center}
\begin{tikzcd}
\widetilde{\mathcal Y}_{\mathfrak t}\arrow{r}\arrow{d}{\mathfrak g_{\mathfrak t}}\arrow[phantom]{rd}{\Box}\arrow[bend right,swap,near start]{dd}{\widetilde{\mathfrak f}_{\mathfrak t}} & \widetilde{\mathcal Y}_{\mathcal Z}\arrow{r}\arrow{d}{\mathfrak g_{\mathcal Z}}\arrow[phantom]{rd}{\Box} \arrow[bend right,swap,near start]{dd}{\widetilde{\mathfrak f}_{\mathcal Z}}& \widetilde{\mathcal Y}\arrow{d}{\mathfrak g} \arrow[phantom]{rd}{\Box}\arrow[bend right,swap,near start]{dd}{\widetilde{\mathfrak f}}& \arrow{l} \widetilde Y\arrow[bend right, swap, near start]{dd}{\widetilde f}\arrow{d}{g}\arrow[phantom]{rd}{\Box}& \arrow{l}\widetilde Y_x\arrow{d}{g_x}\\
\mathcal Y_{\mathfrak t} \arrow{r}\arrow{d}\arrow[phantom]{rd}{\Box}&\mathcal Y_{\mathcal Z}\arrow{r}\arrow{d}{\mathfrak f_{\mathcal Z}}\arrow[phantom]{rd}{\Box}& \mathcal Y\arrow{d}{\mathfrak f}\arrow[phantom]{rd}{\Box} & Y\arrow{l}\arrow{d}{f}\arrow[phantom]{rd}{\Box}&Y_x\arrow{l}\arrow{d}{f_x}\\

\F_q \arrow{r}{\mathfrak t}\arrow[equal]{dr}&\mathcal Z\arrow{d}\arrow{r}&\mathcal X\arrow[phantom]{rd}{\Box}\arrow{d}& X\arrow{l}\arrow{d}&k(x)\arrow[swap]{l}{x}\arrow[bend left]{lll}{\beta}\\
              & \F_q & \arrow{l}\mathcal K & k\arrow{l}{\zeta}

\end{tikzcd}
\end{center}
where $\mathfrak f:\mathcal Y\rightarrow \mathcal X$ is smooth proper with geometrically connected fibres, $\widetilde{\mathfrak f}:\widetilde{\mathcal Y}\rightarrow \mathcal X$ is smooth projective with geometrically connected fibres, $\mathfrak{g}:\widetilde{\mathcal Y}\rightarrow \mathcal Y$ is a dominant generically finite morphism and the base change of $\widetilde{\mathcal Y}_{\mathcal Z}\rightarrow \mathcal Y_{\mathcal Z}\rightarrow \mathcal Z$ along $k(x)\rightarrow \mathcal Z$ identifies with $\widetilde Y_x\rightarrow Y_x\rightarrow k(x)$. 
Replacing $X$ with a finite étale cover (\ref{etale cover}) one can also assume that
\begin{enumerate}
\item $\NS(Y_{\overline x})\otimes \Q=\NS(Y_{x})\otimes \Q$, $\NS(\mathcal Y_{\overline{\mathfrak t}})\otimes \Q=\NS(\mathcal Y_{\mathfrak t})\otimes \Q$, $\NS(\widetilde{\mathcal Y}_{\mathfrak t})\otimes \Q=\NS(\widetilde{\mathcal Y}_{\overline{\mathfrak t}})\otimes \Q$;
\item the Zariski closures of the images of $\pi_1(\mathcal X)\rightarrow \GL(H^{2}(Y_{\overline \eta},\Ql(1)))$ and $\pi_1(\mathcal X_{\F})\rightarrow  \GL(H^{2}(Y_{\overline \eta},\Ql(1)))$ are connected.
\end{enumerate}
Note that, by smooth proper base change, one has the following factorization
\begin{center}
\begin{tikzcd}
\pi_1(x)\arrow{r}\arrow[two heads]{d} &\pi_1(X)\arrow[two heads]{d}\arrow{r} &\GL(H^{2}(Y_{\overline \eta},\Ql(1)))\simeq \GL(H^{2}(\mathcal Y_{\overline t},\Ql(1)))\\
\pi_1(\mathcal Z)\arrow{r}&\pi_1(\mathcal X)\arrow{ru}
\end{tikzcd}
\end{center}
In particular, since $x$ is Galois-generic for $f:Y\rightarrow X$, the image of $\pi_1(\mathcal Z)\rightarrow \pi_1(\mathcal X)\rightarrow \GL(H^{2}(\mathcal Y_{\overline t},\Ql(1)))$ is open in the image of $\pi_1(\mathcal X)\rightarrow \GL(H^{2}(\mathcal Y_{\overline t},\Ql(1)))$.
Hence by (2) and Theorem \ref{ltocrys} the base change map
$$H^0(\mathcal X,R^2\mathfrak f_{\cristalline,*}\mathcal O_{\mathcal Y/W})^{F=q}\otimes \Q\rightarrow H^0(\mathcal Z,R^2\mathfrak f_{\mathcal Z,\cristalline,*}\mathcal O_{\mathcal Y_{\mathcal Z}/W})^{F=q}\otimes \Q$$
is an isomorphism.

\numberwithin{equation}{subsubsection}
\subsubsection{Step 3: Using the Variational Tate conjecture}
Take an element $\varepsilon_x$ in $\NS(Y_x)\otimes \Q$. The goal of this subsection is to prove that there exists a $\widetilde{\varepsilon}_{\eta}\in \NS(\widetilde{Y}_{\overline \eta})\otimes \Q$ such that $\widetilde{\speci}_{\eta,x}(\widetilde{\varepsilon}_{\eta})=g^*(\varepsilon_x)$, where 
$$\widetilde{\speci}_{\eta,x}:\NS(\widetilde{Y}_{\overline \eta})\otimes \Q\rightarrow \NS(\widetilde{Y}_{\overline x})\otimes \Q$$
is the specialization map for $\widetilde f:\widetilde Y\rightarrow X$. Consider the first commutative diagram in \ref{diagram}.
Let $z_x\in \Pic(\mathcal Y_{\mathcal Z})\otimes \Q$ be a lift of $\varepsilon_x$ and write 
$$z_{\mathfrak t}:=i_{\mathfrak t}^*(z_x)\in \Pic(\mathcal Y_{\mathfrak t})\otimes \Q \quad \text{and} \quad \varepsilon_{\mathfrak t}=\speci_{x,\mathfrak t}(\varepsilon_x)\in \NS(\mathcal Y_{\mathfrak t})\otimes \Q.$$
By construction $c_{\mathcal Y_{\mathfrak t}}i_{\mathfrak t}^*(z_x)$ is in the image of 
$$(iii):H^0(\mathcal X,R^2\mathfrak f_{\cristalline,*}\mathcal O_{\mathcal Y/W})^{F=q}\otimes \Q \xrightarrow{\sim} H^0(\mathcal Z,R^2\mathfrak f_{\mathcal Z,\cristalline,*}\mathcal O_{\mathcal Y_{\mathcal Z}/W})^{F=q}\otimes \Q.$$
Since $\mathfrak f:\mathcal Y\rightarrow \mathcal X$ is only assumed to be proper, one cannot apply directly Fact \ref{thmvariationalTate} to it.
However the previous reasoning shows that 
$$H^0(\mathcal X, R^2\widetilde{\mathfrak f}_{*,\cristalline}\mathcal O_{\widetilde{\mathcal Y}/W})\otimes \Q \supseteq \mathfrak g^*(H^0(\mathcal X,R^2\mathfrak f_{\cristalline,*}\mathcal O_{\mathcal Y/W})\otimes \Q)\ni \mathfrak g^*_{\mathfrak t}(c_{\mathcal Y_{\mathfrak t}}i_{\mathfrak t}^*(z_x))=c_{\widetilde{\mathcal Y}_{\mathfrak t}}\widetilde i_{\mathfrak t}^*\mathfrak g^*_{\mathcal Z}(z_x),$$
where the notation is as in the canonical commutative diagram:
\begin{center}
\begin{tikzpicture}[baseline= (a).base]
\node[scale=0.95] (a) at (0,0){
\begin{tikzcd}
\widetilde z\in \Pic(\widetilde{\mathcal Y})\otimes \Q \arrow[swap]{dd}{c_{\widetilde{\mathcal Y}}}\arrow{rr}\arrow{dr}{\widetilde i_{\mathfrak t}^*}& & \Pic(\widetilde{\calY}_{\mathcal Z})\otimes \Q\ni \mathfrak g_{\mathcal Z}^*(z_x)\arrow{dd}{c_{\widetilde{\mathcal Y}_{\mathcal Z}}}\arrow[swap]{dl}{\widetilde i_{\mathfrak t}^*}\\

& \Pic(\widetilde{\mathcal Y}_{\mathfrak t})\otimes \Q\arrow{dd}{c_{\widetilde{\mathcal Y}_{\mathfrak t}}}\ni z_t & \\

H^2_{\cristalline}(\widetilde{\mathcal Y})^{F=q}\arrow{rd}{\widetilde i_{\mathfrak t}^*}\arrow[swap]{dd}{\Leray} & & H^2_{\cristalline}(\widetilde{\mathcal Y}_{\mathcal Z})^{F=q}\arrow{dd}{\Leray}\arrow[near end, swap]{ld}{\widetilde i_{\mathfrak t}^*}\\
& H^2_{\cristalline}(\widetilde{\mathcal Y}_{\mathfrak t})^{F=q} &\\

\mathfrak g^*(\varepsilon_{\mathfrak t})\in H^0(\mathcal X,R^2\widetilde{\mathfrak f}_{\cristalline,*}\mathcal O_{\widetilde{\mathcal Y}/W})^{F=q}\otimes \Q\arrow[hook]{ur}\arrow[hook]{rr}& &H^0(\mathcal Z,R^2\widetilde{\mathfrak f}_{\mathcal Z,\cristalline,*}\mathcal O_{\widetilde{\mathcal Y}_{\mathcal Z}/W})^{F=q}\otimes \Q\arrow[hook]{ul}\ni \mathfrak g_{\mathcal Z}^*(\varepsilon_{\mathfrak t})\\
\varepsilon_{\mathfrak t}\in H^0(\mathcal X,R^2\mathfrak f_{\cristalline,*}\mathcal O_{\mathcal Y/W})^{F=q}\otimes \Q\arrow{u}{\mathfrak g^*}\arrow{rr}& &H^0(\mathcal Z,R^2\mathfrak f_{\mathcal Z,\cristalline,*}\mathcal O_{\mathcal Y_{\mathcal Z}/W})^{F=q}\otimes \Q\arrow{u}{\mathfrak g_{\mathcal Z}^*}\ni \varepsilon_{\mathfrak t}
\end{tikzcd}
};
\end{tikzpicture}
\end{center}
So, by Fact \ref{thmvariationalTate} applied to $\widetilde{\mathfrak f}:\widetilde{\mathcal Y}\rightarrow \mathcal X$ and $\mathfrak t$, there exists $\widetilde z\in \Pic(\widetilde{\mathcal Y})\otimes \Q$ such that $\mathfrak g^*_{\mathfrak t}(c_{\mathcal Y_{\mathfrak t}}i_{\mathfrak t}^*(z_x))=\widetilde i_{\mathfrak t}^*c_{\widetilde{\mathcal Y}}(\widetilde z)$.
Let  $\widetilde{\varepsilon}_{\eta}$ be the image of $\widetilde i_{\eta}^*(\widetilde z)$ in $\NS(\widetilde Y_{\overline \eta})\otimes \Q$.
From the commutative diagram
\begin{center}
\begin{tikzcd}[row sep=tiny]
\Pic(\widetilde{\mathcal Y})\otimes \Q\arrow{dd} \arrow{r} &\Pic(\widetilde Y_{\eta})\otimes \Q\arrow{r}& \NS(\widetilde Y_{\overline \eta})\otimes \Q\arrow{dd}{\widetilde \speci_{\eta,\mathfrak t}}\arrow{rd}{\widetilde \speci_{\eta,x}}\\
&&&\NS(\widetilde Y_{\overline{x}})\arrow{ld}{\widetilde \speci_{x,\mathfrak t}}\\
\Pic(\widetilde{\mathcal Y}_{\mathfrak t})\arrow{rr} && \NS(\widetilde{\mathcal Y}_{\overline{\mathfrak t}})\otimes \Q
\end{tikzcd}
\end{center}
one deduces $$\widetilde{\speci}_{x,\mathfrak t}(\widetilde{\speci}_{\eta,x}(\widetilde{\varepsilon}_{\eta}))=\widetilde{\speci}_{x,\mathfrak t}(g^*\varepsilon_{x}).$$
Since $\widetilde{\speci}_{x,\mathfrak t}$ is injective this implies $\widetilde{\speci}_{\eta,x}(\widetilde{\varepsilon}_{\eta})=g^*(\varepsilon_{x})$.
\numberwithin{equation}{subsubsection}
\subsubsection{Step 4: Trace argument}\label{trace}
To conclude the proof one has to descend from $\widetilde{Y}$ to $Y$. For this we use the trace formalism.
Since $f:Y\rightarrow X$ and $\widetilde f:\widetilde Y\rightarrow X$ are smooth proper morphisms with geometrically connected fibres, by the relative Poincar\'e duality (\cite[Expos\'e XVIII]{SGA4}), there are canonical isomorphisms
$$R^2f_*\Q_{\ell}\simeq (R^{2d-2}f_*\Ql(d))^{\vee}\quad \text{and} \quad R^2\widetilde f_*\Q_{\ell}\simeq (R^{2d-2}\widetilde f_*\Ql(d))^{\vee},$$
where $d=\Dim(Y_x)=\Dim(\widetilde Y_x)$.
Dualizing and twisting the base change map
$$R^{2d-2}f_*\Ql(d)\rightarrow R^{2d-2}\widetilde f_*\Ql(d),$$
one gets a morphism
$$g_*:R^2\widetilde f_*\Ql(1)\simeq (R^{2d-2}\widetilde f_*\Ql(d))^{\vee}(1)\rightarrow (R^{2d-2}f_*\Ql(d))^{\vee}(1)\simeq R^2f_*\Ql(1).$$
By the compatibility of Poincar\'e duality with base change, for every (not necessarily closed) $x\in X$, the fibre of $g_*$ at $\overline x$ is the usual push forward map $g_{x,*}:H^2(\widetilde Y_{\overline x},\Ql(1))\rightarrow H^2(Y_{\overline x},\Ql(1))$
in étale cohomology. In particular it is compatible with the push forward of algebraic cycles $g_{x,*}:\Pic(\widetilde Y_{\overline x})\otimes \Q\rightarrow \Pic(Y_{\overline x})\otimes \Q$. Since $g^*$ and $g_*$ are maps of sheaves, they are compatible with the specialization isomorphisms and hence the following canonical diagram commutes:
\begin{center}
\begin{tikzcd}
\Pic(Y_{\overline{\eta}})\otimes \Q\arrow{r}{g^*}\arrow{d}{c_{Y_{\overline \eta}}} & \Pic(\widetilde Y_{\overline{\eta}})\otimes \Q\arrow{r}{g_*}\arrow{d}{c_{\widetilde Y_{\overline \eta}}} & \Pic(Y_{\overline{\eta}})\otimes \Q\arrow{d}{c_{Y_{\overline \eta}}}\\
H^2(Y_{\overline{\eta}},\Ql(1))\arrow{r}{g^*}\arrow{d}{\speci_{\eta,x}} & H^2(\widetilde Y_{\overline{\eta}},\Ql(1))\arrow{r}{g_*}\arrow{d}{\widetilde{\speci}_{\eta,x}} &H ^2(Y_{\overline{\eta}},\Ql(1))\arrow{d}{\speci_{\eta,x}}\\
H^2(Y_{\overline x},\Ql(1))\arrow{r}{g^*} & H^2(\widetilde Y_{\overline x},\Ql(1))\arrow{r}{g_*} &H ^2(Y_{\overline x},\Ql(1))\\
\Pic(Y_{\overline x})\otimes \Q\arrow{r}{g^*}\arrow[swap]{u}{c_{Y_{\overline x}}} & \Pic(\widetilde Y_{\overline x})\otimes \Q\arrow{r}{g_*}\arrow[swap]{u}{c_{\widetilde Y_{\overline x}}} & \Pic(Y_{\overline x})\otimes \Q\arrow[swap]{u}{c_{Y_{\overline x}}}.
\end{tikzcd}
\end{center}
Let $\widetilde z_{\eta}$ be a lift of $\widetilde{\varepsilon}_\eta$ to $\Pic(\widetilde Y_{\overline \eta})\otimes \Q$.
Since the composition of the horizontal arrows are the multiplication by $n:=\deg(g)$ maps, one concludes observing that the image $\varepsilon_\eta$ of $\frac{g_*(\widetilde z_\eta)}{n}(\in \Pic(Y_{\overline \eta})\otimes \Q)$ in $\NS(Y_{\overline \eta})\otimes \Q$ is an element such that $\speci_{\eta,x}(\varepsilon_\eta)=\varepsilon_x$. 
\endproof
\part{Applications}\label{secapp}
In this part $k$ is an infinite field of characteristic $p>0$, assumed to be finitely generated except in Subsection \ref{ptate}.
\section{Hyperplane sections}\label{sechyperplane}
In this section we apply Theorem \ref{main} to Lefschetz pencils of hyperplane sections. The main result is Corollary \ref{exinstencehyperplane}.
\numberwithin{equation}{subsection}
 \subsection{Geometric versus arithmetic hyperplane sections}
Let $Z$ be a smooth projective $k$-variety and fix a closed embedding $Z \subseteq \mathbb P_k^n$. One can ask whether there exists a smooth hyperplane section $D$ of $Z$ such that the canonical map 
$$i_{D_{\overline k}}: \NS(Z_{\overline k})\otimes \Q\rightarrow \NS(D_{\overline k})\otimes \Q$$
is an isomorphism. If $\Dim(Z)=2$, then $D$ is a curve so that $\NS(D_{\overline k})\otimes \Q=\Q$ hence $i_{D_{\overline k}}$ is not  injective as soon as $\NS(Z_{\overline k})\otimes \Q$ has rank $\geq 2$. Weak Lefschetz (\cite[Thm. 7.1, p. 318]{milne}) and  Grothendieck–Lefschetz (\cite[Exp. XI]{SGA2}) theorems ensure that $i_{D_{\overline k}}$ is injective if $\Dim(Z)\geq 3$, and an isomorphism if $\Dim(Z)\geq 4$. There are smooth projective varieties of dimension $3$ such that the surjectivity of $i_{D_{\overline k}}$ fails for all smooth hyperplane sections. 
\begin{example}\label{counter}
Take $Z=\mathbb P_k^3$ embedded in $\mathbb P_k^9$ via the Veronese embedding:
\begin{align*}
 \mathbb P_k^3&\rightarrow \mathbb P_k^9 \\
 [x:y:z:w]&\mapsto [x^2:y^2:z^2:w^2:xy:xz:xw:yz:yw:zw].
\end{align*}
Then a smooth hyperplane section $D\subseteq Z$ in $\mathbb P_k^9$ is a smooth quadric surface in $\mathbb P_k^3$, so that $D_{\overline k}\simeq \mathbb P^1_{\overline k}\times \mathbb P^1_{\overline k}$. Hence $\NS(D_{\overline k})\simeq \Z\times \Z$, while $\NS(Z_{\overline k})=\Z$.
\end{example}
But things change if one replaces the geometric N\'eron-Severi groups with their arithmetic counterparts.
\begin{example}\label{example}
Let $Z$ and the embedding $Z\hookrightarrow \mathbb P_k^{9}$ be as in Example \ref{counter}. Assume $p\geq 17$ and consider the pencil of hyperplane sections of $Z$ in $\mathbb P_k^9$ given by the hyperplanes 
$$a(x_1+x_2+x_3+x_4)+b(x_1+4x_2+9x_3+16x_4)=0,$$ where $[a:b]\in \mathbb P^1(k)$ and $x_1,...,x_{10}$ are the homogeneous coordinates in $\mathbb P_k^9$. This corresponds in $\mathbb P^3$ to the pencil of quadric surfaces
$$Q_{[a:b]}:a(x^2+y^2+z^2+w^2)+b(x^2+4y^2+9z^2+16w^2)=0.$$
When $Q_{[a:b]}$ is smooth, it a is quadric surface with disriminant $\Delta_{[a:b]}:=(a+b)(a+4b)(a+9b)(a+16b)$. If $\Delta_{[a:b]}$ is not a square in $k$, then the action of $\pi_1(k)$ on
$$\NS(Q_{[a:b],\overline k})\otimes \Q\simeq \Q\times \Q$$
is non-trivial and it factors through the quotient 
$$\pi_1(k)/\pi_1(k(\sqrt{\Delta_{[a:b]}}))\simeq \Z/2\Z.$$
The generator of $\Z/2\Z$ then acts on $\Q\times \Q$ permuting the two factors so that
$$\NS(Q_{[a:b]})\otimes \Q= (\NS(Q_{[a:b],\overline k})\otimes \Q)^{\pi_1(k)}=(\Q^2)^{\pi_1(k)}=\Q.$$ 
So there are ``lots" of $[a:b]\in \mathbb P^1(k)$ such that the canonical map $$\NS(Z)\otimes \Q\rightarrow \NS(Q_{[a:b]})\otimes \Q$$
is an isomorphism. 
\end{example}
The main result is of this subsection is Corollary \ref{exinstencehyperplane} that we now recall:
\begin{customthm}{\ref{exinstencehyperplane}}
\textit{If} $\Dim(Z)\geq 3$ \textit{there are infinitely many} $k$\textit{-rational hyperplane sections} $D$ \textit{such that
the canonical map}
$$\NS(Z)\otimes \Q\rightarrow \NS(D)\otimes \Q$$ \textit{is an isomorphism.}
\end{customthm}
\numberwithin{equation}{subsection} 
\subsection{Proof of Corollary \texorpdfstring{\ref{exinstencehyperplane}}- }\label{ammazzachemazza}
By (\cite[Exp. XVII]{SGA7}) there exists a pencil of hyperplanes $L:=\{H_x\}_{x\in \check{\mathbb P}_k^1}$ such that:
\begin{itemize}
\item For all $x$ in an dense open subscheme $U\subseteq \check{\mathbb P}_k^1$, the intersection $H_x\cap Z$ is smooth;
\item The base locus $B:=\bigcap_{x\in \check{\mathbb P}_k^1} Z\cap H_{x}\subseteq Z$ is smooth.
\end{itemize}
Then one gets a diagram
\begin{center}
\begin{tikzcd}
Z & \arrow[swap]{l}{\pi}\arrow{r}{f}\widetilde Z & \check{\mathbb P}_k^1,
\end{tikzcd}
\end{center}
where $\pi:\widetilde Z\rightarrow Z$ is the blow up of $Z$ along $B$, $f:\widetilde Z\rightarrow \check{\mathbb P}_k^1$ is a projective flat morphism smooth over $U$ and for each $x\in \check{\mathbb P}_k^1$ the fibre $\widetilde Z_x$ of $f:\widetilde Z\rightarrow \check{\mathbb P}_k^1$ at $x$ identifies via $\pi:\widetilde Z\rightarrow Z$ with the hyperplane section $Z\cap H_x\subseteq Z$. Write $E:=\pi^{-1}(B)$ for the exceptional divisor.
Explicitly $\widetilde Z$ is the closed subscheme of $Z\times \check{\mathbb P}_k^1$ defined by
$$\widetilde Z:=\{(z,x)\in Z\times \check{\mathbb P}^1_k \text{ with } z\in Z\cap H_x\}\hookrightarrow Z\times \check{\mathbb P}^1_k,$$
$\pi:\widetilde Z\rightarrow Z$, $f:\widetilde Z\rightarrow \check{\mathbb P}^1_k$ are identified with the canonical projections onto $Z$, $\check{\mathbb P}^1_k$ respectively and $E$ with $B\times \check{\mathbb P}_k^1$. Write $\eta$ for the generic point of $\check{\mathbb P}^1_k$.

Since $k$ is infinite finitely generated (hence Hilbertian), if $S\subseteq |\mathbb P^1|$ is a sparse set, then $\mathbb P^1(k)-S$ is infinite (see for example \cite[Proposition 8.5(d)]{poonen} and \cite[9.5, Remarks, 2)]{Serreweil}). Hence, combining Lemma \ref{hyperplane} below with Corollary \ref{existence} one gets Corollary \ref{exinstencehyperplane}.
\begin{lemma}\label{hyperplane}
If $\Dim(Z)\geq 3$, the canonical map 
$$i^*_{\widetilde Z_{\eta}}:\NS(Z)\otimes \Q\rightarrow \NS(\widetilde Z_{\eta})\otimes \Q$$ is an isomorphism.
\end{lemma}
\proof
This is inspired from \cite[Corollary 1.5]{Varmor}. 
Fix $x\in U$. The natural commutative diagram 
\begin{center}
\begin{tikzcd}
& \widetilde Z_{\eta}\arrow{r}\arrow{d}\arrow[swap]{dl}{i_{\widetilde Z_{\eta}}}\arrow[phantom]{rd}{\Box}\arrow[near start, bend left]{dd}{i_{\eta}} & k(\eta)\arrow{d}{\eta}\arrow[bend left]{dd}\\
Z& \widetilde Z\arrow[swap]{l}{\pi}\arrow{r}{f}\arrow[phantom]{rd} & \check{\mathbb P}_k^1\\
&  \widetilde Z_U\arrow{u}{i_U}\arrow{r}{f_U}\arrow[phantom]{rd}{\Box} & U\arrow{u}\\
& \arrow{uul}{i_{\widetilde Z_x}}\arrow{u}{i_x}\widetilde Z_{x}=Z\cap H_x\arrow{r} &k(x)\arrow{u}{x}
\end{tikzcd}
\end{center}
induces a commutative diagram
\begin{center}
\begin{tikzcd}
\Pic(Z)\otimes \Q\arrow{r}{\pi^*}\arrow[two heads]{d}{q_Z} & \Pic(\widetilde Z)\otimes \Q\arrow[two heads]{r}{i^*_U}\arrow[two heads]{d}{q_{\widetilde Z}} & \Pic(\widetilde Z_{U})\otimes \Q\arrow[two heads, bend left=15]{rr}{i^*_{\eta}}\arrow{r}{i^*_x} & \Pic(\widetilde Z_x)\otimes \Q\arrow[two heads]{d}{q_{\widetilde Z_x}}&\Pic(\widetilde Z_{\eta})\otimes \Q\arrow[two heads]{d}{q_{\widetilde Z_\eta}} \\
\NS(Z)\otimes \Q\arrow{r}{\pi^*}\arrow[bend right=10,swap]{rrr}{i^*_{\widetilde Z_x}} & \NS(\widetilde Z)\otimes \Q\arrow{rr} & & \NS(\widetilde Z_{x})\otimes \Q & \NS(\widetilde Z_{\eta})\otimes \Q \arrow[hook]{l}{\speci^{\mathrm{ar}}_{\eta,x}},
\end{tikzcd}
\end{center}
where $q_Z$, $q_{\widetilde Z}$, $q_{\widetilde Z_x}$ and $q_{\widetilde Z_\eta}$ denote the natural surjections.

Since $\Dim(Z)\geq 3$, by the weak Lefschetz theorem (\cite[Thm. 7.1, p. 318]{milne}), the natural map 
$$i^*_{\widetilde Z_{x}}:H^2(Z_{\overline k},\Ql(1))\rightarrow H^2(\widetilde Z_{\overline x},\Ql(1))$$
is injective. Since there is a commutative diagram
\begin{center}
\begin{tikzcd}
 \NS(Z)\otimes \Q\arrow[hook]{d}\arrow{r}  &  \NS(\widetilde Z_x)\otimes \Q\arrow[hook]{d}\\
H^2(Z_{\overline k},\Ql(1))\arrow[hook]{r} & H^2(\widetilde Z_{\overline x},\Ql(1))
\end{tikzcd}
\end{center}
with injective vertical arrows, this implies that  the map
$$i^*_{\widetilde Z_{x}}:\NS(Z)\otimes \Q\rightarrow  \NS(\widetilde Z_{x})\otimes \Q$$
is injective. Since $i^*_{\widetilde Z_{x}}$ factorizes as
$$\NS(Z)\otimes \Q\xrightarrow{i^*_{\widetilde Z_{\eta}}}  \NS(\widetilde Z_{\eta})\otimes \Q\xrightarrow{\speci^{\ar}_{\eta,x}}  \NS(\widetilde Z_{x})\otimes \Q,$$
also 
$$i^*_{\widetilde Z_{\eta}}:\NS(Z)\otimes \Q\rightarrow  \NS(\widetilde Z_{\eta})\otimes \Q$$
is injective.

To prove the surjectivity, let $\varepsilon\in \NS(\widetilde Z_{\eta})\otimes \Q$ with lift $z\in \Pic(\widetilde Z_{\eta})\otimes \Q$. Since $\speci^{\ar}_{\eta,x}$ is injective, it is enough to show that $\varepsilon_x:=\speci^{\ar}_{\eta,x}(\varepsilon)$ is in the image of $i^*_{\widetilde Z_x}$. 
Since the maps $i^*_U$ and $i^*_{\eta}$ are surjective, $z\in \Pic(\widetilde Z_{\eta})$ lifts to a $\widetilde z\in \Pic(\widetilde Z)\otimes \Q$ and, by the commutativity of the diagram, $\widetilde z$ maps to $\varepsilon_x$ in $\NS(\widetilde Z_{x})\otimes \Q$. Now, since $\pi$ is the blow up of $Z$ along $B$, $q_{\widetilde Z}(\widetilde z)$ can be written as $\pi^*q_{Z}(z')+bq_{\widetilde Z}(E)$, where $z'\in \Pic(Z)\otimes \Q$ and $b\in \Q$. The conclusion follows from the following claim, since it implies that $\varepsilon_x$ is the image of $q_Z(z')+bq_Z(Z\cap H_x)\in \NS(Z)\otimes \Q$.\\
\textcolor{white}{a}\\
\textbf{Claim:} \textit{The restrictions of} $q_{\widetilde Z}(E)$ \textit{and} $\pi^*q_{Z}(Z\cap H_x)$ \textit{to} $\widetilde Z_x$ \textit{coincide.}\proof[Proof of the claim]
By direct computations, one sees that $E=B\times \check{\mathbb P}_k^1$ intersects transversally with $\widetilde Z_{x}$ and that $E\cap \widetilde Z_x=B$, so that the restriction of $q_{\widetilde Z}(E)$ to $\widetilde Z_{x}$ is given by $q_{\widetilde Z_x}(B)\in \NS(\widetilde Z_x)\otimes \Q$. To compute the restriction $\pi^*q_{Z}(Z\cap H_x)$ to $\widetilde Z_x$ observe that it is equal to $i_{\widetilde Z_x}^*(q_{Z}(Z\cap H_x))$. Then take any $y\neq x\in L$ and compute $i_{\widetilde Z_x}^*(q_Z(Z\cap H_x))$ as $q_{\widetilde Z_x}(Z\cap H_x\cap H_y)=q_{\widetilde Z_x}(B)$, since $B=Z\cap H_x\cap H_y$ for any $x\neq y\in L$.
\endproof
\begin{remark}
The key fact that $\Pic(\widetilde Z)\otimes \Q\rightarrow \Pic(\widetilde Z_{\eta})\otimes \Q$ is surjective, does not hold for $\Pic(\widetilde Z_{\overline k})\otimes \Q\rightarrow \Pic(\widetilde Z_{\overline \eta})\otimes \Q$ (see Examples \ref{counter} and \ref{example}).
This is why it is not true in general that for a point $x\in |U|$, which is Galois-generic for $\widetilde Z_{U}\rightarrow U$ the canonical map
$$i^*_{\widetilde Z_x}:\NS(Z)\otimes \Q\rightarrow \NS(\widetilde Z_x)\otimes \Q$$ is an isomorphism and one really needs to restrict to strictly Galois-generic points: during the proof one cannot replace $U$ with a finite étale cover, since any base change destroys the geometry of the pencil. 
\end{remark}
\subsection{Proof of Corollary \texorpdfstring{\ref{corollaryTate}}-}
Replacing $k$ with a finite extension, it is enough to show that the map
$$\NS(Z)\otimes \Ql\rightarrow H^2(Z_{\overline k},\Ql(1))^{\pi_1(k)}$$ is an isomorphism. 
By Corollary \ref{exinstencehyperplane}, there exists $k$-rational hyperplane section $D\rightarrow Z$ such that the canonical map $\NS(Z)\otimes \Ql\xrightarrow{\sim} \NS(D)\otimes \Ql$ is an isomorphism. The conclusion follows from the commutative diagram
\begin{center}
\begin{tikzcd}
\NS(Z)\otimes \Ql\arrow[hook]{r}\arrow{d}{i_{D}^*}& H^2(Z_{\overline k},\Ql(1))^{\pi_1(k)}\arrow[hook]{d}{i_{D}^*}\\
\NS(D)\otimes \Ql\arrow[hook]{r}{(2)}& H^2(D_{\overline k},\Ql(1))^{\pi_1(k)}
\end{tikzcd}
\end{center}
since $i_{D}^*$ is an isomorphism by the choice of $D$ and $(2)$ is an isomorphism by $T(D,\ell)$.
\endproof
\section{Brauer groups in families}\label{secbrauer}
\numberwithin{equation}{subsection} 
\subsection{Specialization of Brauer groups}\label{brauergen}
\subsubsection{Brauer group}
For a smooth proper $k$-variety $Z$ write $H^2(Z_{\overline k},\mathbb G_m):=\Br(Z_{\overline k})$ for the (cohomological) Brauer group of $Z_{\overline k}$, $\Br(Z_{\overline k})[n]$ for its $n$-torsion subgroup and
$$T_{\ell}(\Br(Z_{\overline k})):=\varprojlim_{n}\Br(Z_{\overline k})[\ell^n],\quad \Br(Z_{\overline k})[\ell^{\infty}]:=\varinjlim_{n}\Br(Z_{\overline k})[\ell^n],\quad \Br(Z_{\overline k})[p']:=\varinjlim_{n\nmid p}\Br(Z_{\overline k})[n].$$

Recall that $\Br(Z_{\overline k})$ is a torsion group and that Kummer theory induces, for every $p\nmid n\in \N$, an exact sequence:
$$0\rightarrow \NS(Z_{\overline k})/n\rightarrow H^2(Z_{\overline k},\mu_{n})\rightarrow \Br(Z_{\overline k})[n]\rightarrow 0.$$
It is classically known that if $T(Z,\ell)$ holds, then $\Br(Z_{\overline k})[\ell^{\infty}]^{\pi_1(k)}$ is finite (see e.g. \cite[Proposition 2.1.1]{brauer}). 
\subsubsection{Brauer generic points}\label{brauergeneric}
Let $X$ be a smooth geometrically connected $k$-variety with generic point $\eta$ and $f:Y\rightarrow X$ a smooth proper morphism of $k$-varieties. Taking the direct limit over $p\nmid n$ on the Kummer exact sequence, one gets a commutative specialization exact diagram
\begin{center}
\begin{tikzcd}
0\arrow{r} & \varinjlim_{n\nmid p} \NS(Y_{\overline \eta})/n\arrow[hook]{d}{\speci_{\eta,x}}\arrow{r} & \varinjlim_{n\nmid p} H^2(Y_{\overline \eta},\mu_n)\arrow{r}\arrow{d}{\simeq}& \Br(Y_{\overline \eta})[p']\arrow[two heads]{d}{{\speci^{\Br}_{\eta,x}}}\arrow{r} & 0\\
0\arrow{r} & \varinjlim_{n\nmid p} \NS(Y_{\overline x})/n\arrow{r} & \varinjlim_{n\nmid p} H^2(Y_{\overline x},\mu_n)\arrow{r}& \Br(Y_{\overline x})[p']\arrow{r} & 0
\end{tikzcd}
\end{center}
Since the group $\Ker(\speci_{\eta,x})$ is of $p$-torsion and $\coker(\speci_{\eta,x})_{tors}=\coker(\speci_{\eta,x})[p^{\infty}]$ (see \cite[Proposition 3.6]{poonen}), one sees that a $x\in |X|$ is $\NS$-generic if and only if the map 
$$\speci^{\Br}_{\eta,x}:\Br(Y_{\overline \eta})[p']\rightarrow \Br(Y_{\overline x})[p']$$
is an isomorphism. In particular (Corollary \ref{existence}) the set of $x\in |X|$ such that $\speci^{\Br}_{\eta,x}$ is not an isomorphism is sparse and if $X$ is a curve and $f:Y\rightarrow X$ has the LU-property in degree 2, it contains at most finitely many $k$-rational points. 
\subsection{Uniform boundedness}
Retain the notation and the assumptions as in the previous Section \ref{brauergeneric}. Assume furthermore that $\ell\neq p$.
Taking inverse limit in the Kummer exact sequence, one gets a commutative exact diagram
\begin{center}
\begin{tikzcd}
0\arrow{r} & \NS(Y_{\overline \eta})\otimes \Z_{\ell} \arrow{r}\arrow[hook]{d}{\speci_{\eta,x}} & H^2(Y_{\overline \eta},\Zl(1))\arrow{r}\arrow{d}{\simeq} &T_{\ell}(\Br(Y_{\overline \eta})) \arrow{r}\arrow[two heads]{d}{\speci^{\Br}_{\eta,x}} & 0\\
0\arrow{r} & \NS(Y_{\overline x})\otimes \Z_{\ell} \arrow{r} & H^2(Y_{\overline x},\Zl(1))\arrow{r} &T_{\ell}(\Br(Y_{\overline x})) \arrow{r} & 0.
\end{tikzcd}
\end{center}
The group $\pi_1(x,\overline x)$ acts on $T_{\ell}(\Br(Y_{\overline \eta}))$ by restriction through the map $\pi_1(x,\overline x)\rightarrow \pi_1(X,\overline x)\simeq \pi_1(X,\overline \eta)$ and $\speci^{\Br}_{\eta,x}$ is $\pi_1(x,\overline x)$-equivariant with respect to the natural action of $\pi_1(x,\overline x)$ on $T_{\ell}(\Br(Y_{\overline x}))$. Hence the arguments in Section \ref{brauergen} combined with Theorem \ref{main} show the following
\begin{lemma}\label{lemmaidentification}
Up to replacing $X$ with an open subset, for every Galois-generic $x\in |X|$ and every $\ell\neq p$, the $\pi_1(x,\overline x)$-equivariant specialization morphism
$$\speci^{\Br}_{\eta,x}:T_{\ell}(\Br(Y_{\overline \eta}))\rightarrow T_{\ell}(\Br(Y_{\overline x}))$$
is an isomorphism.
\end{lemma}
Replacing \cite[Proposition 3.2.1]{brauer} with Lemma \ref{lemmaidentification} and \cite[Fact 3.4.1]{brauer} with the main result of \cite{tamagawa}, one can make the arguments in the proof of \cite[Theorem 1.2.1]{brauer} work in positive characteristic and prove Corollary \ref{corollarybrauer}. 
In the same way, using the arguments in the proof of \cite[Theorem 3.5.1]{brauer}, one gets the following unconditional variant:
\begin{corollary}\label{uncondintionalbrauer}
Let $X$ be a curve and assume that the Zariski closure of the image of $\pi_1(X)$ acting on $H^2(Y_{\overline \eta},\Ql(1))$ is connected. If $Y\rightarrow X$ has the $LU$-property in degree 2, then there exists an integer $C:=C(Y\rightarrow X,\ell)$ 
such that $$[\Br(Y_{\overline{x}})^{\pi_1(x,\overline x)}[\ell^\infty]:\Br(Y_{\overline{\eta}})^{\pi_1(X,\overline \eta)}[\ell^\infty]]\leq C $$ 
for all but finitely many $x\in X(k)$.
\end{corollary}
\numberwithin{equation}{subsection} 
\subsection{\texorpdfstring{$p$}--adic Tate module}\label{ptate}
Assume that $X$ is a smooth connected $k$-variety with generic point $\eta$, where $k$ is an algebraically closed field of characteristic $p$ and that $Y\rightarrow X$ is a smooth projective morphism. 
\begin{corollary}\label{corptate}
There exists an $x\in |X|$ such that 
$\Rango(T_p(\Br(Y_{\overline x})))=\Rango(T_p(\Br(Y_{\overline \eta})))$
\end{corollary}
\proof
For every geometric point $\overline t\in X$, one has (\cite[Proposition 5.12]{derhamwitt}):
$$\Dim(\NS(Y_{\overline t})\otimes \Q_p)=\Dim(H^2_{\cristalline}(Y_{\overline t}))-2\Dim(H^2_{\cristalline}(Y_{\overline t}))_{[1]})-\Rango(T_p(\Br(Y_{\overline t})))$$
where $H^2_{\cristalline}(Y_{\overline t})_{[1]}$ is the slope one part of the crystalline cohomology of $Y_{\overline t}$ (see e.g. \cite[Section 3]{notesked} for the definition). By \cite[Theorem 3.12, Corollary 4.2]{notesked} there exists a dense open subset $U$ of $X$ such that for all $x\in |U|$ one has
$$\Dim(H^2_{\cristalline}(Y_{\overline x})_{[1]})=\Dim(H^2_{\cristalline}(Y_{\overline{\eta}})_{[1]})$$
Since $\Dim(H^2_{\cristalline}(Y_{\overline x}))$ is independent of $x\in X$ (smooth proper base change in crystalline cohomology), one concludes applying Corollary \ref{existence} to $Y_U\rightarrow U$.
\endproof
\newpage
\part{Overconvergence of \texorpdfstring{$R^if_{\Ogus,*}\mathcal O_{Y/K}$}- (after A. Shiho)}\label{appendix}
In this part we use the work of Shiho on relative rigid cohomology to prove Theorem \ref{berthelotconjecture}, which is a key ingredient in the proof of Theorem \ref{Shiho}. In Section \ref{secprel}, we recall the definitions of various categories of isocrystals, the relations between them and we state Theorem \ref{berthelotconjecture}. In Section \ref{secproof}, we prove Theorem \ref{berthelotconjecture}.
\section{Preliminaries}\label{secprel}
\numberwithin{equation}{subsubsection} 
\subsection{Notation}
Let $k$ be a perfect field of characteristic $p>0$. Write $K$ for the fraction field of the Witt ring $W:=W(k)$ of $k$ and $|-|:K\rightarrow \mathbb R$ for a norm induced by the ideal $pW\subseteq W$. For any $k$-variety $X$, write $F_X$ for a power of the absolute Frobenius on $X$ and, if there is no danger of confusion, one often drops the lower index and writes just $F$.

Gothic letters ($\mathfrak T,\mathfrak X,\mathfrak U, ...$) denote separated $p$-adic formal schemes topologically of finite type over $W$.
Write $\mathfrak X_1$ for the special fibre of $\mathfrak X$, $\mathfrak X_K$ for its rigid analytic generic fibre and $\speci:\mathfrak X_{K}\rightarrow \mathfrak X_1$ for the specialization map. There is an equivalence between the isogeny category $\mathbf{Coh}(\mathfrak X)\otimes \Q$ of the category $\mathbf{Coh}(\mathfrak X)$ of coherent sheaves on $\mathfrak X$ and the category $\mathbf{Coh}(\mathfrak X_K)$ of coherent sheaves on $\mathfrak X_K$ (see \cite[Remark 1.5]{Ogus}).

If $f:X\rightarrow \mathfrak X_1$ is a closed immersion, one can consider the open tube $]X[_{\mathfrak X}:=\speci^{-1}(X)$ and the closed tube of radius $|p|$,  $[X]_{\mathfrak X,|p|}$ of $X$ in $\mathfrak X$ (see \cite[Definition 1.1.2, Section 1.1.8]{berthelot}). They are admissible open subsets of $\mathfrak X_K$ and there is an inclusion $[X]_{\mathfrak X,|p|}\subseteq ]X[_{\mathfrak X}$.

A pair $(X,\overline X)$ is an open immersion of $k$-varieties $X\rightarrow \overline X$ and a frame $(X,\overline X,\mathfrak X)$ is a pair $(X,\overline X)$ together with a closed immersion of $\overline X$ into a $p$-adic formal scheme $\mathfrak X$. Morphisms of pairs and frames are defined in the obvious way. A pair $(Y,\overline Y)$ over a frame $(X,\overline X,\mathfrak X)$ is a morphism of pairs $(Y,\overline Y)\rightarrow (X,\overline X)$ and a frame $(X,\overline X,\mathfrak X)$ over a pair $(Y,\overline Y)$ is a morphism of pairs $(X,\overline X)\rightarrow (Y,\overline Y)$. If $(X,\overline X,\mathfrak X)$ is a frame, for any sheaf $\mathcal F$ over $]\overline X[_{\mathfrak X}$ one writes $$j^{\dagger}_{X}\calF:=\varinjlim_{V}j_{V*}j_{V}^{*}\mathcal F$$
where the limit runs over all the strict neighbourhoods $V$ of $X$ in $\overline X$ (see \cite[Definition 1.2.1]{berthelot}) and $j_V:V\rightarrow ]\overline X[_{\mathfrak X}$ is the inclusion map.

If $f:Y\rightarrow X$ is a morphism of $k$-varieties, for every morphism $Z\rightarrow X$ write:
\begin{center}
\begin{tikzcd}
Y_{Z}\arrow{r}\arrow{d}{f_Z}\arrow[phantom]{dr}{\Box}& Y\arrow{d}{f}\\
Z\arrow{r}& X
\end{tikzcd}
\end{center}
\numberwithin{equation}{subsection} 
\subsection{Categories of isocrystals}
To a $k$-variety $X$ one can associate the following categories of isocrystals:
\begin{itemize}
\item $\Isoc^{(p)}(X)$, the $p$-adically convergent isocrystals (see \cite[Definition 2.1]{Ogus});
\item $\Isoc^{(1)}(X)$, the convergent isocrystals (see \cite[Definition 2.1]{Ogus});
\end{itemize}
If $(X,\overline X)$ is a pair there is a category $\Isoc^{\dagger}(X,\overline X)$ of isocrystals on $X$ overconvergent along $\overline X-X$ see \cite[Definition 2.3.2]{berthelot}. If $\overline X$ is a compactification of $X$, one writes $\Isoc^{\dagger}(X,\overline X):=\Isoc^{\dagger}(X)$ and calls the object there overconvergent isocrystals on $X$. It is known that $\Isoc^{\dagger}(X)$ does not depend on the choice of the compactification, so that $\Isoc^{\dagger}(X)$ is well defined (see \cite[2.3.6]{berthelot}).
\begin{itemize}
\item $\Isoc^{(p)}(X)$ (resp. $\Isoc^{(1)}(X)$). Write $I^{(p)}_X$ (resp. $I^{(1)}_X$) for the category of $p$-adic enlargements (resp. enlargements). This is the category of pairs $(\mathfrak T,z_{\mathfrak T})$ such that $\mathfrak T$ is a flat $p$-adic formal $W$-scheme and $z_{\mathfrak T}$ is a morphism $\mathfrak T_1\rightarrow X$ (resp. $(\mathfrak T_1)_{\red}\rightarrow X$). A morphism $g:(\mathfrak Z,z_{\mathfrak Z})\rightarrow (\mathfrak T,z_{\mathfrak T})$ between $p$-adic enlargements (resp. enlargements) is a morphism $g:\mathfrak Z\rightarrow \mathfrak T$ such that $z_{\mathfrak Z}\circ g_1= z_{\mathfrak T}$ (resp. $z_{\mathfrak Z}\circ (g_1)_{\red}=z_{\mathfrak T}$), where $g_1:\mathfrak Z_1\rightarrow \mathfrak T_1$ (resp. $(g_1)_{\red}:(\mathfrak Z_1)_{\red}\rightarrow (\mathfrak T_1)_{\red}$) are the natural morphisms induced by $g$.
A $p$-adically convergent isocrystals (resp. a convergent isocrystal) is the following set of data:
\begin{itemize}
\item For every $(\mathfrak T,z_{\mathfrak T})\in \mathbf{Ob}(I^{(p)}_X)$ (resp. $\in \mathbf{Ob}(I^{(1)}_X)$), a $\mathcal M_{(\mathfrak T,z_{\mathfrak T})}\in \mathbf{Coh}(\mathfrak T_K)$;
\item For every morphism $g:(\mathfrak Z,z_{\mathfrak Z})\rightarrow (\mathfrak T,z_{\mathfrak T})$ in $I^{(p)}_X$ (resp. $I^{(1)}_X$) an isomorphism
$$\phi_g:g^*\mathcal M_{(\mathfrak T,z_{\mathfrak T})}\rightarrow \mathcal M_{(\mathfrak Z,z_{\mathfrak Z})}$$
in $\mathbf{Coh}(\mathfrak Z_K)$ such that $\phi_{\mathrm{Id}}=\mathrm{Id}$ and for every other morphism $h:(\mathfrak T,z_{\mathfrak T})\rightarrow (\mathfrak U,z_{\mathfrak U})$ one has $\phi_g\circ g^*(\phi_h)=\phi_{h\circ g}$.
\end{itemize}
A morphism of $p$-adically convergent isocrystals (resp. convergent isocrystals) $\mathcal M\rightarrow \mathcal N$ is a collection of morphisms $\{\mathcal M_{(\mathfrak T,z_{\mathfrak T})}\rightarrow \mathcal N_{(\mathfrak T,z_{\mathfrak T})}\}_{(\mathfrak T,z_{\mathfrak T})\in \mathbf{Ob}(I^{(p)}_X)}$ (resp. $_{(\mathfrak T,z_{\mathfrak T})\in \mathbf{Ob}(I^{(1)}_X)}$) compatible with the isomorphism $\phi_g$ for all morphisms $g$.
\item $\Isoc^{\dagger}(X,\overline X)$. Write $I_{(X,\overline X)}$ for the category of frames over $(X,\overline X)$.
Then an isocrystals on $X$ overconvergent along $\overline X-X$ is the following set of data:
\begin{itemize}
\item For every $(T,\overline T,\mathfrak T)\in \mathbf{Ob}(I_{(X,\overline X)})$ a coherent $j^{\dagger}_T\mathcal O_{]\overline T[_{\mathfrak T}}$-module $\mathcal M_{(T,\overline T,\mathfrak T)}$;
\item For every morphism $g:(Z,\overline Z,\mathfrak Z)\rightarrow (T,\overline T,\mathfrak T)$ in $I_{(X,\overline X)}$ an isomorphism
$$\phi_g:g^*\mathcal M_{(T,\overline T,\mathfrak T)}\rightarrow \mathcal M_{(Z,\overline Z,\mathfrak Z)}$$
of coherent $j^{\dagger}_{Z}\mathcal O_{]\overline Z[_{\mathfrak Z}}$-modules such that $\phi_{\mathrm{Id}}=\mathrm{Id}$ and for every other morphism $h:(T,\overline T,\mathfrak T)\rightarrow (U,\overline U,\mathfrak U)$ one has $\phi_g\circ g^*(\phi_h)=\phi_{h\circ g}$.
\end{itemize}
A morphism $\mathcal M\rightarrow \mathcal N$ in $\Isoc^{\dagger}(X,\overline X)$ is a collection of morphisms $\{\mathcal M_{(T,\overline T,\mathfrak T)}\rightarrow \mathcal N_{(T,\overline T,\mathfrak T)}\}_{(T,\overline T,\mathfrak T)\in \mathbf{Ob}(I_{(X,\overline X)})}$ compatible with the isomorphism $\phi_g$ for all morphisms $g$.
\end{itemize}
There are also enriched versions of the previous categories with Frobenius structure, which we denote $\Fisoc^{(p)}(X)$, $\Fisoc^{(1)}(X)$ and $\Fisoc^{\dagger}(X,\overline X)$. For example, the absolute Frobenius $F_X$ induces an endofunctor $$F_X^*:\Isoc^{(p)}(X)\rightarrow \Isoc^{(p)}(X)$$ and $\Fisoc^{(p)}(X)$ is the category of pairs $(\mathcal M,\Phi)$, where $\mathcal M\in \Isoc^{(p)}(X)$ and $\Phi$ is a Frobenius structure  on $\mathcal M$, i.e. an isomorphism $F_X^*\mathcal M\rightarrow \mathcal M$. A morphism in $\Fisoc^{(p)}(X)$ is a morphism in $\Isoc^{(p)}(X)$ compatible with the Frobenius structures. The constructions of $\Fisoc^{(1)}(X)$ and $\Fisoc^{\dagger}(X,\overline X)$ are similar.
\numberwithin{equation}{subsection} 
\subsection{Functors between the categories}\label{functors}
For every pair $(X,\overline X)$ there is a canonical commutative diagram of functors:
\begin{center}
\begin{tikzcd}
\Fisoc^{(p)}(X)\arrow{d}& \arrow{l}{\mathrm{F1}\mbox{-}\mathrm{Fp}} \Fisoc^{(1)}(X)\arrow{d} & \Fisoc^{\dagger}(X,X)\arrow{l}{\mathrm{Fconv}\mbox{-}\mathrm{F1}}\arrow{d}& \Fisoc^{\dagger}(X,\overline X)\arrow{d} \arrow{l}{\mathrm{Fov}\mbox{-}\mathrm{Fconv}}\\
\Isoc^{(p)}(X)& \arrow{l}{\mathrm{1}\mbox{-}\mathrm{p}} \Isoc^{(1)}(X)& \Isoc^{\dagger}(X,X)\arrow{l}{\mathrm{conv}\mbox{-}\mathrm{1}}& \Isoc^{\dagger}(X,\overline X) \arrow{l}{\mathrm{ov}\mbox{-}\mathrm{conv}}
\end{tikzcd}
\end{center}
and 
\begin{itemize}
\item $\mathrm{F1}\mbox{-}\mathrm{Fp}, \mathrm{conv}\mbox{-}\mathrm{1}, \mathrm{Fconv}\mbox{-}\mathrm{F1}$ are equivalences of categories (\cite[Proposition 2.18]{Ogus}, \cite[2.3.4]{berthelot});
\item $\mathrm{Fov}\mbox{-}\mathrm{Fconv}$ is fully faithful if $X$ is smooth (\cite[Theorem 1.1]{Kedfull}).
\end{itemize}
All the functors are easy to construct from the definitions. For example, to construct $\mathrm{conv}\mbox{-}\mathrm{1}$, to an enlargement $(\mathfrak T,z_{\mathfrak T})$ one associates the frame
$((\mathfrak T_1)_{\red},(\mathfrak T_1)_{\red},\mathfrak T)$ over $(X,X)$ and so, for every $\mathcal M\in \Foi(X,X)$, one defines
$$\mathrm{conv}\mbox{-}\mathrm{1}(\mathcal M)_{(\mathfrak T,z_{\mathfrak T})}:=\mathcal M_{((\mathfrak T_1)_{\red},(\mathfrak T_1)_{\red},\mathfrak T)}.$$
The constructions of $\mathrm{1}\mbox{-}\mathrm{p},\mathrm{ov}\mbox{-}\mathrm{conv}$ are similar. In view of these functors, if $\mathcal M$ is in $\Fisoc^{\dagger}(X,\overline X)$ and $(\mathfrak T,z_{\mathfrak T})$ is a $p$-adic enlargement of $X$, write 
$$\mathcal M_{(\mathfrak T,z_{\mathfrak T})}:=\mathcal M_{(\mathfrak T_1,\mathfrak T_1,\mathfrak T)}.$$
\numberwithin{equation}{subsection}
\subsection{Stratification}
Assume that $X$ admits a closed immersion into a $p$-adic formal scheme $\mathfrak X$ formally smooth over $W$. Then the categories of isocrystals on $X$ admit a more concrete description in term of modules with a stratification. We now recall the notion of universal $p$-adic enlargement and we use it to define modules with a stratification.
\subsubsection{Universal \texorpdfstring{$p$}--adic enlargements}
By \cite[Proposition 2.3]{Ogus}, there exists a universal $p$-adic enlargement $(\mathfrak T(X), z_{\mathfrak T(X)})$ of $X$ in $\mathfrak X$. The $p$-adic enlargement $(\mathfrak T(X), z_{\mathfrak T(X)})$ of $X$ is endowed with a map $g:\mathfrak T(X)\rightarrow \mathfrak X$ making the following diagram commutative:
\begin{center}
\begin{tikzcd}
& \mathfrak T(X)_1\arrow[hook]{r}\arrow[swap]{dl}{z_{\mathfrak T(X)}}\arrow{d}{g_1} &\mathfrak T(X)\arrow{d}{g}\\
X\arrow[hook]{r} & \mathfrak X_1\arrow{r}& \mathfrak X
\end{tikzcd}
\end{center}
which is universal for all the $p$-adic enlargements $(\mathfrak Y, z_{\mathfrak Y})$ of $X$ in $\mathfrak X$, i.e for all the $p$-adic enlargements $(\mathfrak Y,z_{\mathfrak Y})$ admitting a map $g:\mathfrak Y\rightarrow \mathfrak X$ making the previous diagram commutative.

Write $\mathfrak T(X)(1)$ for the universal $p$-adic enlargements of $X$ in $\mathfrak X\times \mathfrak X$, where one considers $X$ embedded in $\mathfrak X_1\times \mathfrak X_1$ via the diagonal immersion. 
The $p$-adic formal schemes $\mathfrak T(X)$ and $\mathfrak T(X)(1)$ are such that $\mathfrak T(X)_K=[X]_{\mathfrak X,|p|}$ and $\mathfrak T(X)(1)_K=[X]_{\mathfrak X\times \mathfrak X,|p|}$ (see \cite[1.1.10]{berthelot}).
\subsubsection{Stratifications}
Let  $\Strat^{(p)}(X,\mathfrak X)$ be the category of modules with a $p$-adically convergent stratifications (\cite[Proposition 2.11]{Ogus}). An object $(\mathcal M,\varepsilon)$ in $\Strat^{(p)}(X,\mathfrak X)$ is a coherent module $\mathcal M$ over $[X]_{\mathfrak X,|p|}=\mathfrak T(X)_K$ together with an isomorphism $\varepsilon:p_1^*\calM\rightarrow p_2^*\calM$ satisfying a natural cocycle condition, where the $p_i's$ are the two projections $\mathfrak T(X)(1)_K=[X]_{\mathfrak X\times \mathfrak X,|p|}\rightarrow \mathfrak T(X)_K=[X]_{\mathfrak X,|p|}$.
The projections $p_1,p_2:\mathfrak X\times \mathfrak X\rightarrow  \mathfrak X$
give morphisms of enlargements $p_1,p_2:(\mathfrak T(X)(1),z_{\mathfrak T(X)(1)})\rightarrow (\mathfrak T(X),z_{\mathfrak T(X)})$.
Hence, if $\mathcal M$ is in $\Isoc^{(p)}(X)$, there is an isomorphism
$$\varepsilon_{\mathcal M,\mathfrak X}:p_1^*(\mathcal M_{(\mathfrak T(X),z_{\mathfrak T(X)})})\simeq \mathcal M_{(\mathfrak T(X)(1),z_{\mathfrak T(X)(1)})} \simeq p_2^*(\mathcal M_{(\mathfrak T(X),z_{\mathfrak T(X)})}). $$
This gives a functor
$$(-_{(\mathfrak T(X),z_{\mathfrak T(X)})},\varepsilon_{-,\mathfrak X}):\Isoc^{(p)}(X)\rightarrow \Strat^{(p)}(X,\mathfrak X)$$ that sends $\mathcal M$ to $(\mathcal M_{(\mathfrak T(X),z_{\mathfrak T(X)})},\varepsilon_{\mathcal M,\mathfrak X})$. By the universal property of $\mathfrak T(X)$, this functor is an equivalence of categories (\cite[Proposition 2.11]{Ogus}).

Given a frame $(X,\overline X,\mathfrak X)$, one can define the category $\Strat(X,\overline X,\mathfrak X)$ of modules with a stratification on $X$ overconvergent along $\overline X-X$, see \cite[P. 50]{ShihoI} where it is denoted by $I^{\dagger}((X,\overline X)/W,\mathfrak X)$). An object ($\calM,\varepsilon$) in $\Strat(X,\overline X,\mathfrak X)$ is a coherent $j^{\dagger}_{X}\mathcal O_{]\overline X[_{\mathfrak X}}$-module together with a $j^{\dagger}_X\mathcal O_{]\overline X[_{\mathfrak X\times \mathfrak X}}$-linear isomorphism $\varepsilon:p_1^*\calM\simeq p_2^*\mathcal M$ satisfying a natural cocycle condition, where the $p_i's$ are the two projection maps 
$]\overline X[_{\mathfrak X\times \mathfrak X}\rightarrow ]\overline X[_{\mathfrak X}$. 
As in the $p$-adically convergent situation, one constructs a functor 
$$(-_{(X,\overline X,\mathfrak X)},\varepsilon_{-,\mathfrak X}):\Isoc^{\dagger}(X,\overline X)\rightarrow \Strat(X,\overline X,\mathfrak X)$$
which is an equivalence of categories (see \cite[Propositions 7.2.2 and 7.3.11]{Lestum}).
\subsection{Relative \texorpdfstring{$p$}--adic cohomology theories} 
\numberwithin{equation}{subsubsection}
\subsubsection{Relative \texorpdfstring{$p$}--adic cohomology theories}\label{relativepadic}
Fix a smooth proper morphism of $k$-varieties $f:Y\rightarrow X$ and a closed immersion $i:X\rightarrow \mathfrak X$, where $\mathfrak X$ is a flat $p$-adic formal scheme. Assume that $f:Y\rightarrow X$ has (log-) smooth parameter in the sense of \cite[Definition 3.4]{ShihoI}.
\begin{remark}\label{logsmooth}
If $f:Y\rightarrow X$ has (log-) smooth parameter, for every morphism of $k$-varieties $Z\rightarrow X$, the base change $Y_Z\rightarrow Z$ has (log-) smooth parameter (\cite[Remark 3.5]{ShihoI}). Moreover, if $X$ is smooth, every smooth proper morphism $f:Y\rightarrow X$ of $k$-varieties has (log-) smooth parameter.
\end{remark}
Depending on the nature of $i:X\rightarrow \mathfrak X$ one defines different $p$-adic cohomology theories:
\begin{itemize}
\item If $X=\mathfrak X_1$ and   $i:X\rightarrow \mathfrak X$ is the canonical inclusion, then one can define the crystalline higher direct image $R^if_{\mathfrak X,\cristalline,*}\mathcal O_{Y/\mathfrak X}$, that is the higher direct image in the relative crystalline site of $X$ in $\mathfrak X$, well defined since $X\numberwithin{equation}{subsection} \subseteq \mathfrak X$ is defined by the ideal $(p)$. It lives in $\mathbf{Coh}(\mathfrak X_K)$, see e.g \cite[Section 1]{ShihoI}.
\item If $i:X\rightarrow \mathfrak X$ is an homeomorphism, then one can define the convergent higher direct image $R^if_{\mathfrak X,\conver,*}\mathcal O_{Y/\mathfrak X}$, that is the higher direct image in the relative convergent site of $X$ in $\mathfrak X$. It lives in $\mathbf{Coh}(\mathfrak X_K)$. See e.g \cite[Sections 2-3]{ShihoI}.
\item If $i:X\rightarrow \mathfrak X$ is an arbitrary  closed immersion, one can define the analytic higher direct image $R^if_{\mathfrak X,\an,*}\mathcal O_{Y/\mathfrak X}$. It is defined via descent using De Rham cohomology, and it lives in $\mathbf{Coh}(]X[_{\mathfrak X})$. For details see \cite[Section 4]{ShihoI}.
\end{itemize}
We complete the picture discussing higher direct images for pairs and frames, in the context of overconvergent isocrystals.
If $(Y,\overline Y)$ is a pair, write $\mathcal O^{\dagger}_{(Y,\overline Y)}\in \oi(Y,\overline Y)$ for the unique overconvergent isocrystal such that, for every frame $(Z,\overline Z,\mathfrak Z)$ over $(Y,\overline Y)$, the restriction of $\mathcal O^{\dagger}_{(Y,\overline Y)}$ to $(Z,\overline Z,\mathfrak Z)$ is given by $j_Z^{\dagger}\mathcal O_{]\overline Z[_{\mathfrak Z}}$.
If $(Y,\overline Y)$ is a pair over a frame $(X,\overline X,\mathfrak X)$ and the first arrow $f:Y\rightarrow X$ is smooth and proper, one can define the overconvergent higher direct image $R^if_{(Y,\overline Y)/\mathfrak X, \rig,*}\mathcal O^{\dagger}_{(Y,\overline Y)}$. It is again defined using De Rham cohomology and descent. See \cite[Section 5]{Shiho} for the definition. It is a  $j^{\dagger}_X\mathcal O_{]\overline X[_{\mathfrak X}}$ module. If $f:(Y,\overline Y)\rightarrow (X,\overline X)$ is cartesian (i.e. $Y=X\times_{\overline X}\overline Y$), it is still an open question whether $R^if_{(Y,\overline Y)/\mathfrak X, \rig,*}\mathcal O^{\dagger}_{(Y,\overline Y)}$ is a coherent $j^{\dagger}_X\mathcal O_{]\overline X[_{\mathfrak X}}$-module.\numberwithin{equation}{subsubsection}
\subsubsection{Comparison}\label{comparisonhigher}
In some particular situation one can compare the various higher direct images defined in Section \ref{relativepadic}.
Assume that $i:X\rightarrow \mathfrak X$ is a closed immersion, $\mathfrak X$ is formally smooth over $W$ and $f:Y\rightarrow X$ is smooth proper with (log-) smooth parameter. 
Using $f$ one considers $(Y,Y)$ as a pair over the frame $(X,X,\mathfrak X)$. The universal $p$-adic enlargement $\mathfrak T(X)$ of $X$ in $\mathfrak X$ induces a commutative diagram:
\begin{center}
\begin{tikzcd}
Y_{\mathfrak T(X)_1}\arrow{r}{f_{\mathfrak T(X)_1}}\arrow{dd}\arrow[phantom]{rdd}{\Box}& \mathfrak T(X)_1\arrow{rr}\arrow{dd}{u}\arrow{rd} && \mathfrak T(X)\arrow{dd}{u}\\
&& \mathfrak X_1\arrow{rd}\\
Y\arrow{r}{f}& X\arrow{rr}\arrow{ru} && \mathfrak X
\end{tikzcd}
\end{center}
By Remark \ref{logsmooth}, the morphism $f_{\mathfrak T(X)_1}:Y_{\mathfrak T(X)_1}\rightarrow \mathfrak T(X)_1$ has (log-) smooth parameter.
In this situation one has $j^{\dagger}_X\mathcal O_{]X[_{\mathfrak X}}=\mathcal O_{]X[_{\mathfrak X}}$, so that
$R^if_{(Y,Y)/\mathfrak X,\rig,*}\mathcal O^{\dagger}_{(Y,Y)}$ and $R^if_{\mathfrak X,\an,*}\mathcal O_{Y/\mathfrak X}$ are coherent $\mathcal O_{]X[_{\mathfrak X}}$-modules (see \cite[Theorem 5.13]{ShihoI}), while $R^if_{\mathfrak T(X)_1,\mathfrak T(X),\an,*}\mathcal O_{Y_{\mathfrak T(X)_1}/\mathfrak T(X)}$, $R^if_{\mathfrak T(X)_1,\mathfrak T(X),\conver,*}\mathcal O_{Y_{\mathfrak T(X)_1}/\mathfrak T(X)}$ and $R^if_{\mathfrak T(X)_1,\mathfrak T(X),\cristalline,*}\mathcal O_{Y_{\mathfrak T(X)_1}/\mathfrak T(X)}$ are coherent $[X]_{\mathfrak X, |p|}=\mathfrak T(X)_K$-modules. 
Write
$$u:[X]_{\mathfrak X, |p|}\rightarrow ]X[_{\mathfrak X}$$
for the natural inclusion.
Essentially by definition (see \cite[Theorem 5.13]{ShihoI} for a much more general statement), one has a canonical isomorphism of $\mathcal O_{]X[_{\mathfrak X}}$-modules
$$R^if_{(Y,Y)/\mathfrak X,\rig,*}\mathcal O^{\dagger}_{(Y,Y)}\simeq R^if_{\mathfrak X,\an,*}\mathcal O_{Y/\mathfrak X}.$$
Pulling back along $u$, one finds canonical isomorphisms of coherent $[X]_{\mathfrak X,|p|}$-modules 
 \numberwithin{equation}{subsubsection}

\begin{equation}\label{equation}
\begin{split}
u^*R^if_{(Y,Y)/\mathfrak X,\rig,*}\mathcal O^{\dagger}_{(Y,Y)}\simeq u^*R^if_{\mathfrak X,\an,*}\mathcal O_{Y/\mathfrak X}\simeq R^if_{\mathfrak T(X)_1,\mathfrak T(X),\an,*}\mathcal O_{Y_{\mathfrak T(X)_1}/\mathfrak T(X)}\\
\simeq R^if_{\mathfrak T(X)_1,\mathfrak T(X),\conver,*}\mathcal O_{Y_{\mathfrak T(X)_1}/\mathfrak T(X)}\simeq R^if_{\mathfrak T(X)_1,\mathfrak T(X),\cristalline,*}\mathcal O_{Y_{\mathfrak T(X)_1}/\mathfrak T(X)}.
\end{split}
\end{equation}
where the second isomorphism comes from \cite[Remark 4.2]{ShihoI}, the third from \cite[Theorem 4.6]{ShihoI} and the last one from \cite[Theorem 2.36]{ShihoI}.

These isomorphisms are functorial in the following sense. Assume that there is a $k$-variety $Z$, a closed embedding $Z\rightarrow \mathfrak Z$ into a  $p$-adic formal scheme $\mathfrak Z$ formally smooth over $W$ and a commutative diagram
\begin{center}
\begin{tikzcd}
Z\arrow{r}\arrow{d}{g}& \mathfrak Z\arrow{d}{g}\\
X\arrow{r} & \mathfrak X.
\end{tikzcd}
\end{center}
By the universal property of $\mathfrak T(X)$, there is an induced map $\mathfrak T(Z)\rightarrow \mathfrak T(X)$ that fits into a commutative diagram
\begin{center}
\begin{tikzcd}[column sep=large, row sep=large]
&&& Y_{Z}\arrow{r}\arrow[near end]{d}{f_Z}& Y\arrow{d}{f}\\
Y_{\mathfrak T(Z)_1}\arrow{r}\arrow[swap]{d}{f_{\mathfrak T(Z)_1}}\arrow{rrru} &Y_{\mathfrak T(X)_1}\arrow[swap, very near start]{d}{f_{\mathfrak T(X)_1}}\arrow{rrru}&& Z\arrow{r}{g}\arrow{d} & X\arrow{d}\\
\mathfrak T(Z)_1\arrow{r}\arrow{d}\arrow{rrru} &\mathfrak T(X)_1\arrow{rrru}\arrow{d}&& \mathfrak Z\arrow{r}{g} & \mathfrak X\\
\mathfrak T(Z)\arrow[swap]{r}{g}\arrow[swap]{rrru}{u} &\mathfrak T(X)\arrow[swap]{rrru}{u},
\end{tikzcd}
\end{center}
where all the squares in the diagram
\begin{center}
\begin{tikzcd}
&& Y_{Z}\arrow{r}\arrow[near end]{d}& Y\arrow{d}\\
Y_{\mathfrak T(Z)_1}\arrow{r}\arrow[swap]{d}\arrow{rru} &Y_{\mathfrak T(X)_1}\arrow[swap, near start]{d}\arrow{rru}& Z\arrow{r}& X\\
\mathfrak T(Z)_1\arrow{r}\arrow{rru} &\mathfrak T(X)_1\arrow{rru}
\end{tikzcd}
\end{center}
are cartesian.
Then the following diagram is commutative
\begin{equation}\label{gigadiagram}
{\begin{tikzcd}
u^*g^*R^if_{(Y,Y)/\mathfrak X,\rig,*}\mathcal O^{\dagger}_{(Y,Y)}\arrow{r}\arrow{d}{\simeq}& u^*R^if_{(Y_{Z},Y_{Z})/\mathfrak Z,\rig,*}\mathcal O^{\dagger}_{(Y_{Z},Y_{Z})}\arrow{d}{\simeq}\\
g^*u^*R^if_{(Y,Y)/\mathfrak X,\rig,*}\mathcal O^{\dagger}_{(Y,Y)}\arrow{r}\arrow{d}{\simeq}& u^*R^if_{(Y_{Z},Y_{Z})/\mathfrak Z,\rig,*}\mathcal O^{\dagger}_{(Y_{Z},Y_{Z})}\arrow{d}{\simeq}\\
g^*u^*R^if_{\mathfrak X,\an,*}\mathcal O_{Y/\mathfrak X}\arrow{r}\arrow{d}{\simeq} & u^*R^if_{\mathfrak Z,\an,*}\mathcal O_{Y_{Z}/\mathfrak Z}\arrow{d}{\simeq}\\
g^*R^if_{\mathfrak T(X)_1,\mathfrak T(X),\an,*}\mathcal O_{Y_{\mathfrak T(X)_1}/\mathfrak T(X)}\arrow{d}{\simeq} \arrow{r}& R^if_{\mathfrak T(Z)_1,\mathfrak T(Z),\an,*}\mathcal O_{Y_{\mathfrak T(Z)_1}/\mathfrak T(Z)}\arrow{d}{\simeq}\\
g^*R^if_{\mathfrak T(X)_1,\mathfrak T(X),\conver,*}\mathcal O_{Y_{\mathfrak T(X)_1}/\mathfrak T(X)}\arrow{r}\arrow{d}{\simeq} & R^if_{\mathfrak T(Z)_1,\mathfrak T(Z),\conver,*}\mathcal O_{Y_{\mathfrak T(Z)_1}/\mathfrak T(Z)}\arrow{d}{\simeq}\\
g^*R^if_{\mathfrak T(X)_1,\mathfrak T(X),\cristalline,*}\mathcal O_{Y_{\mathfrak T(X)_1}/\mathfrak T(X)}\arrow{r} & R^if_{\mathfrak T(Z)_1,\mathfrak T(Z),\cristalline,*}\mathcal O_{Y_{\mathfrak T(Z)_1}/\mathfrak T(Z)}
\end{tikzcd}}
\end{equation}
where the vertical arrows are the isomorphisms in (\ref{equation}) and the horizontal arrows are the base change maps.
\numberwithin{equation}{subsection} 
\subsubsection{Ogus higher direct image}\label{Ogushigher}
Fix a smooth proper morphism $f:Y\rightarrow X$ of $k$-varieties. Write $R^if_{\Ogus,*}\mathcal O_{Y/K}$ in $\Fisoc^{(p)}(X)$ for the Ogus higher direct image (\cite[Section 3, Theorem 3.1]{Ogus}) and recall that its formation is compatible with base change (\cite[Proposition 3.5]{Ogus}).

As an object in $\Isoc^{(p)}(X)$, $R^if_{\Ogus,*}\mathcal O_{Y/K}$ is characterized by the property that for every $p$-adic enlargement $(\mathfrak T,z_{\mathfrak T})$ one has 
$$(R^if_{\Ogus,*}\mathcal O_{Y/K})_{(\mathfrak T,z_{\mathfrak T})}=R^if_{\mathfrak T_1,\mathfrak T,\cristalline,*}\mathcal O_{Y_{\mathfrak T_1}/\mathfrak T}$$
and if $g:(\mathfrak T,z_{\mathfrak T})\rightarrow (\mathfrak Z,z_{\mathfrak Z})$ if a morphism of $p$-adic enlargements, the map
\begin{center}
\begin{tikzcd}
g^*(R^if_{\Ogus,*}\mathcal O_{Y/K})_{(\mathfrak Z,z_{\mathfrak Z})}\arrow{r}{\phi_g}\arrow[equal]{d}&(R^if_{\Ogus,*}\mathcal O_{Y/K})_{(\mathfrak T,z_{\mathfrak T})}\arrow[equal]{d}\\
g^*R^if_{\mathfrak Z_1,\mathfrak Z,\cristalline,*}\mathcal O_{Y_{\mathfrak Z_1}/\mathfrak Z}\arrow{r}{\phi_g} & R^if_{\mathfrak T_1,\mathfrak T,\cristalline,*}\mathcal O_{Y_{\mathfrak T_1}/\mathfrak T}
\end{tikzcd}
\end{center}
is the base change morphism induced by $g$ (see the proof of \cite[Theorem 3.1]{Ogus}).
In particular if $X$ admits a closed immersion into a $p$-adic formal scheme $\mathfrak X$ formally smooth over $W$, the image of $R^if_{\Ogus,*}\mathcal O_{Y/K}$ in $\Strat^{(p)}(X,\mathfrak X)$ is given by the pair $$(R^if_{\mathfrak T(X)_1,\mathfrak T(X),\cristalline,*}\mathcal O_{Y_{\mathfrak T(X)_1}/\mathfrak T(X)},\varepsilon_{R^if_{\Ogus,*}\mathcal O_{Y/K},\mathfrak X}),$$ where $\varepsilon_{R^if_{\Ogus,*}\mathcal O_{Y/K},\mathfrak X}$ is induced by the base change morphisms
$$
\begin{forcedcentertikzcd}[column sep=tiny]
p_1^*R^if_{\mathfrak T(X)_1,\mathfrak T(X),\cristalline,*}\mathcal O_{Y_{\mathfrak T(X)_1}/\mathfrak T(X)}\arrow{r}{\simeq}&
R^if_{\mathfrak T(X)(1)_1,\mathfrak T(X)(1),\cristalline,*}\mathcal O_{Y_{\mathfrak T(X)(1)_1}/\mathfrak T(X)(1)}&
p_2^*R^if_{\mathfrak T(X)_1,\mathfrak T(X),\cristalline,*}\mathcal O_{Y_{\mathfrak T(X)_1}/\mathfrak T(X)}\arrow[swap]{l}{\simeq}.
\end{forcedcentertikzcd}
$$
The Frobenius structure 
$$F_X^*R^if_{\Ogus,*}\mathcal O_{Y/K}\rightarrow R^if_{\Ogus,*}\mathcal O_{Y/K}$$ 
is constructed in the following way (see the proof of \cite[Theorem 3.7]{Ogus}). 
Consider the commutative cartesian diagram
\begin{center}
\begin{tikzcd}
Y'\arrow{d}{f'}\arrow{r}\arrow[phantom]{rd}{\Box}& Y\arrow{d}{f}\\

X\arrow{r}{F_X}& X
\end{tikzcd}
\end{center}
and for every $p$-adic enlargement $(\mathfrak T,z_{\mathfrak T})$ of $X$, consider the following diagram
\begin{center}
\begin{tikzcd}
Y_{\mathfrak T_1} \arrow[bend left]{rrr}{F_{Y_{\mathfrak T_1} }}\arrow[swap]{rr}{\Fr_{Y_{\mathfrak T_1}/\mathfrak T_1}}\arrow{d}{f_{\mathfrak T_1}}&&Y'_{\mathfrak T_1}\arrow{r}\arrow{d}{f'_{\mathfrak T_1}}\arrow[phantom]{rd}{\Box} & Y_{\mathfrak T_1}\arrow{d}{f_{\mathfrak T_1}}\\
\mathfrak T_{1}\arrow{rr}\arrow{d}\arrow[equal]{rr}&&\mathfrak T_1\arrow{r}{F_{\mathfrak T_1}}\arrow{d} & \mathfrak T_1\\
\mathfrak T\arrow[equal]{rr}&& \mathfrak T
\end{tikzcd}
\end{center}
where $\Fr_{Y_{\mathfrak T_1}/\mathfrak T_1}$ is the relative Frobenius morphism.
By the compatibility of $R^if_{\Ogus,*}\mathcal O_{Y/K}$ with base change, there is a canonical isomorphism
$$F_X^*R^if_{\Ogus,*}\mathcal O_{Y/K}\simeq R^if'_{\Ogus,*}\mathcal O_{Y'/K}$$
and hence
$$(F_X^*R^if_{\Ogus,*}\mathcal O_{Y/K})_{(\mathfrak T,z_{\mathfrak T})}\simeq (R^if'_{\Ogus,*}\mathcal O_{Y'/K})_{(\mathfrak T,z_{\mathfrak T})}=R^if'_{\mathfrak T_1,\mathfrak T,\cristalline,*}\mathcal O_{Y'_{\mathfrak T_1}/\mathfrak T}.$$
Then the Frobenius structure is constructed as the base change map
\begin{center}
\begin{tikzcd}
(F_X^*R^if_{\Ogus,*}\mathcal O_{Y/K})_{(\mathfrak T,z_{\mathfrak T})}\arrow{r}\arrow{d}{\simeq}&(R^if_{\Ogus,*}\mathcal O_{Y/K})_{(\mathfrak T,z_{\mathfrak T})}\arrow{d}{\simeq}\\
R^if'_{\mathfrak T_1,\mathfrak T,\cristalline,*}\mathcal O_{Y'_{\mathfrak T_1}/\mathfrak T}\arrow{r} & R^if_{\mathfrak T_1,\mathfrak T,\cristalline,*}\mathcal O_{Y_{\mathfrak T_1}/\mathfrak T}
\end{tikzcd}
\end{center}
induced by $\Fr_{Y_{\mathfrak T_1}/\mathfrak T_1}$.
\numberwithin{equation}{subsubsection} 

\subsubsection{Statement of Theorem \texorpdfstring{\ref{berthelotconjecture}}-}
The aim of this part is to prove the following theorem. 
\begin{theorem}\label{berthelotconjecture}
Assume that $X$ is a smooth $k$-variety and $f:Y\rightarrow X$ a smooth proper morphism. Then 
$R^if_{\Ogus,*}\mathcal O_{Y/K}$ is in the essential image of 
$\Fisoc^{\dagger}(X)\rightarrow \Fisoc^{\dagger}(X,X)\simeq \Fisoc^{(p)}(X)$. 
\end{theorem}
\begin{remark}
Theorem \ref{berthelotconjecture} already appears in the literature as \cite[Corollary 6.2]{berthelotconjecture}, but, as pointed out to us by T.Abe, there might be a gap in the proof. The problem is in the gluing process in \cite[Corollary 6.1]{berthelotconjecture}. The author uses the theory of arithmetic $D$-modules and he tries to compare the higher direct image in that world with $R^i\mathfrak f_{\cristalline,*}\mathcal O_{\mathcal Y/K}(j)$. Locally they coincide, but it is not so clear that the gluing data are compatible, since the isomorphism is defined not at level of complex but only on the level of the derived category. 
So, following a suggestion of T.Abe, we give another proof of Theorem \ref{berthelotconjecture}, using the work of Shiho on the relative log crystalline cohomology (\cite{Shiho}).
\end{remark}
\begin{remark}
The proof actually works more generally for every $E\in \Fisoc^{(p)}(Y)$. The construction of $R^if_{\Ogus,*}E$ does not appear in the literature, so we decided to restrict ourself to $R^if_{\Ogus,*}\mathcal O_{Y/K}$.
\end{remark}

\section{Proof of Theorem \texorpdfstring{\ref{berthelotconjecture}}-} \label{secproof}
\numberwithin{equation}{subsection} 
\subsection{Construction of an overconvergent \texorpdfstring{$F$}--isocrystal}\label{shihosssss}
Fix compactifications $Y\subseteq \overline Y$ and $X\subseteq \overline X$ such that the morphism $f:Y\rightarrow X$ extends to a map of pairs $(Y,\overline Y)\rightarrow (X,\overline X)$ and $X$ (resp. $Y$) is dense in $\overline X$ (resp. $\overline Y$). We start recalling the main result of \cite{Shiho}. This gives a $\mathcal M$ in $\Fisoc^{\dagger}(X,\overline X)$ which, after a base change and on appropriate frames, looks like $R^if_{\Ogus,*}\mathcal O_{Y/K}$. 
To recall the statement, it is helpful to give the following definition.
\begin{definition}
If $(Z,\overline Z, \mathfrak Z)\rightarrow  (X',\overline X',\mathfrak X')$ is a morphism of frames over $(X,\overline X)$ we say that $(Z,\overline Z, \mathfrak Z)$ has $(P_{(X',\overline X',\mathfrak X')})$ if $\overline Z=X\times_{\overline X}\overline Z$ and $\mathfrak Z\rightarrow \mathfrak X'$ is formally smooth.
\end{definition}

By \cite[Theorems 7.6 and 7.9]{Shiho} (and its proof) there exists a frame $(X',\overline X',\mathfrak X')$ over $(X,\overline X)$ such that 
\begin{itemize}
\item $X':=X\times_{\overline X} \overline X'$;
\item $\mathfrak X'$ is formally smooth over $W$;
\item the map $\overline X'\rightarrow \overline X$ is a composition of a surjective proper map followed by a surjective étale map;
\end{itemize}
and an object $\mathcal M$ in $\Foi(X,\overline X)$ with the following properties:
\begin{enumerate}
\item Let $(Z,\overline Z, \mathfrak Z)$ be a frame over $(X',\overline X',\mathfrak X')$ that has $(P_{(X',\overline X',\mathfrak X')})$, so that there is a commutative diagram
 \begin{center}
\begin{tikzcd}
(Y,\overline Y)\arrow{d}\arrow[phantom]{rd}{\Box}&(Y_{X'},\overline Y_{\overline X'})\arrow[phantom]{rd}{\Box}\arrow{d}\arrow{l}& (Y_{Z},\overline Y_{\overline Z})\arrow{l}\arrow{d}\\
(X,\overline X)& (X',\overline X',\mathfrak X')\arrow{l}& (Z,\overline Z,\mathfrak Z)\arrow{l}.
\end{tikzcd}
\end{center}
Then, the image of $\mathcal M$ in $\Strat(Z,\overline Z, \mathfrak Z)$ is given by $$(R^if_{(Y_{Z},\overline Y_{\overline Z})/\mathfrak Z, \rig,*}\mathcal O^{\dagger}_{(Y_Z,\overline Y_{\overline Z})},\varepsilon)$$
where $\varepsilon$ is an isomorphism:
 $$p_1^*R^if_{(Y_{Z},\overline Y_{\overline Z})/\mathfrak Z, \rig,*}\mathcal O^{\dagger}_{(Y_Z,\overline Y_{\overline Z})}\xrightarrow{\simeq} R^if_{(Y_{Z},\overline Y_{\overline Z})/\mathfrak Z\times \mathfrak Z, \rig,*}\mathcal O^{\dagger}_{(Y_Z,\overline Y_{\overline Z})} \xleftarrow{\simeq} p_2^*R^if_{(Y_{Z},\overline Y_{\overline Z})/\mathfrak Z, \rig,*}\mathcal O^{\dagger}_{(Y_Z,\overline Y_{\overline Z})}$$
 and $p_1,p_2:]\overline Z[_{\mathfrak Z\times_W\mathfrak Z}\rightarrow ]\overline Z[_{\mathfrak Z} $ are the projection maps. 
 If moreover $Z=\overline Z$, then $\varepsilon$ is induced by the base change morphisms (\cite[Last paragraph of page 75]{Shiho} and \cite[Theorem 5.14]{ShihoI});
\item Let $h:(Z,\overline Z, \mathfrak Z)\rightarrow (T,\overline T, \mathfrak T)$ be a morphism of frames over $(X',\overline X',\mathfrak X')$ that have $(P_{(X',\overline X',\mathfrak X')})$, so that there is a commutative diagram
 \begin{center}
\begin{tikzcd}
(Y_{Z},\overline Y_{\overline Z})\arrow{d}\arrow{r}& (Y_{T},\overline Y_{\overline T})\arrow{d}\\
(Z,\overline Z,\mathfrak Z)\arrow{r}& (T,\overline T,\mathfrak T).
\end{tikzcd}
\end{center}
 Then, the isomorphism 
\begin{center}
\begin{tikzcd}
h^*\mathcal M_{(Z,\overline Z,\mathfrak Z)}\arrow{r}{\phi_h}\arrow{d}{\simeq} & \mathcal M_{(T,\overline T,\mathfrak T)}\arrow{d}{\simeq}\\
h^*R^if_{(Y_Z,\overline Y_{\overline Z})/\mathfrak Z,\rig,*}\mathcal O^{\dagger}_{(Y_Z,\overline Y_{\overline Z})}\arrow{r} & R^if_{(Y_T,\overline Y_{\overline T})/\mathfrak T,\rig,*}\mathcal O^{\dagger}_{(Y_T,\overline Y_{\overline T})}
\end{tikzcd}
\end{center}
 given by the isocrystals structure is the base change morphism (This is the functoriality in the statement of \cite[Theorem 7.9]{Shiho}, see \cite[Proof of Theorem 4.8]{ShihoI});
\item Let $(Z,\overline Z, \mathfrak Z)$ be a frame over $(X,\overline X)$ that has $(P_{(X',\overline X',\mathfrak X')})$ and assume that $\mathfrak Z$ admits a lifting $\sigma_{\mathfrak Z}$ of $F_{\mathfrak Z_1}$, so that there is a commutative diagram
 \begin{center}
\begin{tikzcd}
(Y_{Z},\overline Y_{\overline Z})\arrow{d}\arrow{r}{(F_{Y_Z},F_{\overline Y_{\overline Z}})}& (Y_{Z},\overline Y_{\overline Z})\arrow{d}\\
(Z,\overline Z,\mathfrak Z)\arrow{r}{(F_Z,F_{\overline Z},\sigma_{\mathfrak Z})}& (Z,\overline Z,\mathfrak Z).
\end{tikzcd}
\end{center}
Then, the isomorphism induced by the Frobenius structure 
 \begin{center}
\begin{tikzcd}
\sigma_{\mathfrak Z}^*\mathcal M_{(Z,\overline Z,\mathfrak Z)}\simeq (F_X^*\mathcal M)_{(Z,\overline Z,\mathfrak Z)}\arrow{r}\arrow{d}{\simeq} & \mathcal M_{(Z,\overline Z,\mathfrak Z)}\arrow{d}{\simeq}\\
\sigma^*_{\mathfrak Z}R^if_{(Y_Z,\overline Y_{\overline Z})/\mathfrak Z,\rig,*}\mathcal O^{\dagger}_{(Y_Z,\overline Y_{\overline Z})}\arrow{r} &R^if_{(Y_Z,\overline Y_{\overline Z})/\mathfrak Z,\rig,*}\mathcal O^{\dagger}_{(Y_Z,\overline Y_{\overline Z})}
\end{tikzcd}
\end{center}
is given by the base change morphism induced by $\sigma_{\mathfrak Z}$ and $F_{Y_{Z}}$ (\cite[Proof of Theorem 7.9]{Shiho}).
\end{enumerate}

\numberwithin{equation}{subsection} 
\subsection{Strategy}\label{strategyp}
To prove Theorem \ref{berthelotconjecture}, it is enough to show that the image of $\mathcal M$ in $\Fisoc^{(p)}(X)$ is isomorphic to $R^{i}f_{\Ogus,*}\mathcal O_{Y/K}$. Since $X'$ does not admit a closed immersion directly in $\mathfrak X'$ and $\mathfrak X'$ does not admits a lifting of the absolute Frobenius of $\mathfrak X'_1$, one can't use directly the description of $\mathcal M$ given in previous Section \ref{shihosssss}. But there exists\footnote{To construct it, consider a finite covering $\{\mathrm{Spf}(A_i)\}$  of $\mathfrak X'$ by formal affine open sub schemes such that every $\mathrm{Spf}(A_i)$ admits a formally étale morphism to an affine formal space. Then $\{\Spec(A_{i,1})\}$ is a covering of $\mathfrak X'_1$ by affine open sub schemes and $\{V_i:=\Spec(A_{i,1})\times_{\mathfrak X'_1} X'\}$ is a Zariski open covering of $X'$. Consider a finite covering $\{U_{i,j}\}$ of $V_i$ by affine open sub schemes. Then the maps $U_{i,j}\rightarrow \Spec(A_{i,1})$ are affine and of finite type, so that there are closed immersions $U_{i,j}\rightarrow \mathbb A^{n_{i,j}}_{\Spec(A_{i,1})  }$. Write $\mathfrak U_{i,j}$ for the formal affine space of dimension $n_{i,j}$ over $\mathrm{Spf}(A_i)$. Then $U:=\coprod_{i,j}U_{i,j}$ admits a closed immersion into  $\mathfrak U:=\coprod_{i,j}\mathfrak U_{i,j}$ and $\mathfrak U$ is formally smooth over $\mathfrak X'$. To show that $\mathfrak U$ admits a lifting of $F_{\mathfrak U_1}$ it is enough to show that each $\mathfrak U_{i,j}$ admits a lifting of $F_{\mathfrak U_{i,j,1}}$. This follows from the fact that $\mathfrak U_{i,j}$ is formally affine admitting a formally étale morphism to an formal affine space.} an étale surjective morphism $U\rightarrow X'$ such that $U$ admits a closed immersion into a $p$-adic formal scheme $\mathfrak U$ which is formally smooth over $\mathfrak X'$ and it is endowed with a lifting $\sigma_{\mathfrak U}$ of $F_{\mathfrak U_1}$.
Write $g$ for the composition $U\rightarrow X'\rightarrow X$. 

To prove Theorem \ref{berthelotconjecture}, first one constructs an isomorphism 
$$\psi:g^*\mathcal M\simeq g^*R^{i}f_{\Ogus,*}\mathcal O_{Y/K}\simeq R^{i}f_{U,\Ogus,*}\mathcal O_{Y_U/K}\text{ in } \Fisoc^{(1)}(U)$$ where the isomorphism on the right comes from the fact that the formation of $R^{i}f_{\Ogus,*}\mathcal O_{Y/K}$ is compatible with base change, see Section \ref{Ogushigher}. Then one uses étale and proper descent for convergent isocrystals to deduce that $\psi$ descent to $\Fisoc^{(p)}(X)$.
More precisely the proof decomposes as follows:
\begin{enumerate}
\item One constructs an isomorphism $$\psi:g^*\mathcal M\simeq R^{i}f_{U,\Ogus,*}\mathcal O_{Y_U/K}$$ in $\Isoc^{(p)}(U)\simeq \Strat^{(p)}(U,\mathfrak U)$. This is done in Section \ref{comparisonisocrystal},  using that $(U,U,\mathfrak U)$ has $(P_{(X',\overline X',\mathfrak X')})$ (so that one can apply the property $(1)$ of $\mathcal M$) and the comparison isomorphisms in \ref{comparisonhigher};
\item One verifies that the $\psi$ commutes with the Frobenius structures i.e. that $\psi$ makes the following diagram commutative
\begin{center}
\begin{tikzcd}
F_U^*g^*\mathcal M\arrow{r}\arrow{d}{F_U^*\psi} & g^*\mathcal M\arrow{d}{\psi}\\
F_U^*R^if_{U,\Ogus,*}\mathcal O_{Y_U/K}\arrow{r} & R^if_{U,\Ogus,*}\mathcal O_{Y_U/K}
\end{tikzcd}
\end{center}
in $\Isoc^{(p)}(U)\simeq \Strat^{(p)}(U,\mathfrak U)$.  This is done in Section \ref{comparisonfrobenius}, using that $\mathfrak U$ has a lifting of $F_{\mathfrak U_1}$ (so that one can apply the property $(3)$ of $\mathcal M$) and the comparison isomorphisms in \ref{comparisonhigher};
\item By the equivalence $\mathrm{F1}\mbox{-}\mathrm{Fp}$ in Section \ref{functors}, the first two steps imply that there is an isomorphism 
$$\psi:g^*\mathcal M\simeq R^{i}f_{U,\Ogus,*}\mathcal O_{Y_U/K}$$ in $\Fisoc^{(1)}(U)$;
\item To apply descent for convergent isocrystals, one has to check that $\psi$ makes the following diagram in $\Fisoc^{(1)}(U\times_X U)$ commutative:
\begin{center}
\begin{tikzcd}
q_1^*g^*\mathcal M\arrow{r}\arrow{d}{q_1^*\psi} & q_2^*g^*\mathcal M\arrow{d}{q^*_2\psi}\\
q_1^*g^*R^if_{\Ogus,*}\mathcal O_{Y/K}\arrow{r} & q_2^*g^*R^if_{\Ogus,*}\mathcal O_{Y/K}
\end{tikzcd}
\end{center}
where $q_1,q_2:U\times_XU\rightarrow U$ are the projections. To check this, by the equivalence $\mathrm{F1}\mbox{-}\mathrm{Fp}$, it is enough to show that it is commutative in $\Fisoc^{(p)}(U\times_X U)$ or equivalently in $\Isoc^{(p)}(U\times_X U)\simeq \Strat^{(p)}(U\times_X U,\mathfrak U\times_W \mathfrak U)$. This is done in Section \ref{descent}, using that $q_1,q_2:(U\times_X U,U\times_X U,\mathfrak U\times_W \mathfrak U)\rightarrow (U,U,\mathfrak U)$ are morphisms of frames that have $(P_{(X',\overline X',\mathfrak X')})$ (so that one can apply the property $(2)$ of $\mathcal M$) and the comparison isomorphisms in \ref{comparisonhigher}. 
\end{enumerate}  
\begin{remark}
The reason why one needs to pass back and forth between $\Fisoc^{(p)}(U)$ and $\Fisoc^{(1)}(U)$ is that proper descent is not known for the category $\Isoc^{(p)}(U)$, while proper descent for the category $\Isoc^{(1)}(U)$ (and hence for $\Fisoc^{(1)}(U)$) is proved in \cite{Ogus}. On the other hand one knows the value of $R^if_{U,\Ogus,*}\mathcal O_{Y_U/K}$ only on $p$-adic enlargements. The equivalences of categories in Section \ref{functors} allow to combine these informations. 
\end{remark}
\subsection{Comparison of isocrystals}\label{comparisonisocrystal}
In this section we construct an isomorphism
$$\psi:g^*\mathcal M\simeq R^if_{U,\Ogus,*}\mathcal O_{Y_{U}/K} \text{ in } \Isoc^{(p)}(U).$$
Consider the universal $p$-adic enlargements $\mathfrak T(U)$ and $\mathfrak T(U)(1)$ of $U$ in $\mathfrak U$ and $\mathfrak U\times \mathfrak U$ and write $u,p_{1},p_2$ for the natural morphisms
\begin{center}
\begin{tikzcd}
(\mathfrak T(U)(1)_{1},\mathfrak T(U)(1)_{1},\mathfrak T(U)(1))\arrow{r}{u}\arrow[d, shift left=1.5ex,"p_1"]\arrow[swap]{d}{p_2} & (U,U,\mathfrak U\times \mathfrak U)\arrow[d, shift left=1.5ex,"p_1"]\arrow[swap]{d}{p_2}\\
(\mathfrak T(U)_{1},\mathfrak T(U)_{1},\mathfrak T(U))\arrow{r}{u} & (U,U,\mathfrak U).
\end{tikzcd}
\end{center}
Since $(U,U,\mathfrak U)$ has $(P_{(X',\overline X',\mathfrak X')})$, by the property $(1)$ of $\mathcal M$ in \ref{shihosssss}, one has:
$$\mathcal M_{(U,U,\mathfrak U)}=R^if_{(Y_{U},Y_{U})/\mathfrak U, \rig,*}\mathcal O^{\dagger}_{(Y_{U},Y_{U})} \text{ in } \mathbf{Coh}(]U[_{\mathfrak U}).$$
Since $\mathcal M$ is an isocrystal, one gets $\mathcal M_{(\mathfrak T(U)_1,\mathfrak T(U)_1,\mathfrak T(U))}\simeq u^*\mathcal M_{(U,U,\mathfrak U)} \text{ in } \mathbf{Coh}(\mathfrak T(U)_K).$
Then as in (\ref{equation}): 
$$\mathcal M_{(\mathfrak T(U),z_{\mathfrak T(U)})}\simeq u^*R^if_{(Y_{U},Y_{U})/\mathfrak U, \rig,*}\mathcal O^{\dagger}_{(Y_{U},Y_{U})}\simeq u^*R^if_{U,\mathfrak U,\an,*}\mathcal O_{Y_{U}/\mathfrak U}\simeq R^if_{\mathfrak T(U)_1,\mathfrak T(U),\an,*}\mathcal O_{Y_{\mathfrak T(U)_{1}}/\mathfrak T(U)}$$
$$\simeq R^if_{\mathfrak T(U)_{1},\mathfrak T(U),\conver,*}\mathcal O_{Y_{\mathfrak T(U)_{1}}/\mathfrak T(U)}\simeq R^if_{\mathfrak T(U)_{1},\mathfrak T(U),\cristalline,*}\mathcal O_{Y_{\mathfrak T(U)_{1}}/\mathfrak T(U)} \text{ in } \mathbf{Coh}(\mathfrak T(U)_K).$$

Since, by construction (\ref{Ogushigher}), $$R^if_{\mathfrak T(U)_{1},\mathfrak T(U),\cristalline,*}\mathcal O_{Y_{\mathfrak T(U)_{1}}/\mathfrak T(U)}=(R^if_{U,\Ogus,*}\mathcal O_{Y_{U}/K})_{(\mathfrak T(U),z_{\mathfrak T(U)})}$$
one has an isomorphism 
$$\psi:\mathcal M_{(\mathfrak T(U),z_{\mathfrak T(U)})}\simeq (R^if_{U,\Ogus,*}\mathcal O_{Y_{U}/K})_{(\mathfrak T(U),z_{\mathfrak T(U)})} \text{ in } \mathbf{Coh}(\mathfrak T(U)_K).$$
To promote $\psi$ to an isomorphism in $\Strat^{(p)}(U,\mathfrak U)\simeq \Isoc^{(p)}(U)$ one has to check that $\psi$ is compatible with the stratifications $\varepsilon_{g^*\mathcal M,\mathfrak U}$ on $g^*\mathcal M$ and $\varepsilon_{R^if_{U,\Ogus,*}\mathcal O_{Y_{U}/K},\mathfrak U}$ on $(R^if_{U,\Ogus,*}\mathcal O_{Y_{U}/K})_{(\mathfrak T(U),z_{\mathfrak T(U)})}$.

Since $(U,U,\mathfrak U)$ has $(P_{(X',\overline X',\mathfrak X')})$, by the property (1) in \ref{shihosssss}, the stratification $\varepsilon_{g^*\mathcal M,\mathfrak U}$ is given by the base change morphisms:
$$p_1^*R^if_{(Y_{U},Y_{U})/\mathfrak U, \rig,*}\mathcal O^{\dagger}_{(Y_{U},Y_{U})}\rightarrow R^if_{(Y_{U},Y_{U})/\mathfrak U\times \mathfrak U, \rig,*}\mathcal O^{\dagger}_{(Y_{U},Y_{U})} \leftarrow p_2^*R^if_{(Y_{U},Y_{U})/\mathfrak U, \rig,*}\mathcal O^{\dagger}_{(Y_{U},Y_{U})}.$$
As in (\ref{gigadiagram}) pulling back to $u^*$, one has a commutative diagram
\begin{center}
\begin{tikzpicture}[baseline= (a).base]
\node[scale=0.95] (a) at (0,0){
\begin{tikzcd}[column sep=tiny]
u^*p_1^*R^if_{(Y_{U},Y_{U})/\mathfrak U, \rig,*}\mathcal O^{\dagger}_{(Y_{U},Y_{U})}\arrow{r}{\simeq}\arrow{d}{\simeq} &\arrow{d}{\simeq} u^*R^if_{(Y_{U},Y_{U})/\mathfrak U\times \mathfrak U, \rig,*}\mathcal O^{\dagger}_{(Y_{U},Y_{U})} & \arrow{d}{\simeq}u^*p_2^*R^if_{(Y_{U},Y_{U})/\mathfrak U, \rig,*}\mathcal O^{\dagger}_{(Y_{U},Y_{U})} \arrow{l}{\simeq}\arrow{d}{\simeq}\\
p_1^*u^*R^if_{(Y_{U},Y_{U})/\mathfrak U, \rig,*}\mathcal O^{\dagger}_{(Y_{U},Y_{U})}\arrow{r}{\simeq}\arrow{d}{\simeq} &\arrow{d}{\simeq} u^*R^if_{(Y_{U},Y_{U})/\mathfrak U\times \mathfrak U, \rig,*}\mathcal O^{\dagger}_{(Y_{U},Y_{U})} & \arrow{d}{\simeq}p_2^*u^*R^if_{(Y_{U},Y_{U})/\mathfrak U, \rig,*}\mathcal O^{\dagger}_{(Y_{U},Y_{U})} \arrow{l}{\simeq}\arrow{d}{\simeq}\\
p_1^*u^*R^if_{U,\mathfrak U,\an,*}\mathcal O_{Y_{U}/\mathfrak U}\arrow{r}{\simeq}\arrow{d}{\simeq} & u^*R^if_{U,\mathfrak U\times \mathfrak U,\an,*}\mathcal O_{Y_{U}/\mathfrak U\times \mathfrak U}\arrow{d}{\simeq}& p_2^*u^*R^if_{U,\mathfrak U,\an,*}\mathcal O_{Y_{U}/\mathfrak U} \arrow{l}{\simeq}\arrow{d}{\simeq}\\
p_1^*R^if_{\mathfrak T(U)_{1},\mathfrak T(U),\cristalline,*}\mathcal O_{Y_{\mathfrak T(U)_1}/\mathfrak T(U)}\arrow{r}{\simeq}&R^if_{\mathfrak T(U)(1)_{1},\mathfrak T(U)(1),\cristalline,*}\mathcal O_{Y_{\mathfrak T(U)(1)_1}/\mathfrak T(U)(1)} &\arrow{l}{\simeq}p_2^*R^if_{\mathfrak T(U)_{1},\mathfrak T(U),\cristalline,*}\mathcal O_{Y_{\mathfrak T(U)_1}/\mathfrak T(U)} \arrow{l}{\simeq},
\end{tikzcd}
};
\end{tikzpicture}
\end{center}

where the horizontal maps are the natural base change maps. So, the stratification $\varepsilon_{g^*\mathfrak M,\mathfrak U}$ on $g^*\mathcal M$
\begin{center}
\begin{tikzpicture}[baseline= (a).base]
\node[scale=0.95] (a) at (0,0){
\begin{tikzcd}[column sep=tiny]
p_1^*\mathcal M_{(\mathfrak T(U),z_{\mathfrak T(U)})}\arrow{r}\arrow{d}{\simeq}&\mathcal M_{(\mathfrak T(U)(1),z_{\mathfrak T(U)(1)})}\arrow{d}{\simeq}& \arrow{l}p_2^*\mathcal M_{(\mathfrak T(U),z_{\mathfrak T(U)})}\arrow{d}{\simeq}\\
p_1^*R^if_{\mathfrak T(U)_{1},\mathfrak T(U),\cristalline,*}\mathcal O_{Y_{\mathfrak T(U)_1}/\mathfrak T(U)}\arrow{r}{\simeq}&R^if_{\mathfrak T(U)(1)_{1},\mathfrak T(U)(1),\cristalline,*}\mathcal O_{Y_{\mathfrak T(U)(1)_1}/\mathfrak T(U)(1)} &\arrow{l}{\simeq}p_2^*R^if_{\mathfrak T(U)_{1},\mathfrak T(U),\cristalline,*}\mathcal O_{Y_{\mathfrak T(U)_1}/\mathfrak T(U)} \arrow{l}{\simeq}
\end{tikzcd}
};
\end{tikzpicture}
\end{center}
is induced by the base change morphisms.
Since $\varepsilon_{R^if_{U,\Ogus,*}\mathcal O_{Y_{U}/K},\mathfrak U}$ is induced by the base change morphisms by construction (\ref{Ogushigher}), one concludes that 
$$\psi:\mathcal M_{(\mathfrak T(U),z_{\mathfrak T(U)})}\simeq (R^if_{U,\Ogus,*}\mathcal O_{Y_{U}/K})_{(\mathfrak T(U),z_{\mathfrak T(U)})}\text{ in } \mathbf{Coh}(\mathfrak T(U)_K)$$ is compatible with the stratifications and hence induces an isomorphism 
$$\psi: g^*\mathcal M\simeq R^if_{U,\Ogus,*}\mathcal O_{Y_{U}/K}  \text{ in } \Strat^{(p)}(U,\mathfrak U)\simeq \Isoc^{(p)}(U).$$
\subsection{Comparison of Frobenius structures}\label{comparisonfrobenius}
We now check that $\psi$ is compatible with the Frobenius structures, i.e. that the following diagram in $\Isoc^{(p)}(U)$
is commutative:
\begin{center}
\begin{tikzcd}
F_U^*g^*\mathcal M\arrow{r}\arrow{d}{F_U^*\psi} & g^*\mathcal M\arrow{d}{\psi}\\
F_U^*R^if_{U,\Ogus,*}\mathcal O_{Y_U/K}\arrow{r} & R^if_{U,\Ogus,*}\mathcal O_{Y_U/K}
\end{tikzcd}
\end{center}
Since 
$$(-_{(\mathfrak T(U),z_{\mathfrak T(U)})},\varepsilon_{-,\mathfrak U}):\Isoc^{(p)}(U)\rightarrow \Strat^{(p)}(U,\mathfrak U)$$
is an equivalence of categories, it is enough to show that 
\begin{center}
\begin{tikzcd}
(F_U^*g^*\mathcal M)_{(\mathfrak T(U),z_{\mathfrak T(U)})}\arrow{r}\arrow{d} & g^*\mathcal M_{(\mathfrak T(U),z_{\mathfrak T(U)})}\arrow{d}\\
(F_U^*R^if_{U,\Ogus,*}\mathcal O_{Y_U/K})_{(\mathfrak T(U),z_{\mathfrak T(U)})}\arrow{r} & (R^if_{U,\Ogus,*}\mathcal O_{Y_U/K})_{(\mathfrak T(U),z_{\mathfrak T(U)})}
\end{tikzcd}
\end{center}
is commutative.
Since ($U,U,\mathfrak U)$ has $(P_{(X',\overline X',\mathfrak X')})$ and it is endowed with a morphism $\sigma_{\mathfrak U}$ lifting $F_{\mathfrak U_1}$, by the property $(3)$ of $\mathcal M$ in \ref{shihosssss}, the Frobenius structure on $\mathcal M_{(U,U,\mathfrak U)}$ is given by the base change map induced by $\sigma_{\mathfrak U}$ and $F_{Y_{U}}:$
$$\sigma_{\mathfrak U}^*R^if_{(Y_{U},Y_{U})/\mathfrak U, \rig,*}\mathcal O^{\dagger}_{(Y_{U},Y_{U})}\rightarrow R^if_{(Y_{U},Y_{U})/\mathfrak U, \rig,*}\mathcal O^{\dagger}_{(Y_{U},Y_{U})}.$$
By the universal property of the universal $p$-adic enlargement one gets a commutative diagram:
\begin{center}
\begin{tikzcd}
U\arrow{r}\arrow{d}{F_U} & \mathfrak T(U)\arrow{r}{u}\arrow{d}{\sigma_{\mathfrak T(U)}} & \mathfrak U\arrow{d}{\sigma_{\mathfrak U}}\\
U\arrow{r} & \mathfrak T(U)\arrow{r}{u}& \mathfrak U
\end{tikzcd}
\end{center}
Pulling back via $u$, there is a commutative diagram (\ref{gigadiagram}):
\begin{center}
\begin{tikzcd}
u^*\sigma_{\mathfrak U}^*R^if_{(Y_{U},Y_{U})/\mathfrak U, \rig,*}\mathcal O^{\dagger}_{(Y_{U},Y_{U})}\arrow{r}\arrow{d}{\simeq} & u^*R^if_{(Y_{U},Y_{U})/\mathfrak U, \rig,*}\mathcal O^{\dagger}_{(Y_{U},Y_{U})}\arrow{d}{\simeq}\\
\sigma_{\mathfrak T(U)}^*u^*R^if_{(Y_{U},Y_{U})/\mathfrak U, \rig,*}\mathcal O^{\dagger}_{(Y_{U},Y_{U})}\arrow{r}\arrow{d}{\simeq} & u^*R^if_{(Y_{U},Y_{U})/\mathfrak U, \rig,*}\mathcal O^{\dagger}_{(Y_{U},Y_{U})}\arrow{d}{\simeq}\\
\sigma_{\mathfrak T(U)}^*R^if_{\mathfrak T(U)_1,\mathfrak T(U),\an,*}\mathcal O_{Y_{\mathfrak T(U)_1}/\mathfrak T(U)}\arrow{r}\arrow{d}{\simeq} & R^if_{\mathfrak T(U)_1\mathfrak T(U),\an,*}\mathcal O_{Y_{\mathfrak T(U)_1}/\mathfrak T(U)}\arrow{d}{\simeq}\\
\sigma_{\mathfrak T(U)}^* R^if_{\mathfrak T(U)_1,\mathfrak T(U),\cristalline,*}\mathcal O_{Y_{\mathfrak T(U)_1}/\mathfrak T(U)}\arrow{r} & R^if_{\mathfrak T(U)_1,\mathfrak T(U),\cristalline,*}\mathcal O_{Y_{\mathfrak T(U)_1}/\mathfrak T(U)}
\end{tikzcd}
\end{center}
where the horizontal morphisms are the base change morphisms. 
So the morphism
\begin{center}
\begin{tikzcd}
\sigma_{\mathfrak T}^*\mathcal M_{(\mathfrak T(U),z_{\mathfrak T(U)})}\simeq (F_U^*g^*\mathcal M)_{(\mathfrak T(U),z_{\mathfrak T(U)})}\arrow{r}\arrow{d}{\simeq} & g^*\mathcal M_{(\mathfrak T(U),z_{\mathfrak T(U)})}\arrow{d}{\simeq}\\
\sigma_{\mathfrak T(U)}^* R^if_{\mathfrak T(U)_1,\mathfrak T(U),\cristalline,*}\mathcal O_{Y_{\mathfrak T(U)_1}/\mathfrak T(U)}\arrow{r} & R^if_{\mathfrak T(U)_1,\mathfrak T(U),\cristalline,*}\mathcal O_{Y_{\mathfrak T(U)_1}/\mathfrak T(U)}
\end{tikzcd}
\end{center}
is induced by the base change morphism for $F_{Y_{\mathfrak T(U)_1}}$ and $\sigma_{\mathfrak T(U)}$.

We check that the same is true for $R^if_{U,\Ogus,*}\mathcal O_{Y_U/K}$. Consider the commutative diagram
\begin{center}
\begin{tikzcd}
Y_{\mathfrak T(U)_1} \arrow[bend left]{rrr}{F_{Y_{\mathfrak T(U)_1} }}\arrow[swap]{rr}{\Fr_{Y_{\mathfrak T(U)_1}/\mathfrak T(U)_1}}\arrow{d}{f_{\mathfrak T(U)_1}}&&Y'_{\mathfrak T(U)_1}\arrow{r}\arrow{d}{f'_{\mathfrak T(U)_1}}\arrow[phantom]{rd}{\Box} & Y_{\mathfrak T(U)_1}\arrow{d}{f_{\mathfrak T(U)_1}}\\
\mathfrak T(U)_{1}\arrow{rr}\arrow{d}\arrow[equal]{rr}&&\mathfrak T(U)_1\arrow{r}{F_{\mathfrak T(U)_1}}\arrow{d} & \mathfrak T(U)_1\arrow{d}\\
\mathfrak T(U)\arrow[equal]{rr}&& \mathfrak T(U)\arrow{r}{\sigma_{\mathfrak T(U)}} &\mathfrak T(U)
\end{tikzcd}
\end{center}
and recall (\ref{Ogushigher}) that the Frobenius structure is defined by the base change map induced by $\Fr_{Y_{\mathfrak T(U)_1}/\mathfrak T(U)_1}$:
$$(F_U^*R^if_{U,\Ogus,*}\mathcal O_{Y_U/K})_{(\mathfrak T(U)_1,\mathfrak T(U)_1,\mathfrak T(U))}\simeq R^if'_{\mathfrak T(U)_1,\mathfrak T(U),\cristalline,*}\mathcal O_{Y'_{\mathfrak T(U)_1}/\mathfrak T(U)}\rightarrow  R^if_{\mathfrak T(U)_1,\mathfrak T(U),\cristalline,*}\mathcal O_{Y_{\mathfrak T(U)_1}/\mathfrak T(U)}$$
Since there is a lifting $\sigma_{\mathfrak T(U)}$ of $F_{\mathfrak T(U)_1}$, there is a morphism of enlargements
$$\sigma_{\mathfrak T(U)}:(\mathfrak T(U),z_{\mathfrak T(U)})\rightarrow (\mathfrak T(U),z_{\mathfrak T(U)}\circ F_{\mathfrak T(U)_1}).$$ Since $F^*_{U}R^{i}f_{\Ogus,*}\mathcal O_{\mathcal Y_U/K} $ is a crystal, there is an isomorphism
$$\sigma_{\mathfrak T(U)}^*R^if_{\mathfrak T(U)_1,\mathfrak T(U),\cristalline,*}\mathcal O_{Y_{\mathfrak T(U)_1}/\mathfrak T(U)}\rightarrow R^if'_{\mathfrak T(U)_1,\mathfrak T(U),\cristalline,*}\mathcal O_{Y'_{\mathfrak T(U)_1}/\mathfrak T(U)}.$$
that, as recalled in section \ref{Ogushigher}, identifies with the base change map induced by $\sigma_{\mathfrak T(U)}$ and $F_{\mathfrak T(U)_1}$.
So the Frobenius structure

\begin{center}
\begin{tikzpicture}[baseline= (a).base]
\node[scale=0.95] (a) at (0,0){
\begin{tikzcd}
\sigma_{\mathfrak T(U)}^*(R^if_{U,\Ogus,*}\mathcal O_{Y_U/K})_{(\mathfrak T(U),z_{\mathfrak T(U)})}\simeq (F_U^*R^if_{U,\Ogus,*}\mathcal O_{Y_U/K})_{(\mathfrak T(U),z_{\mathfrak T(U)})}\arrow{r}\arrow{d}{\simeq} & (R^if_{U,\Ogus,*}\mathcal O_{Y_U/K})_{(\mathfrak T(U),z_{\mathfrak T(U)})}\arrow{d}{\simeq}\\
\sigma_{\mathfrak T(U)}^* R^if_{\mathfrak T(U)_1,\mathfrak T(U),\cristalline,*}\mathcal O_{Y_{\mathfrak T(U)_1}/\mathfrak T(U)}\arrow{r} & R^if_{\mathfrak T(U)_1,\mathfrak T(U),\cristalline,*}\mathcal O_{Y_{\mathfrak T(U)_1}/\mathfrak T(U)}
\end{tikzcd}
};
\end{tikzpicture}
\end{center}
is given by the composition of the base change maps induced by $\Fr_{Y_{\mathfrak T(U)_1}/\mathfrak T(U)_1}$ followed by the base change map induced by $\sigma_{\mathfrak T(U)}$ and $F_{\mathfrak T(U)_1}$, hence it is given by the base change map induced by $F_{Y_{\mathfrak T(U)_1}}$ and $\sigma_{\mathfrak T(U)}$.

In conclusion, $\psi$ is compatible with the Frobenius structures of $\mathcal M$ and $R^if_{U,\Ogus,*}\mathcal O_{Y_U/K}$, so that $\psi$ is gives an isomorphism 
$$\psi: R^if_{U,\Ogus,*}\mathcal O_{Y_U/K}\simeq g^*\mathcal M \text{ in } \Fisoc^{(p)}(U)$$ and hence an isomorphism 
$$\psi: R^if_{U,\Ogus,*}\mathcal O_{Y_U/K}\simeq g^*\mathcal M \text{ in } \Fisoc^{(1)}(U).$$
\subsection{Descent}\label{descent}
 \numberwithin{equation}{subsection}

Now one has to descend from $U$ to $X$.
To do this, consider the closed immersion
$$U\times_XU\rightarrow U\times_{k} U\rightarrow \mathfrak U_1\times_{k}\mathfrak U_1\rightarrow \mathfrak U\times_{W}\mathfrak U$$
where the first map is a closed immersion by \cite[Tag 01KR]{stacks-project} since $X$ is separated.
Write $\mathfrak T(U\times_X U)$ for the universal $p$-adic enlargement of $U\times_XU$ in $\mathfrak U\times_{W}\mathfrak U$ and $q_1,q_2$ for the projections $$U\times_X U\rightarrow U \quad \text{and}\quad \mathfrak U\times_{W}\mathfrak U\rightarrow \mathfrak U.$$ Finally write $u_{\mathfrak T(U\times_X U)},q_{\mathfrak T(U\times_X U),1}$ and $q_{\mathfrak T(U\times_X U),2}$ for the natural morphisms:
\begin{center}
\begin{tikzcd}
(\mathfrak T(U\times_X U)_{1},\mathfrak T(U\times_X U)_{1},\mathfrak T(U\times_X U))\arrow{r}{u_{\mathfrak T(U\times_X U)}}\arrow[d, shift left=1.5ex,"q_{\mathfrak T(U\times_X U),1}"]\arrow[swap]{d}{q_{\mathfrak T(U\times_X U),2}} & (U\times_X U,U\times_X U, \mathfrak U\times_{W} \mathfrak U)\arrow[d, shift left=1.5ex,"q_1"]\arrow[swap]{d}{q_2}\\
(\mathfrak T(U)_{1},\mathfrak T(U)_{1},\mathfrak T(U))\arrow{r}{u} & (U,U,\mathfrak U)
\end{tikzcd}
\end{center} 
and $g'$ for $U\times_X U\rightarrow X$.
By the equivalence $\mathrm{F1}\mbox{-}\mathrm{Fp}$ in Section \ref{functors}, to show that the descent diagram in $\Fisoc^{(1)}(U\times_X U)$
\begin{equation}\label{descentdiagram}
{
\begin{tikzcd}
q_1^*g^*\mathcal M\arrow{r}\arrow{d} & q_2^*g^*\mathcal M\arrow{d}\\
q_1^*g^*R^if_{\Ogus,*}\mathcal O_{Y/K}\arrow{r} & q_2^*g^*R^if_{\Ogus,*}\mathcal O_{Y/K}
\end{tikzcd}}
\end{equation}
is commutative, it is enough to show that it is commutative in $\Fisoc^{(p)}(U\times_XU)$. Then one decomposes \ref{descentdiagram} as follows:
\begin{center}
\begin{tikzcd}
q_1^*g^*\mathcal M\arrow{r}\arrow{d} &g^{',*}\mathcal M\arrow{d} & q_2^*g^*\mathcal M\arrow{d}\arrow{l}\\
q_1^*R^if_{\Ogus,*}\mathcal O_{Y/K}\arrow{r} &g^{',*}R^if_{\Ogus,*}\mathcal O_{Y/K}\arrow{r} & q_2^*R^if_{\Ogus,*}\mathcal O_{Y/K}\arrow{l}
\end{tikzcd}
\end{center}
So it is enough to show that, for $?\in \{1,2\}$, the following diagram is commutative
\begin{center}
\begin{tikzcd}
q_{?}^*g^{*}\mathcal M\arrow{r}\arrow{d} & g^{',*}\mathcal M\arrow{d}\\
q_{?}^*g^{*}R^if_{\Ogus,*}\mathcal O_{Y/K}\arrow{r} & g^{',*}R^if_{\Ogus,*}\mathcal O_{Y/K}
\end{tikzcd}
\end{center}
in $\Fisoc^{(p)}(U\times_X U)$ or, equivalently, in $\Isoc^{(p)}(U\times_X U)$.
Since $$(-_{(\mathfrak T(U\times_X U),z_{\mathfrak T(U\times_X U)})},\varepsilon_{-,\mathfrak U\times_W \mathfrak U}):\Isoc^{(p)}(U\times_X U)\rightarrow \Strat^{(p)}(U\times_X U,\mathfrak U\times_W\mathfrak U)$$
is an equivalence of categories, it is enough to show that  
\begin{center}
\begin{tikzcd}
(q_{?}^*g^{*}\mathcal M)_{(\mathfrak T(U\times_X U),z_{\mathfrak T(U\times_X U)})}\arrow{r}\arrow{d} & (g^{',*}\mathcal M)_{(\mathfrak T(U\times_X U),z_{\mathfrak T(U\times_X U)})}\arrow{d}\\
(q_{?}^*g^{*}R^if_{\Ogus,*}\mathcal O_{Y/K})_{(\mathfrak T(U\times_X U),z_{\mathfrak T(U\times_X U)})}\arrow{r} & (g^{',*}R^if_{\Ogus,*}\mathcal O_{Y/K})_{(\mathfrak T(U\times_X U),z_{\mathfrak T(U\times_X U)})}
\end{tikzcd}
\end{center}
commutes.
Since $q_?:(U\times_X U,U\times_X U,\mathfrak U\times_{W} \mathfrak U)\rightarrow (U,U,\mathfrak U)$ is a morphism of frame over $(X',\overline X', \mathfrak X')$ that have ($P_{(X',\overline X', \mathfrak X')}$), by the property $(2)$ of $\mathcal M$ in \ref{shihosssss}, the morphism given by the isocrystals structure 
$$q_{?}^*R^if_{(Y_{U},Y_{U})/\mathfrak U, \rig,*}\mathcal O^{\dagger}_{(Y_{U},Y_{U})}\rightarrow R^if_{(Y_{U\times_{X}U},Y_{U\times_{X}U})/\mathfrak U\times \mathfrak U, \rig,*}\mathcal O^{\dagger}_{(Y_{U\times_{X}U},Y_{U\times_{X}U})}$$
is the natural base change map. Pulling back via $u_{\mathfrak T(U\times_X U)}$ one finds a commutative diagram (\ref{gigadiagram})
\begin{center}
\begin{tikzcd}
u_{\mathfrak T(U\times_X U)}^*q_{?}^*R^if_{(Y_{U},Y_{U})/\mathfrak U, \rig,*}\mathcal O^{\dagger}_{(Y_{U},Y_{U})}\arrow{r}\arrow{d}{\simeq} & u_{\mathfrak T(U\times_X U)}^* R^if_{(Y_{U\times_{X}U},Y_{U\times_{X}U})/\mathfrak U\times_{W} \mathfrak U, \rig,*}\mathcal O^{\dagger}_{(Y_{U\times_{X}U},Y_{U\times_{X}U})}\arrow{d}{\simeq}\\

q_{\mathfrak T(U\times_X U),?}^*u^*R^if_{(Y_{U},Y_{U})/\mathfrak U, \rig,*}\mathcal O^{\dagger}_{(Y_{U},Y_{U})}\arrow{r}\arrow{d}{\simeq} & u_{\mathfrak T(U\times_X U)}^* R^if_{(Y_{U\times_{X}U},Y_{U\times_{X}U})/\mathfrak U\times_{W} \mathfrak U, \rig,*}\mathcal O^{\dagger}_{(Y_{U\times_{X}U},Y_{U\times_{X}U})}\arrow{d}{\simeq}\\

q_{\mathfrak T(U\times_X U),?}^*R^if_{\mathfrak T(U)_1,\mathfrak T(U),\cristalline,*}\mathcal O_{Y_{\mathfrak T(U)_1}/\mathfrak T(U)}\arrow{r} & R^if_{\mathfrak T(U\times_X U)_1,\mathfrak T(U\times_X U),\cristalline,*}\mathcal O_{Y_{\mathfrak T(U\times_X U)_1}/\mathfrak T(U\times_X U)}
\end{tikzcd}
\end{center}
such that the horizontal morphism are the base change morphism. So the morphism
\begin{center}
\begin{tikzcd}
(q_{?}^*g^{*}\mathcal M)_{(\mathfrak T(U\times_X U),z_{\mathfrak T(U\times_X U)})}\arrow{r}\arrow{d}{\simeq} & (g^{',*}\mathcal M)_{(\mathfrak T(U\times_X U),z_{\mathfrak T(U\times_X U)})}\arrow{d}{\simeq}\\
q_{\mathfrak T(U\times_X U),?}^*R^if_{\mathfrak T(U)_1,\mathfrak T(U),\cristalline,*}\mathcal O_{Y_{\mathfrak T(U)_1}/\mathfrak T(U)}\arrow{r} & R^if_{\mathfrak T(U\times_X U)_1,\mathfrak T(U\times_X U),\cristalline,*}\mathcal O_{Y_{\mathfrak T(U\times_X U)_1}/\mathfrak T(U\times_X U)}
\end{tikzcd}
\end{center}
is induced by the base change morphism.
Since the same is true for $R^if_{\Ogus,*}\mathcal O_{Y/K}$ by its construction in Section \ref{Ogushigher}, this shows that the descent diagram \ref{descentdiagram} is commutative, hence $$\psi:g^*\mathcal M\simeq g^*R^if_{\Ogus,*}\mathcal O_{Y/K}$$ gives an isomorphism in the category of descent data for the category $\Fisoc^{(1)}(U)$ of $U$ over $X$. 

By étale and proper descent for convergent isocrystals (\cite[Theorems 4.5 and 4.6]{Ogus}), this implies that $\psi$ descends to an isomorphism 
$$\mathcal M\simeq R^if_{\Ogus,*}\mathcal O_{Y/K} \text{ in } \Fisoc^{(1)}(X)$$
and concludes the proof of Theorem \ref{berthelotconjecture}.
\endproof

\vskip 1\baselineskip
 \textit{eambrosi@unistra.fr}\\
Institut de Recherche Mathématique Avancée, Université de Strasbourg\\

\begin{thebibliography}{99}
\urlstyle{same}
\bibitem[Abe18]{abe}  T. Abe, Langlands correspondence for isocrystals and the existence of crystalline companions for curves, J. Amer. Math. Soc. 31, p. 921-1057, 2018.
\bibitem[AC18]{AbeCaro} T. Abe and D. Caro, Theory of weights in $p$-adic cohomology, American Journal of Mathematics 140, p. 879-975, 2018.
\bibitem[Amb17]{mioUOIp}
E. Ambrosi, A uniform open image theorem in positive characteristic, arXiv:1711.06132 (2017).
\bibitem[Amb19]{miotesi}
E. Ambrosi, $\ell$-adic, $p$-adic and geometric invariants in families of varieties, Ph.D. thesis, École polytechnique, 2019. Available at \url{http://www.theses.fr/2019SACLX019}.
\bibitem[And96]{AndreIHES} Y. André, 
Pour une théorie inconditionnelle des motifs, Publ. Math. IHES 83, p. 5-49, 1996. 
\bibitem[Ber96]{berthelot} 
P. Berthelot, Cohomologie rigide et cohomologie ridige à supports propres, première partie, Prépublication IRMAR 96-03, 1996.
\bibitem[Cad12]{MC} A. Cadoret, 
Motivated cycles under specialization, Groupes de Galois Géométriques et différentiels, Séminaires et Congrès, S.M.F., p. 25-55, 2012.
\bibitem[Cad19]{uniformadelic}
A. Cadoret, An open adelic image theorem for motivic representations over function fields, Math. Res. Lett. 26, p. 1-8, 2019. 
\bibitem[CC20]{brauer} A. Cadoret, F. Charles, A remark on uniform boundedness of Brauer, Algebraic Geometry 7, p. 512-522, 2020
\bibitem[CK16]{cadoretkret}
A. Cadoret, A.Kret, Galois-generic points on Shimura varieties, Algebra and Number Theory 10, p. 1893-1934, 2016.
\bibitem[Chr18]{Chr} A. Christensen, Specialization of Néron-Severi groups in characteristic p, arXiv:1810.06550 (2018).
\bibitem[Cre92]{Crew} R. Crew, $F$-isocrystals and their monodromy groups, Ann. Sci. Ecole Norm. Sup. 25, p. 429-464, 1992.
\bibitem[CLS98]{purityproper} B. Charellotto, B. Le Stum, Sur la pureté de la cohomologie cristalline, C. R. Acad. Sci. Paris 326, p. 961–963, 1998.
\bibitem[dJ]{djtate} A.J. de Jong, Tate conjecture for divisors, unpublished note.
\bibitem[dJ96]{alteration} A.J. de Jong, Smoothness, semi-stability and alterations. Publ. Math. I.H.E.S. 83, p. 51-93, 1996.
\bibitem[Del71]{HodgeII}
P. Deligne, Th\'eorie de Hodge, II, Publ. Math. I.H.E.S. 40, p. 5-57, 1971.
\bibitem[Del80]{Weil2}
 P. Deligne,  La Conjecture de Weil, II, Publ. Math. I.H.E.S. 52, p. 137-252, 1980.
 \bibitem[EElsHKo09]{mordel} J.  Ellenberg,  C. Elsholtz, C. Hall, E.Kowalski,  Non-simple abelian varieties in a family: geometric and analytic approaches, J. London Math. Soc. 80, p. 135-154, 2009.
\bibitem[ES93]{Trace}
Etesse, J.-Y., Le Stum, B.: Fonctions $L$ associées aux $F$-isocristaux surconvergent I. Interprétation cohomologique, Math. Ann. 296, p. 557-576, 1993.
\bibitem[Il79]{derhamwitt} L.Illusie, Complexe de de Rham-Witt et cohomologie cristalline, Ann. Sci. Ecole Norm. Sup. 12, p. 501-661, 1979.
\bibitem[K79]{Katz} N. M. Katz,, Slope filtrations of F-crystals, Astérisque 63, p. 113-164, 1979
\bibitem[KM74]{KatzMessing}  N. M.Katz, W. Messing, Some consequences of the Riemann hypothesis for varieties over finite fields, Invent. Math. 23, p. 73-77, 1974.
\bibitem[Ked04]{Kedfull}
K. S. Kedlaya, Full faithfulness for overconvergent $F$-isocrystals, in (refereed proceedings) Adolphson et al (eds.), Geometric Aspects of Dwork Theory (Volume II), de Gruyter (Berlin), p. 819-83, 2004.
\bibitem[Ked06]{finiteness}
K. S. Kedlaya, Finiteness of rigid cohomology with coefficients, Duke Mathematical Journal 134, p. 15-97, 2006.
\bibitem[Ked17]{notesked}
K. S. Kedlaya, Notes on isocrystals, arXiv:1606.01321v5 (2017).
\bibitem[Laf02]{lafforgue} L. Lafforgue, Chtoucas de Drinfeld et correspondance de Langlands, Invent. Math. 147, p. 1-241 2002.
\bibitem[Laz16]{berthelotconjecture}
C. Lazda, Incarnations of Berthelot's conjecture, J. Number Theory 166, p. 137-157, 2016.
\bibitem[LeS07]{Lestum}
 B. Le Stum, Rigid cohomology, Cambridge Tracts in Mathematics, Cambridge University Press, 2007.
\bibitem[LF15]{Leifu} L.Fu, \'Etale cohomology theory. Revised edition. Nankai Tracts in Mathematics, 14. World Scientific Publishing Co. Pte. Ltd., Hackensack, NJ, 2015.
\bibitem[MP12]{poonen}
D. Maulik, B. Poonen, Néron-Severi groups under specialization, Duke Mathematical Journal 161, p. 2167-2206, 2012.
\bibitem[Mil80]{milne}
J.S. Milne, Étale cohomology, Princeton Mathematical Series, Princeton University Press, 1980.
\bibitem[Mo77]{Mori}
S. Mori, On Tate conjecture concerning endomorphisms of Abelian varieties, International symposium of Algebraic Geometry, Kyoto, p. 219-230, 1977.
\bibitem[Mor19]{Varmor}
M. Morrow,  A Variational Tate Conjecture in crystalline cohomology, Journal of the European Mathematical Society 21, p. 3467-3511, 2019.
\bibitem[Og84]{Ogus}
A. Ogus, $F$-isocrystals and de Rham cohomology II - Convergent isocrystals, Duke Math. J. 51, p. 765-850, 1984.
\bibitem[Ser89]{Serreweil}
J. P. Serre, Lectures on the Mordell-Weil theorem, Aspects of Mathematics, Friedr. Vieweg and Sohn, 1989.
\bibitem[Sh08a]{ShihoI}
A.Shiho, Relative log convergent cohomology and relative rigid cohomology I, arXiv:0707.1742 (2008).
\bibitem[Sh08b]{Shiho}
A.Shiho, Relative log convergent cohomology and relative rigid cohomology II,  arXiv:0707.1743 (2008).
\bibitem[SGA1]{SGA1} Revêtements étales et groupe gondamental, Lecture Notes in Mathematics, Springer-Verlag, Berlin-New York, 1971, Séminaire de géométrie algébrique du Bois Marie 1960–1961, Dirige par Alexandre Grothendieck. Augmente de deux exposes de M. Raynaud. 
\bibitem[SGA2]{SGA2} A.Grothendieck, Cohomologie locale des faisceaux cohérents et théorèmes de Lefschetz locaux et globaux, North-Holland Publishing Co., Amsterdam, Masson and Cie, 1968.
\bibitem[SGA4]{SGA4} Théorie des topos et cohomologie étale des schemas. Tome 3, Seminaire de Geometrie Algébrique du Bois-Marie 1963–1964. Dirige par M. Artin, A. Grothendieck et J. L. Verdier. Avec la collaboration de P. Deligne et B. Saint-Donat. Lecture Notes in Mathematics, Springer, 1973.
\bibitem[SGA7]{SGA7} P. Deligne and N. Katz, Groupe de monodromie en géométrie algébrique (SGA7) - 2éme partie, Lecture Notes in Math., Springer-Verlag, Berlin-New York, 1973.
\bibitem[SP]{stacks-project}
The Stacks Project Authors, Stacks Project, available at \url{http://stacks.math.columbia.edu}.
\bibitem[Tam21]{tamagawa} A. Tamagawa, Uniform open image in positive characteristic, in preparation, 2021.
\bibitem[Ta65]{tateconj} J. Tate, Algebraic cohomology classes and poles of zeta functions, in Schilling, O. F. G., Arithmetical Algebraic Geometry (Proc. Conf. Purdue Univ., 1963), New York: Harper and Row, p. 93-110, 1965.
\bibitem[VAV17]{VAV} A. Varilly-Alvarado, B. Viray, Abelian $n$-division fields of elliptic curves and Brauer groups of product of Kummer and abelian surfaces, Forum of Mathematics Sigma 5, 2017.
\end{thebibliography}
\end{document}